\documentclass[]{report}

\usepackage{amssymb}
\usepackage{amsmath}
\usepackage{amsfonts}
\usepackage{amsthm}

\usepackage{graphicx}
\usepackage{subfig}
\usepackage{enumerate}
\usepackage{verbatim}
\usepackage{color}
\usepackage{float}
\usepackage{tikz}
\usepackage{wrapfig}
\usepackage{listings} 

\usepackage{color}
\usepackage{float}

\usepackage{setspace}
\usepackage[left=2.5cm,top=3cm,right=3cm,bottom=2.5cm]{geometry}

\usepackage[nottoc]{tocbibind}
\newtheorem{thrm}{Theorem}[chapter]
\newtheorem{definition}[thrm]{Definition}
\newtheorem{lemma}[thrm]{Lemma}
\newtheorem{corollary}[thrm]{Corollary}
\newtheorem{claim}[thrm]{Claim}
\newtheorem{conjecture}[thrm]{Conjecture}

\newcommand{\PP}{P}
\newcommand{\NPP}{NP}

\usepackage{pdfpages}

\begin{document}
\begin{titlepage}
\begin{center}
\begin{minipage}{0.1\textwidth}
\includegraphics[height=40pt]{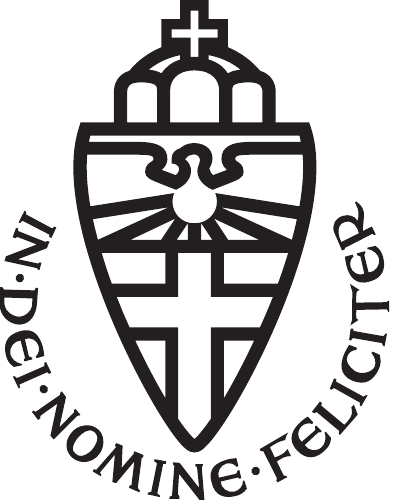}
\end{minipage}
{\Large Radboud Universiteit Nijmegen}\\[1cm]

{ \huge \bfseries Barnette's Conjecture}\\[3.5cm]

\begin{tikzpicture}
    [scale=.2,auto=left,every node/.style={circle,draw, fill=black, scale=.3}]
	\node (O1) at (17,19){};
	\node (O2) at (15,20){};
	\node (O3) at (15,18){};
	\node (O4) at (17,17){};
	\node (O5) at (19,18){};
	\node (O6) at (19,20){};
	\node (O7) at (17,21){};
	\node (O8) at (13,21){};
	\node (O9) at (13,17){};
	\node (O10) at (17,15){};
	\node (O11) at (21,17){};
	\node (O12) at (21,21){};
	\node (O13) at (17,23){};

	\node (P1) at (37,19){};
	\node (P2) at (35,20){};
	\node (P3) at (35,18){};
	\node (P4) at (37,17){};
	\node (P5) at (39,18){};
	\node (P6) at (39,20){};
	\node (P7) at (37,21){};
	\node (P8) at (33,21){};
	\node (P9) at (33,17){};
	\node (P10) at (37,15){};
	\node (P11) at (41,17){};
	\node (P12) at (41,21){};
	\node (P13) at (37,23){};

	\node (Q1) at (17,-5){};
	\node (Q2) at (15,-6){};
	\node (Q3) at (15,-4){};
	\node (Q4) at (17,-3){};
	\node (Q5) at (19,-4){};
	\node (Q6) at (19,-6){};
	\node (Q7) at (17,-7){};
	\node (Q8) at (13,-7){};
	\node (Q9) at (13,-3){};
	\node (Q10) at (17,-1){};
	\node (Q11) at (21,-3){};
	\node (Q12) at (21,-7){};
	\node (Q13) at (17,-9){};

	\node (R1) at (37,-5){};
	\node (R2) at (35,-6){};
	\node (R3) at (35,-4){};
	\node (R4) at (37,-3){};
	\node (R5) at (39,-4){};
	\node (R6) at (39,-6){};
	\node (R7) at (37,-7){};
	\node (R8) at (33,-7){};
	\node (R9) at (33,-3){};
	\node (R10) at (37,-1){};
	\node (R11) at (41,-3){};
	\node (R12) at (41,-7){};
	\node (R13) at (37,-9){};

	\node (A1) at (5,29){};
	\node (A2) at (5,21){};
	\node (A3) at (7,23){};
	\node (A4) at (9,25){};
	\node (A5) at (11,27){};
	\node (A6) at (13,29){};
	\node (A7) at (7,27){};

	\node (B1) at (49,29){};
	\node (B2) at (49,21){};
	\node (B3) at (47,23){};
	\node (B4) at (45,25){};
	\node (B5) at (43,27){};
	\node (B6) at (41,29){};
	\node (B7) at (47,27){};

	\node (C1) at (5,-15){};
	\node (C2) at (5,-7){};
	\node (C3) at (7,-9){};
	\node (C4) at (9,-11){};
	\node (C5) at (11,-13){};
	\node (C6) at (13,-15){};
	\node (C7) at (7,-13){};

	\node (D1) at (49,-15){};
	\node (D2) at (49,-7){};
	\node (D3) at (47,-9){};
	\node (D4) at (45,-11){};
	\node (D5) at (43,-13){};
	\node (D6) at (41,-15){};
	\node (D7) at (47,-13){};

  \foreach \from/\to in {
	A1/A7, A4/A5, A4/A3, A3/A2, A5/A6,
	B1/B7,  B4/B5, B4/B3, B3/B2, B5/B6,
	C1/C7, C4/C5, C4/C3, C3/C2, C5/C6,
	D1/D7, D4/D5, D4/D3, D3/D2, D5/D6,
	O1/O2, O1/O6, O1/O4, O2/O7, O2/O3, O7/O6, O6/O5, O5/O4, O4/O3, O7/O13, O11/O5, O3/O9, O8/O9, O9/O10, 	
	O10/O11, O11/O12, O12/O13, O13/O8,
	O12/P8,
	P1/P2, P1/P6, P1/P4, P2/P7, P2/P3, P7/P6, P6/P5, P5/P4, P4/P3, P7/P13, P11/P5, P3/P9, P8/P9, P9/P10, 	
	P10/P11, P11/P12, P12/P13, P13/P8,
Q1/Q2, Q1/Q6, Q1/Q4, Q2/Q7, Q2/Q3, Q7/Q6, Q6/Q5, Q5/Q4, Q4/Q3, Q7/Q13, Q11/Q5, Q3/Q9, Q8/Q9, Q9/Q10, 	
	Q10/Q11, Q11/Q12, Q12/Q13, Q13/Q8,
	Q12/R8,
	R1/R2, R1/R6, R1/R4, R2/R7, R2/R3, R7/R6, R6/R5, R5/R4, R4/R3, R7/R13, R11/R5, R3/R9, R8/R9, R9/R10, 	
	R10/R11, R11/R12, R12/R13, R13/R8,
	O10/Q10, P10/R10}
    	\draw (\from) -- (\to);

  \path[every node/.style={}]

	(A6) edge [very thick, purple] node [black] {} (A1)
	(A1) edge [very thick, purple] node [black] {} (A2)
	(A5) edge [very thick, purple] node [black] {} (A7)
	(A7) edge [very thick, purple] node [black] {} (A3)
	(B6) edge [very thick, purple] node [black] {} (B1)
	(B1) edge [very thick, purple] node [black] {} (B2)
	(B5) edge [very thick, purple] node [black] {} (B7)
	(B7) edge [very thick, purple] node [black] {} (B3)
	(A6) edge [] node [black] {} (B6)

	(C6) edge [very thick, purple] node [black] {} (C1)
	(C1) edge [very thick, purple] node [black] {} (C2)
	(C5) edge [very thick, purple] node [black] {} (C7)
	(C7) edge [very thick, purple] node [black] {} (C3)
	(D6) edge [very thick, purple] node [black] {} (D1)
	(D1) edge [very thick, purple] node [black] {} (D2)
	(D5) edge [very thick, purple] node [black] {} (D7)
	(D7) edge [very thick, purple] node [black] {} (D3)
	(C6) edge [] node [black] {} (D6)

	(A2) edge [] node [black] {} (C2)
	(B2) edge [] node [black] {} (D2);

  \foreach \from/\to in {A3/A4, A4/O8, O8/O13, O13/O7, O7/O2, O2/O1, O1/O6, O6/O5, O5/O4, O4/O3, O3/O9, O9/O10, O10/O11, O11/O12, O12/P8, P8/P9, P9/P10, P10/P11, P11/P5, P5/P4, P4/P3, P3/P2, P2/P1, P1/P6, P6/P7, P7/P13, P13/P12, P12/B4, B4/B3, A3/A4, B3/B4, A5/A6, B5/B6,
C4/Q8, Q8/Q13, Q13/Q7, Q7/Q2, Q2/Q1, Q1/Q6, Q6/Q5, Q5/Q4, Q4/Q3, Q3/Q9, Q9/Q10, Q10/Q11, Q11/Q12, Q12/R8, R8/R9, R9/R10, R10/R11, R11/R5, R5/R4, R4/R3, R3/R2, R2/R1, R1/R6, R6/R7, R7/R13, R13/R12, R12/D4, C3/C4, D3/D4, C5/C6, D5/D6, A2/C2, B2/D2}
    	\draw [very thick, purple](\from) -- (\to);

\end{tikzpicture}
\\[3cm]
\begin{minipage}{0.4\textwidth}
\begin{flushleft} \large
\emph{Authors:}\\  
Lean Arts\\
Meike Hopman\\
Veerle Timmermans
\end{flushleft}
\end{minipage}
\begin{minipage}{0.4\textwidth}
\begin{flushright} \large
\emph{Supervisor:} \\
Wieb Bosma\\
\vspace{1cm}
\end{flushright}
\end{minipage}
\vfill
{\large \today}
\end{center}
\end{titlepage}

\begin{abstract}
This report provides an overview of theorems and statements related to a conjecture stated by D. W. Barnette in 1969. Barnette’s conjecture is an open problem in graph theory about Hamiltonicity of graphs.

\begin{conjecture}[Barnette’s conjecture]Every cubic, bipartite, polyhedral graph contains a Hamiltonian cycle. \end{conjecture}

We will have a look Steinitz’s theorem, which allows to reformulate the conjecture. We look at the complexity of previous, related, conjectures and determine some necessary conditions for a Barnette graph to contain a Hamiltonian cycle. We will also determine some sufficient conditions for a graph to contain a Hamiltonian cycle.
We will give a sketch of the proof for Barnette’s conjecture for graphs up to 64 vertices and we will describe two operations by which all Barnette graphs can be constructed. Finally, we will describe some programs in C++ that will help to check some properties in graphs, like bipartiteness, 3-connectedness and Hamiltonicity.
\end{abstract}

\tableofcontents

\chapter{Introduction}
\label{introduction}
This report provides an overview of theorems and statements related to Barnette's conjecture stated by D. W. Barnette in 1969. Barnette's conjecture is an open problem in graph theory about Hamiltonicity of graphs. 
\begin{conjecture}[Barnette's conjecture]
Every cubic $3$-connected bipartite planar graph contains a Hamiltonian cycle.
\end{conjecture}

Before we can discuss related theorems and statements we first have to understand what type of graphs Barnette's conjecture is about. Chapter \ref{Steinitz} is on Steinitz's theorem, which tells us something about the correspondence between polyhedral graphs and $3$-dimensional polytopes. Here we will conclude that a graph is polyhedral if and only if it is 3-connected and planar. We will only give a sketch of the proof of Steinitz' theorem, but this will be enough to understand the relation between graphs considered in Barnette's conjecture and $3$-dimensional polytopes.

Before looking at Barnette's conjecture we consider two weaker conjectures that are closely related to Barnette's conjecture. 
In  1880 P.G. Tait stated the following open problem:
\begin{conjecture}[Tait's conjecture]
Every cubic $3$-connected planar graph contains a Hamiltonian cycle.
\end{conjecture}
Using Steinitz' theorem of Chapter \ref{Steinitz}, it is similar to Barnette's conjecture, except the graph does not have to be bipartite. In 1946 W.T. Tutte was the first one who found a counterexample on 46 vertices, known as Tutte's graph \cite{Tutte1946}. After Tait's conjecture was proven to be false, people kept looking for smaller counterexamples. In 1986 it was proven that the smallest counterexample consists of 38 vertices \cite{Holton1988}. In Chapter \ref{Tait} on Tait's conjecture we consider Tutte's graph, the smallest counterexample and its relation to the four color theorem.

Tutte noticed that none of the counterexamples that were found are bipartite, so he adapted the conjecture:
\begin{conjecture}[Tutte's conjecture]
Every cubic 3-connected bipartite graph contains a Hamiltonian cycle.
\end{conjecture}
Again, this is similar to Barnette's conjecture, except that the graph does not have to be planar. In 1976 J.D. Horton found the first counterexample, on 96 vertices, now called the Horton graph \cite{BM1976}. The smallest currently known counterexample is Georges's graph on 50 vertices \cite{Georges1989}. In Chapter \ref{Tutte}, on Tutte's conjecture, we consider the Horton fragment and the Ellingham fragment, which both play a role in the counterexamples.

In Chapter \ref{complexity} we will discuss the complexity of some problems related to Barnette's conjecture. We will prove that the problems corresponding to Tait's and Tutte's conjecture are \NPP-complete. This leads us to an idea of how we could prove something about the complexity of Barnette's conjecture. If we are able to find a certain type of subgraph that is cubic, $3$-connected, planar and bipartite, then we are able to prove that Barnette's conjecture is also \NPP-complete. In Tait's and Tutte's conjecture we have seen what happens if we remove one of the properties of Barnette graphs. In Chapter \ref{complexity} on complexity we have proven that these problems become \NPP-complete. In Chapter \ref{4connected} we change the statement into: {\it Every $4$-connected planar graph contains a Hamiltonian cucle.} We will prove that this statement, in which we have replaced $3$-connected by the stronger property $4$-connected,  is true. If we combine these two facts we conclude that Barnette's conjecture is on the boundary between \NPP-completeness and \PP. In these two chapters we thus prove that bipartiteness and planarity are necessary conditions in Barnette's conjecture and that 4-connectedness is sufficient for a Barnette graph to contain a Hamiltonian cycle.

Although Barnette's conjecture remains unsolved, Holton, Manvel and McKay proved in 1984 \cite{HMM1985} that Barnette's conjecture is true up to 64 vertices. In Chapter \ref{Proofonvertices} we will have a look at this proof and try to sketch their main idea. This proof consists of two parts: they reduce Barnette graphs to graphs with fewer vertices and use the computer to check all 'smaller' graphs for certain properties.
To check those smaller graphs, Holton, Manvel and KcKay developed a way to generate all Barnette graphs up to isomorphism in a fairly simple way. In Chapter \ref{Thetwooperations} we will have a look at how all Barnette graphs can be generated and prove why this indeed works. This will give us some new insights in the difficulty of proving the conjecture.
The figure on the title page is a Barnette graph. It is generated from the graph of a cube by using the two operations that are treated in Chapter~\ref{Thetwooperations}. One Hamiltonian cycle is indicated in purple.

Apart from finding a way to generate all Barnette graphs, we also need the computer to check certain properties. For example: we wanted to verify the counterexamples given for Tait's and Tutte's conjecture. We wrote some programs in C++ to check if they are indeed bipartite, cubic and 3-connected. We also wrote some programs to check for Hamiltonian cycles and Hamiltonian paths. In Chapter \ref{Programsinc} we will give and explain the programs written in C++. 

At the end of this project we thought about possible ways to construct a counterexample to Barnette's conjecture. Chapter \ref{Counterexample} shows our way of thinking. 


\section{Preliminaries}
\begin{definition}
A {\bf graph} $G$  is an ordered pair $(V(G), E(G))$ of a vertex set $V(G)$ and an edge set $E(G)$. An edge $e\in E(G)$ is of the form $e=\{v_1, v_2\}$ with $v_1, v_2 \in V(G)$. We will use $v_1v_2$ or $v_2v_1$ to denote edge $\{v_1, v_2\}$. We say that $v_1$ and $v_2$ are the endpoints of $e$.
\end{definition}

\begin{definition}
A {\bf simple graph} $G$ is a graph in which there is no more than one edge between any two different vertices and in which every edge has two different endpoints.
\end{definition}

\begin{definition}
A {\bf subgraph} $H$ of a graph $G$ is a graph for which $V(H) \subseteq V(G)$, $E(H) \subseteq E(G)$ and $v_1, v_2 \in V(H)$ for every $v_1v_2 \in E(H)$. 
 \end{definition}

\begin{definition}
The {\bf degree} of a vertex $v \in V(G)$ is the number of edges of which $v$ is an endpoint. For every edge of the form $vv$, one extra is added to the degree of $v$. \emph{Notation: $\deg(v)$.}
\end{definition}

\begin{definition}
A {\bf path} of length  $m$ in a graph is a sequence of vertices $(v_0, v_1, \dots, v_m)$ such that $v_iv_{i+1} \in E(G) \forall i\in\{0, 1, \dots, m-1\}$. 
\end{definition}

\begin{definition}
A {\bf cycle} of length $m$ in a graph is a path $(v_0, v_1, \dots, v_m)$ where $v_0=v_m$ and $v_i\not=v_j$ for each other pair $i,j\in\{0, 1, \dots, m\}$.
\end{definition}

\begin{definition}
A {\bf Hamiltonian path} is a path that contains each vertex of the graph exactly once. 
\end{definition}

\begin{definition}
A {\bf Hamiltonian cycle} is a cycle that contains each vertex of the graph exactly once.
\end{definition}

\begin{definition}
A {\bf cubic} graph is a graph in which all vertices have degree three.
\end{definition}

\begin{definition}
A graph is {\bf connected} if there is a path from each $v_1 \in V(G)$ to all other $v_2 \in V(G)$. 
\end{definition}

\begin{definition}
A graph is {\bf $k$-connected} if it has more than $k$ vertices and the result of deleting any set of fewer than $k$ vertices is a connected graph.
\end{definition}

\begin{definition}
 In a {\bf complete graph} on $t$ vertices there exists an edge between every pair of vertices. \emph{Notation: $K_t$.}
 \end{definition}

\begin{definition}
A {\bf bipartite} graph is a graph $G$ for which $V(G)=V_1 \cup V_2$, with $V_1 \cap V_2 = \emptyset$, such that every edge has one endpoint in $V_1$  and one in $V_2$. 
\end{definition}

\begin{definition}
A {\bf complete bipartite graph} is a bipartite graph $G$ with $E(G)=\{v_1v_2 | \, v_1\in V_1 \, \mbox{and} \ v_2\in V_2\}$. \emph{Notation: $K_{|V_1|, |V_2|}$.}
 \end{definition}

\begin{definition}
A graph is {\bf planar} if it can be drawn on the plane in such a way that its edges intersect only at their endpoints.
\end{definition}

\chapter{Steinitz's Theorem}
\label{Steinitz}
Barnette's conjecture is a statement on graphs that are bipartite, cubic and polyhedral. The properties bipartite and cubic are well known. To understand what it means that a graph is polyhedral we will study Steinitz's theorem about the correspondence between polyhedral graphs and $3$-dimensional convex polytopes. A polyhedral graph is a simple, $3$-connected, planar graph.
\begin{thrm}[Steinitz's theorem]
$G$ is the graph of a $3$-dimensional convex polytope if and only if $G$ is simple, planar and $3$-connected.
\end{thrm}
In this chapter we will only give a sketch of the proof to get an idea of what polyhedral graphs look like. The entire proof can be found in \cite{Ziegler} in Chapter 4. In this chapter graphs are not simple, unless stated otherwise. So graphs can have parallel edges and loops. A loop is an edge of the form $xx$, with $x$ a vertex of the graph.

Suppose $G$ is the graph of a $3$-dimensional convex polytope. It is clear that $G$ is a simple, planar graph. Intuitively this follows from radial projection of the graph onto a sphere from a point inside the polytope. The projection onto the sphere can easily be transformed into a planar graph by opening the sphere at a face and make this face the outer face of the planar representation of the graph by flattening the sphere. Balinski's theorem (Theorem 3.14 in Chapter 3 of \cite{Ziegler}) states that the graph of a $d$-dimensional polytope is $d$-connected. This completes the first part of the proof. To show that every simple, planar and $3$-connected graph is the graph of a $3$-dimensional convex polytope is a bit more complicated. We use two reductions and two transformations for this.

\noindent \begin{minipage}{0.5\textwidth}
\begin{definition}
A {\bf serial reduction} consists of removing a subdivision point. A subdivision point is a point $x\in V(G)$ that is added to a graph $G'$, with $yz \in E(G')$, such that we obtain graph $G$ with $V(G)=V(G')\cup \{x\}$ and $E(G)=(E(G') - \{yz\}) \cup \{xy, xz\}$.
\end{definition}
\end{minipage}
\begin{minipage}{0.5\textwidth}
\begin{figure}[H]
\centering
  \includegraphics[scale= 0.5]{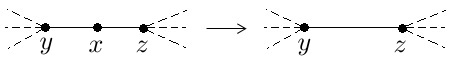} 
\caption{Serial reduction}
\label{onderverdeling}
\end{figure}
\end{minipage}

\noindent \begin{minipage}{0.5\textwidth}
\begin{definition}
A {\bf parallel reduction} consists of removing a parallel edge. 
\end{definition}
\end{minipage}
\begin{minipage}{0.5\textwidth}
\begin{figure}[H]
\centering
  \includegraphics[scale= 0.5]{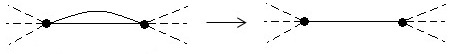} 
\caption{Parallel reduction}
\label{onderverdeling}
\end{figure}
\end{minipage}

In the remainder of this chapter we will use the terminology {\it SP reductions} to refer to a sequence of serial and/or parallel reductions.

\begin{minipage}{0.6\textwidth}
\begin{definition}
In a {\bf $\Delta$-to-$Y$ transformation} a triangle $T=(\{x, y, z\}, \{xy, xz, yz\})$, which can be a subgraph, is replaced by a $3$-star $S=(\{w, x, y, z\}, \{wx, wy, wz\})$, with the restriction that $w$ is a vertex of degree $3$.\\
In a {\bf $Y$-to-$\Delta$ transformation} a $3$-star $S$, which can be a subgraph, is replaced by a triangle $T$.\\
\end{definition}
\end{minipage}
\begin{minipage}{0.4\textwidth}
\begin{figure}[H]
\centering
  \includegraphics[scale= 0.5]{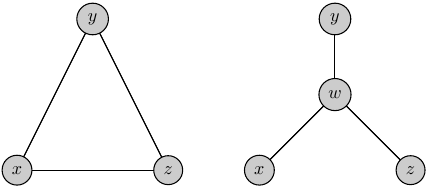} 
\caption{Triangle ($\Delta$) and $3$-star ($Y$)}
\label{transformation}
\end{figure}
%
%
 %
%
%
 %
\end{minipage}

 Note that between vertices $x$, $y$ and $z$ more edges may exist; those are preserved after the transformation.

The dual graph $G'=(V(G'), E(G'))$ of a planar graph $G$ is defined as follows: $V(G')=\{v\ | \ v \mbox{ is a face of} \ G\}$ and $E(G')=\{vw| v, w \in V(G'), \ v \ \mbox{and} \ w \ \mbox{share an edge in} \ G\}$. This construction is shown in Figure \ref{dual}, where the red graph is the dual of the blue graph and vice versa. 

\noindent \begin{minipage}{0.6\textwidth}
Clearly the dual of a graph depends on the planar embedding of the graph. Note that the dual graph of a planar graph is still a planar graph. This holds because in every face of the planar graph there is exactly one vertex of the dual graph and every edge of this face is crossed by exactly one edge in the dual graph. So edges in the dual graph can be drawn without crossing, thus the dual graph is planar. Note that if $H$ is the dual graph of $G$, then $G$ is the dual graph of $H$, so the dual graph of the dual graph of $G$ is $G$.
\end{minipage}
\begin{minipage}{0.4\textwidth}
\begin{figure}[H]
\centering
  \includegraphics[scale= 0.34]{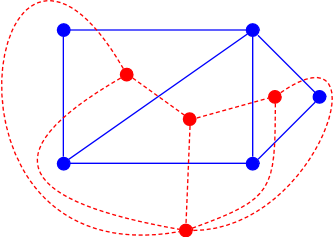} 
\caption{Dual graph}
\label{dual}
\end{figure}
\end{minipage}

The following lemma states that $3$-connectedness is preserved under $Y$-to-$\Delta$ transformations after deletion of parallel edges. For $\Delta$-to-$Y$ transformations the same holds. This can be proved by looking at the dual graph, in which $\Delta$-to-$Y$ transformations correspond to $Y$-to-$\Delta$ transformations and serial reductions correspond to parallel reductions.

\begin{lemma}
\label{connected1}
Given are a graph $G$ that is $3$-connected, but not $K_4$, and a vertex $w$ of $G$ of degree three. We can do a $Y$-to-$\Delta$ transformation on $w$ and its neighbours. If we delete all parallel edges after a $Y$-to-$\Delta$ transformation, then the result is still $3$-connected. If $v$ has three different neighbors, the result of a $Y$-to-$\Delta$ transformation is still $3$-connected.
\end{lemma}

\begin{proof} Consider a $3$-connected graph $G$ (not $K_4$) and let $G'$ be the result of a $Y$-to-$\Delta$ transformation after deletion of parallel edges. Suppose $G'$ is no longer $3$-connected. Consider the separating set $S$ of $G^{'}$ of one or two vertices. It is clear that this is a separating set for $G$ if $w$ has three different neighbours. So by contradiction $3$-connectedness is preserved. \end{proof}

We will use the terminology {\it simple $\Delta Y$ reduction} to refer to a $Y$-to-$\Delta$ or $\Delta$-to-$Y$ transformation directly followed by all possible SP reductions. Lemma \ref{connected1} and its dual version state that simple $\Delta Y$ reductions preserve $3$-connectedness.

\begin{lemma}
\label{polytope}
Let $G$ be a $3$-connected planar graph and let $G'$ be the result of a simple $\Delta Y$ reduction. If $G'$ is the graph of a $3$-dimensional convex polytope, then so is $G$.
\end{lemma}

This lemma is intuitively true by comparing $Y$-to-$\Delta$ transformations to cutting off a vertex of a $3$-dimensional convex polytope and $\Delta$-to-$Y$ transformations to the opposite operation. A formal proof can be found in \cite{Ziegler} as the proof of Lemma 4.3 in Chapter 4 on page 108.

\noindent \begin{minipage}{0.5\textwidth}
\begin{definition} Graph $H$ is a {\bf minor} of graph $G$ if $H$ can be obtained from $G$ by deleting edges and vertices and by contracting edges.
\end{definition}
\end{minipage}
\begin{minipage}{0.5\textwidth}
\begin{figure}[H]
\centering
  \includegraphics[scale= 0.5]{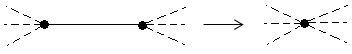} 
\caption{Contraction of an edge}
\label{contraction}
\end{figure}
\end{minipage}

The following three lemma's are used in the proof of Steinitz's theorem. They are stated without a proof, which can be found in \cite{Ziegler} as the proof of Lemma 4.4, Lemma 4.5 and Lemma 4.6 in Chapter 4 on pages 110, 111 and 112.

\begin{lemma}
\label{minor}
If a planar $G$ can be reduced by simple $\Delta Y$ reductions to two vertices with two parallel edges, then so can every $2$-connected minor $H$ of $G$.
\end{lemma}

\begin{lemma}
\label{grid}
If $G$ is planar, then it is a minor of a grid graph.
\end{lemma}

\noindent \begin{minipage}{0.6\textwidth}
A grid graph $G(m,n)$ is a planar graph on $mn$ vertices. It has a vertex corresponding to every pair of integers $(a,b)$ for $a \in \{1, \dots, n\}$ and $b\in \{1, \dots, m\}$ and an edge connecting $(a,b)$ to $(a+1, b)$ for $a \in \{1, \dots, n-1\}$ and $b\in \{1, \dots, m\}$ and an edge connecting $(a,b)$ to $(a, b+1)$ for $a \in \{1, \dots, n\}$ and $b\in \{1, \dots, m-1\}$.

The degree of any vertex in the graphs we are considering in Barnette's conjecture is $3$. Now we can place every vertex of graph $G$ on a coordinate of a grid. From this representation of $G$ it is clear that it is a minor of a grid graph, because we can subdivide edges to obtain some of the missing vertices of the grid graph until we have a subgraph of a grid graph.
\end{minipage}
\begin{minipage}{0.4\textwidth}
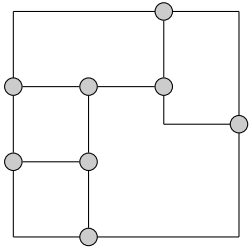
\begin{figure}[H]
\centering
\begin{tikzpicture}
  [scale=.5,auto=left,every node/.style={circle, draw,fill=black!20, scale=0.7}]

\draw  plot coordinates {   (0,6)   (4,6)  };
\draw  plot coordinates {   (0,6)   (0,4)  };
\draw  plot coordinates {   (6,6)   (4,6)  };
\draw  plot coordinates {   (6,6)  (6,3)   };
\draw  plot coordinates {   (4,3)  (4,4)   };
\draw  plot coordinates {   (4,3)  (6,3)   };
\draw  plot coordinates {   (6,0)  (2,0)   };
\draw  plot coordinates {   (6,0)  (6,3)   };
\draw  plot coordinates {   (0,0)  (0,2)   };
\draw  plot coordinates {   (0,0)  (2,0)   };

  	\node (B) at (0,2) {};
  	\node (C) at (0,4) {};

  	\node (E) at (2,0) {};
  	\node (F) at (2,2) {};
  	\node (G) at (2,4) {};

  	\node (I) at (4,4) {};
  	\node (J) at (4,6) {};

  	\node (L) at (6,3) {};

  \foreach \from/\to in {B/C, B/F, C/G, E/F, F/G, G/I, I/J}
    	\draw (\from) -- (\to);
\end{tikzpicture}
\caption{Minor of a grid graph}
\end{figure}
\end{minipage}

\begin{lemma}
\label{reducible}
All grid graphs $G(m, n)$ with $m, n \geq 3$ are reducible to $K_4$ using simple $\Delta Y$ reductions.
\end{lemma}

This can be checked easily using $\Delta$-to-$Y$ transformations and SP reductions. Now we are ready to prove the second part of Steinitz's theorem.

\begin{proof}[Proof of the second part of Steinitz's theorem]
Let graph $G$ be simple, planar and $3$-connected. By Lemma \ref{grid} graph $G$ is a minor of a grid graph $G'$. By Lemma \ref{reducible} graph $G'$ is reducible to $K_4$ using simple $\Delta Y$ reductions, so it is also reducible to two vertices with two parallel edges. This can be shown by first applying a $Y$-to-$\Delta$ transformation to $K_4$, followed by parallel reductions on the three sides of the formed triangle. Apply a serial reduction on two edges of this triangle and the result will be two vertices with two parallel edges.

The $3$-connected graph $G$ is a minor of the graph $G'$. So by Lemma \ref{minor} graph $G$ can also be reduced to two vertices with two parallel edges using simple $\Delta Y$ reductions. By Lemma \ref{connected1} and its dual version we know that $3$-connectedness is preserved. We follow this reduction until serial or parallel edges have to be reduced. This results in a graph that has fewer edges. By induction on the number of edges we now know that every $3$-connected planar graph $G$ can be reduced to $K_4$ by a sequence of simple $\Delta Y$ transformations. We can complete the proof of Steinitz's theorem using Lemma \ref{polytope}. Graph $G$ is the graph of a $3$-dimensional convex polytope, because $K_4$ is the graph of a tetrahedron, which is a $3$-dimensional convex polytope.
\end{proof}

\chapter{Tait's Conjecture}
\label{Tait}
In 1880 P.G. Tait formulated a statement that is related to Barnette's conjecture and which is known as Tait's conjecture:

\begin{conjecture}[Tait's conjecture]
Every cubic 3-connected planar graph contains a Hamiltonian cycle.
\end{conjecture}

In 1946 it was disproven by W.T. Tutte, who constructed a couterexample on 46 vertices \cite{Tutte1946}. Later smaller counterexamples were given, in many cases based on Grinberg's theorem \cite[page 479-481]{BM2008}.  In 1965 Lederberg found the Barnette-Bos\'ak-Lederberg graph on 38 vertices, in 1988 Holton and McKay proved that no smaller counterexample exists \cite{Holton1988}.

The existence of these counterexamples shows that bipartiteness in Barnette's conjecture is a necessary condition. 

\section{Tutte's counterexample}

The counterexample Tutte gave to Tait's conjecture is known as Tutte's graph, see Figure~\ref{Tutte's graph}. It first  appeared in \cite{Tutte1946}. The graph consists of three Tutte fragments; one such fragment is shown in Figure~\ref{Tutte's fragment}.

\begin{definition}
\label{requirededgefragment}
A {\bf required edge fragment} is a fragment of a graph with three half edges with the properties that follow. A half edge is an edge with only one endpoint. These three half edges are connected to vertices of the graph that are not contained in the fragment. The {\bf required edge} is one of these half edges which must be used in any Hamiltonian cycle of the graph which contains the fragment. A Hamiltonian path goes through the fragment by using the required edge and one of the other half edges and it visits every vertex in the fragment.
\end{definition}
\begin{minipage}{.5\textwidth}
\begin{figure}[H]
\centering
\begin{tikzpicture}
  [scale=.145,auto=left,every node/.style={circle,draw, fill=black!100, scale=0.5}]
  	\node (A1) at (16,4) {};
  	\node (B1) at (18,6) {};
  	\node (C1) at (16,8) {};
  	\node (D1) at (14,7) {};
  	\node (E1) at (14,5) {};
  	\node (G1) at (21,7) {};
  	\node (I1) at (12,8) {};
  	\node (J1) at (11,5) {};
  	\node (K1) at (21,2) {};
  	\node (L1) at (16,1) {};
  	\node (M1) at (11,2) {};
  	\node (N1) at (19,11) {};
  	\node (O1) at (16,10) {};
  	\node (P1) at (13,11) {};
  	\node (Q1) at (16,13) {};

  	\node (R) at (16,16) {};

  \foreach \from/\to in {A1/B1, B1/C1, C1/D1, D1/E1, E1/A1, A1/L1, B1/G1, C1/O1, D1/I1, E1/J1, O1/P1, O1/N1, M1/J1, J1/I1, I1/P1, P1/Q1, Q1/R, Q1/N1, N1/G1, G1/K1, K1/L1, L1/M1}
    	\draw (\from) -- (\to);

  	\node (A2) at (26.3,22) {};
  	\node (B2) at (23.6,22.8) {};
  	\node (C2) at (23,20) {};
  	\node (D2) at (24.8,18.8) {};
  	\node (E2) at (26.6,19.8) {};
  	\node (G2) at (21.3,25) {};
  	\node (I2) at (24.8,16.5) {};
  	\node (J2) at (28,17.2) {};
  	\node (K2) at (25.5,27.3) {};
  	\node (L2) at (29,23.6) {};
  	\node (M2) at (30.8,18.8) {};
  	\node (N2) at (19,21) {};
  	\node (O2) at (21,19) {};
  	\node (P2) at (21.8,16) {};
  	\node (Q2) at (18.8,17.6) {};

  \foreach \from/\to in {A2/B2, B2/C2, C2/D2, D2/E2, E2/A2, A2/L2, B2/G2, C2/O2, D2/I2, E2/J2, O2/P2, O2/N2, M2/J2, J2/I2, I2/P2, P2/Q2, Q2/R, Q2/N2, N2/G2, G2/K2, K2/L2, L2/M2}
    	\draw (\from) -- (\to);

	\node (A3) at (5.6,22) {};
  	\node (B3) at (6.3,19.2) {};
  	\node (C3) at (9.2,20) {};
  	\node (D3) at (9.2,22.2) {};
  	\node (E3) at (7.6,23.2) {};
  	\node (G3) at (5.8,16.2) {};
  	\node (I3) at (11.2,23.4) {};
  	\node (J3) at (9,25.8) {};
  	\node (K3) at (1.2,18.4) {};
  	\node (L3) at (3,23.4) {};
  	\node (M3) at (6.4,27.2) {};
  	\node (N3) at (10.2,16) {};
  	\node (O3) at (10.8,19) {};
  	\node (P3) at (13.2,21) {};
  	\node (Q3) at (13.4,17.5) {};

 \foreach \from/\to in {A3/B3, B3/C3, C3/D3, D3/E3, E3/A3, A3/L3, B3/G3, C3/O3, D3/I3, E3/J3, O3/P3, O3/N3, M3/J3, J3/I3, I3/P3, P3/Q3, Q3/R, Q3/N3, N3/G3, G3/K3, K3/L3, L3/M3}
	\draw (\from) -- (\to);

	\draw (K1) to [bend right=20] (M2); 
  	\draw (K2) to [bend right=20] (M3); 
	\draw (K3) to [bend right=20] (M1); 
\end{tikzpicture}
\caption{Tutte's graph}
\label{Tutte's graph}

\end{figure}
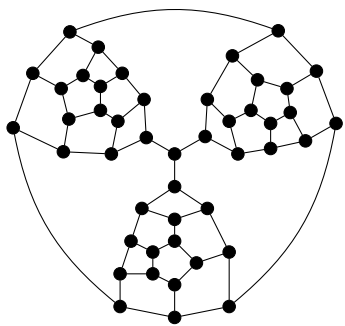
\end{minipage}
\begin{minipage}{.5\textwidth}
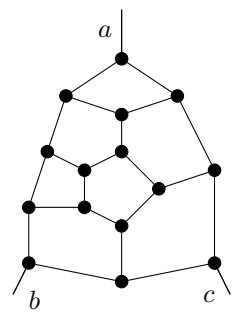
\begin{figure}[H]
\centering
\begin{tikzpicture}
  [scale=.247,auto=left,every node/.style={circle,draw, fill=black!100, scale=0.5}]
  	\node (A) at (6,4) {};
  	\node (B) at (8,6) {};
  	\node (C) at (6,8) {};
  	\node (D) at (4,7) {};
  	\node (E) at (4,5) {};
  	\node (G) at (11,7) {};
  	\node (I) at (2,8) {};
  	\node (J) at (1,5) {};
  	\node (K) at (11,2) {};
  	\node (L) at (6,1) {};
  	\node (M) at (1,2) {};
  	\node (N) at (9,11) {};
  	\node (O) at (6,10) {};
  	\node (P) at (3,11) {};
  	\node (Q) at (6,13) {};
  	\node[white](X) at (6,16) {};
	\node [white](Y) at (12,0) {};
	\node [white](Z) at (0,0) {};

  \foreach \from/\to in {A/B, B/C, C/D, D/E, E/A, A/L, B/G, C/O, D/I, E/J, O/P, O/N, M/J, J/I, I/P, P/Q, Q/N, N/G, G/K, K/L, L/M, Q/X, M/Z, K/Y}
    	\draw (\from) -- (\to);

  \path[every node/.style={font=\sffamily}]
	(Q) edge [] node [black] {$a$} (X)
	(M) edge [] node [black] {$b$} (Z)
	(Y) edge [] node [black] {$c$} (K);

\end{tikzpicture}
\caption{Tutte's fragment}
\label{Tutte's fragment}
\end{figure}
\end{minipage}

\begin{lemma}
Tutte's fragments is a required egde fragment where edge $a$ is the required edge.
\end{lemma}
\begin{proof}
Consider a pentagonal prism, given in Figure~\ref{Graph 1}.

It is easy to see that a Hamiltionian cycle in Figure~\ref{Graph 1} can not contain both $AF$ and $CH$. This is similar to saying that no Hamiltonian cycle in Figure~\ref{Graph 2} can contain both $AL$ and $CO$. 

When a Hamiltonian cycle in Figure~\ref{Graph 2} contains $KM$ then it contains either $KL$ or $LM$ and therefore it contains $AL$. For the same reason, if it contains $NP$ then it also contains $OC$. So no Hamiltonian cycle in Figure~\ref{Graph 2} can contain both $KM$ and $NP$. 

Now consider the graph in Figure~\ref{Graph 3}. If there is a Hamiltonian cycle in Figure~\ref{Graph 3}, it must contain $QR$, because otherwise we obtain a Hamiltonian cycle in Figure~\ref{Graph 2} which contains both $KM$ and $NP$.

\begin{figure}[h]
\centering
\begin{minipage}{.32\textwidth}

\begin{center}
 \begin{tikzpicture}
  [scale=.36,auto=left,every node/.style={circle,draw, fill=black!20,scale=0.5}]
  	\node (A) at (5,7) {$A$};
  	\node (B) at (8,9) {$B$};
  	\node (C) at (11,7) {$C$};
  	\node (D) at (10,4) {$D$};
  	\node (E) at (6,4) {$E$};
  	\node (G) at (8,12) {$G$};
  	\node (I) at (12,1) {$I$};
  	\node (J) at (4,1) {$J$};
  	\node (H) at (14,9) {$H$};
  	\node (F) at (2,9) {$F$};

  \foreach \from/\to in {A/B, B/C, C/D, D/E, E/A, A/F, B/G, C/H, D/I, E/J, J/I, I/H, H/G, G/F, F/J}
    	\draw (\from) -- (\to);
 
\end{tikzpicture}
\caption{}
\label{Graph 1}\
\end{center}
\end{minipage}
\begin{minipage}{.32\textwidth}

\begin{center}
\begin{tikzpicture}
  [scale=.33,auto=left,every node/.style={circle, draw ,fill=black!20, scale=0.5}]
  	\node (A) at (5,7) {$A$};
  	\node (B) at (8,9) {$B$};
  	\node (C) at (11,7) {$C$};
  	\node (D) at (10,4) {$D$};
  	\node (E) at (6,4) {$E$};
  	\node (G) at (8,12) {$G$};
  	\node (I) at (12,1) {$I$};
  	\node (J) at (4,1) {$J$};
  	\node (K) at (3,10) {$K$};
  	\node (L) at (3,8) {$L$};
  	\node (M) at (1,8) {$M$};
  	\node (N) at (13,10) {$N$};
  	\node (O) at (13,8) {$O$};
  	\node (P) at (15,8) {$P$};

  \foreach \from/\to in {A/B, B/C, C/D, D/E, E/A, A/L, B/G, C/O, D/I, E/J, O/P, O/N, M/J, J/I, I/P, P/N, N/G, G/K, K/M,  K/L, L/M}
    	\draw (\from) -- (\to);

\end{tikzpicture}      
\caption{}
\label{Graph 2}
\end{center}
\end{minipage}
\begin{minipage}{.32\textwidth}

\begin{center}
\begin{tikzpicture}
  [scale=.28,auto=left,every node/.style={circle,draw, fill=black!20, scale=0.5}]
  	\node (A) at (5,7) {$A$};
  	\node (B) at (8,9) {$B$};
  	\node (C) at (11,7) {$C$};
  	\node (D) at (10,4) {$D$};
  	\node (E) at (6,4) {$E$};
  	\node (G) at (8,12) {$G$};
  	\node (I) at (12,1) {$I$};
  	\node (J) at (4,1) {$J$};
  	\node (K) at (3,10) {$K$};
  	\node (L) at (3,8) {$L$};
  	\node (M) at (1,8) {$M$};
  	\node (N) at (13,10) {$N$};
  	\node (O) at (13,8) {$O$};
  	\node (P) at (15,8) {$P$};
  	\node (Q) at (14,9) {$Q$};
  	\node (R) at (2,9) {$R$};

  \foreach \from/\to in {A/B, B/C, C/D, D/E, E/A, A/L, B/G, C/O, D/I, E/J, O/P, O/N, M/J, J/I, I/P, P/Q, Q/N, N/G, G/K, K/R, R/M,  K/L, L/M}
    	\draw (\from) -- (\to);

  	\draw (R) to [bend left=100] (Q); 
    
\end{tikzpicture}

\caption{}
\label{Graph 3}
\end{center}
\end{minipage}
\end{figure}

Now it is easy to see that the graph in Figure~\ref{Graph 3} is isomorphic to the graph in Figure~\ref{Graph 4} which contains Tutte's fragment. A Hamiltonian path through Tutte's fragment must be a Hamiltonian cycle in Figure~\ref{Graph 4}, thus it must use edge $a$. 
\end{proof}

Now we put three of these fragments together like in Figure~\ref{Tutte's graph} . A Hamiltonian cycle through the whole graph must use the three 'top-edges'  of the fragments, but then the three edges of the center vertex must be contained in the cycle, which is not possible. Note that a Hamiltionian cycle must visit every vertex of a fragment before going to the next fragment. So there can not exist a Hamiltionian cycle in Tutte's graph.

\begin{figure}[h]
\centering
\begin{minipage}{.45\textwidth}

\begin{center}
\begin{tikzpicture}
  [scale=.3,auto=left,every node/.style={circle, draw ,fill=black!20, scale=0.5}]
  	\node (A) at (6,4) {$A$};
  	\node (B) at (8,6) {$B$};
  	\node (C) at (6,8) {$C$};
  	\node (D) at (4,7) {$D$};
  	\node (E) at (4,5) {$E$};
  	\node (G) at (11,7) {$G$};
  	\node (I) at (2,8) {$I$};
  	\node (J) at (1,5) {$J$};
  	\node (K) at (11,2) {$K$};
  	\node (L) at (6,1) {$L$};
  	\node (M) at (1,2) {$M$};
  	\node (N) at (9,11) {$N$};
  	\node (O) at (6,10) {$O$};
  	\node (P) at (3,11) {$P$};
  	\node (Q) at (6,13) {$Q$};
  	\node (R) at (6,16) {$R$};

  \foreach \from/\to in {A/B, B/C, C/D, D/E, E/A, A/L, B/G, C/O, D/I, E/J, O/P, O/N, M/J, J/I, I/P, P/Q, Q/R, Q/N, N/G, G/K, K/L, L/M}
    	\draw (\from) -- (\to);

  	\draw (K) to [bend right=45] (R); 
  	\draw (M) to [bend left=45] (R);

\end{tikzpicture}
\caption{}
\label{Graph 4}
\end{center}
\end{minipage}
\begin{minipage}{.45\textwidth}

\begin{center}
\begin{tikzpicture}
  [scale=.22,auto=left,every node/.style={circle, draw ,fill=black!20, scale=0.5}]

	\draw (14,16) -- (14,8) [dashed];
	\draw (7,4) -- (14,8) [dashed];
	\draw (21,4) -- (14,8) [dashed];

  	\node (X) at (0,0) {$X$};
  	\node (Y) at (28,0) {$Y$};
  	\node (Z) at (14,8) {$Z$};

  	\node (G) at (17,8) {$G$};
  	\node (I) at (12.5,9) {$I$};
  	\node (J) at (10,6.1) {$J$};
  	\node (K) at (21,4) {$K$};
  	\node (L) at (14,6) {$L$};
  	\node (M) at (7,4) {$M$};
  	\node (N) at (15,13) {$N$};
  	\node (P) at (13.5,13) {$P$};
  	\node (Q) at (14,16) {$Q$};
  	\node (R) at (14,24) {$R$};

  \foreach \from/\to in {X/Y, Y/R, R/X, X/M, Y/K, M/J, J/I, I/P, P/Q, Q/R, Q/N, N/G, G/K, K/L, L/M}
    	\draw (\from) -- (\to);

\end{tikzpicture}
\caption{Tetrahedron $RXYZ$ (where $RXY$ lies in the plane of the page and vertex $Z$ lies above the plane of the page). }
\label{Tetrahedron}
\end{center}
\end{minipage}
\end{figure}

We check that Tutte's graph really is a counterexample tot Tait's conjecture. It is easy to see that it is cubic. Because of Steinitz's theorem it is enough to show that the graph corresponds to a convex 3-dimensional polytope. In this case it fits in a tetrahedron. The three faces between each pair of fragments correspond to three faces of the tetrahedron, the exterior of the graph is the fourth face. In  Figure~\ref{Tetrahedron} it is shown how some of the vertices of one Tutte fragment can be placed on the faces of the tetrahedron $XYZR$. 

When one of the tops of the tetrahedron is cut off, there  arises a triangle. The pentagon $ABCDE$ is placed on the interior of this triangle and the vertices are connected to the vertices on the original faces of the tetrahedron. Note that the vertices of the triangle are not part of the polyhedron we are constructing. Only the edge leaving from $C$ should not be connected directly to vertices $N$ and $P$, but has to make an extra snap and then in point $O$ be split into two edges.

This procedure should be repeated in two other tops, the fourth vertex of the original tetrahedron corresponds to the center vertex of Tutte's graph. When all vertices are placed in the tetrahedron in a convex way the graph is polyhedral. This shows that Tutte's graph really is a counterexample to Tait's conjecture.

\section{Smaller counterexamples}

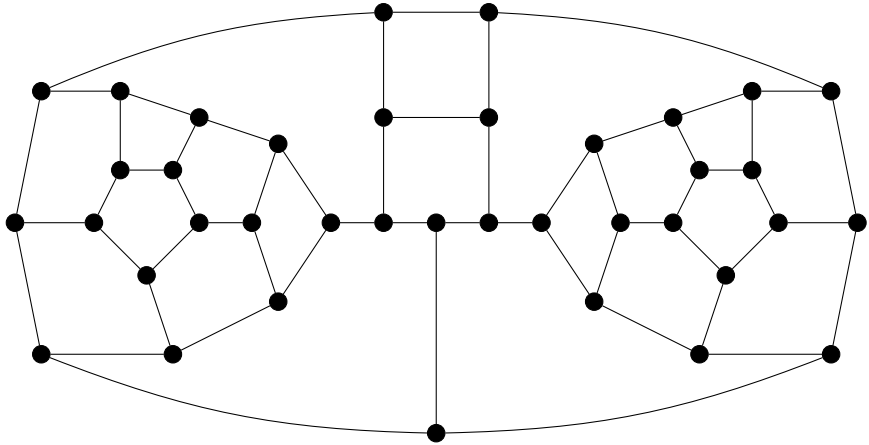
\begin{figure}[H]
\begin{center}
\begin{tikzpicture}
  [scale=.35,auto=left,every node/.style={circle,draw, fill=black!100,scale=0.7}]
  	\node (A1) at (4,9) {};
  	\node (B1) at (6,7) {};
  	\node (C1) at (8,9) {};
  	\node (D1) at (7,11) {};
  	\node (E1) at (5,11) {};
  	\node (G1) at (7,4) {};
  	\node (I1) at (8,13) {};
  	\node (J1) at (5,14) {};
  	\node (K1) at (2,4) {};
  	\node (L1) at (1,9) {};
  	\node (M1) at (2,14) {};
  	\node (N1) at (11,6) {};
  	\node (O1) at (10,9) {};
  	\node (P1) at (11,12) {};
  	\node (Q1) at (13,9) {};
  	\node (R1) at (15,9) {};

  \foreach \from/\to in {A1/B1, B1/C1, C1/D1, D1/E1, E1/A1, A1/L1, B1/G1, C1/O1, D1/I1, E1/J1, O1/P1, O1/N1, M1/J1, J1/I1, I1/P1, P1/Q1, Q1/N1, N1/G1, G1/K1, K1/L1, L1/M1, Q1/R1}
    	\draw (\from) -- (\to);

	\node (A2) at (30,9) {};
  	\node (B2) at (28,7) {};
  	\node (C2) at (26,9) {};
  	\node (D2) at (27,11) {};
  	\node (E2) at (29,11) {};
  	\node (G2) at (27,4) {};
  	\node (I2) at (26,13) {};
  	\node (J2) at (29,14) {};
  	\node (K2) at (32,4) {};
  	\node (L2) at (33,9) {};
  	\node (M2) at (32,14) {};
  	\node (N2) at (23,6) {};
  	\node (O2) at (24,9) {};
  	\node (P2) at (23,12) {};
  	\node (Q2) at (21,9) {};
  	\node (R2) at (19,9) {};

  \foreach \from/\to in {A2/B2, B2/C2, C2/D2, D2/E2, E2/A2, A2/L2, B2/G2, C2/O2, D2/I2, E2/J2, O2/P2, O2/N2, M2/J2, J2/I2, I2/P2, P2/Q2, Q2/N2, N2/G2, G2/K2, K2/L2, L2/M2, Q2/R2}
    	\draw (\from) -- (\to);

  	\node (X1) at (15,13) {};
  	\node (X2) at (19,13) {};
  	\node (Y1) at (15,17) {};
  	\node (Y2) at (19,17) {};
  	\node (Z1) at (17,9) {};
  	\node (Z2) at (17,1) {};

  \foreach \from/\to in {R1/Z1, R1/X1, X1/Y1, Y1/Y2, X2/Y2, R2/X2, Z1/R2, Z1/Z2, X1/X2}
    	\draw (\from) -- (\to);

  	\draw (M1) to [bend left=10] (Y1); 
  	\draw (M2) to [bend right=10] (Y2); 
  	\draw (K1) to [bend right=10] (Z2); 
  	\draw (K2) to [bend left=10] (Z2); 
     
\end{tikzpicture}
\caption{Barnette-Bos\'ak-Lederberg Graph}
\label{Barnette-Bosak-Lederberg Graph}
\end{center}
\end{figure}

Several smaller counterexamples were found. In 1965 Lederberg found a counterexample on 38 vertices, around the same time as Barnette and Bos\'ak. This graph is known as the Barnette-Bos\'ak-Lederberg graph and is given in Figure~\ref{Barnette-Bosak-Lederberg Graph}. Two Tutte fragments can be recognized. 


In 1988 Holton and McKay proved no smaller counterexample exists.\cite{Holton1988}

\section{Four color theorem for cubic 3-connected planar graphs}

If Tait's conjecture would have been true, it would imply a proof of the four color theorem for cubic 3-connected planar graphs \cite{TheIcosianGame}.

Suppose we have a cubic polyhedral graph with a Hamiltonian cycle; an example is given in Figure~\ref{fct Hamiltonian cycle}. We can alternately color the edges of the Hamiltionian cycle blue and purple, because the number of vertices is even for a Tait graph. Every vertex has one uncolored edge, which can be colored red. See Figure~\ref{fct Colored Hamiltonian cycle}.

\begin{figure}
\centering
\begin{minipage}{.45\textwidth}

\begin{center}
 \begin{tikzpicture}
  [scale=.35,auto=left,every node/.style={circle,draw, fill=black!100,scale=0.7}]
  	\node (A) at (5,7) {};
  	\node (B) at (8,9) {};
  	\node (C) at (11,7) {};
  	\node (D) at (10,4) {};
  	\node (E) at (6,4) {};
  	\node (G) at (8,12) {};
  	\node (I) at (12,1) {};
  	\node (J) at (4,1) {};
  	\node (H) at (14,9) {};
  	\node (F) at (2,9) {};

  \foreach \from/\to in {A/B, C/D, D/E, E/A, B/G, C/H, J/I, I/H, G/F, F/J}
    	\draw [very thick](\from) -- (\to);

 \foreach \from/\to in {B/C, A/F,  D/I, E/J, H/G}
    	\draw [very thin, dashed](\from) -- (\to);

\end{tikzpicture}
\caption{Hamiltonian cycle}
\label{fct Hamiltonian cycle}
\end{center}
\end{minipage}
\begin{minipage}{.45\textwidth}

\begin{center}
 \begin{tikzpicture}
  [scale=.35,auto=left,every node/.style={circle,draw, fill=black!100,scale=0.7}]
  	\node (A) at (5,7) {};
  	\node (B) at (8,9) {};
  	\node (C) at (11,7) {};
  	\node (D) at (10,4) {};
  	\node (E) at (6,4) {};
  	\node (G) at (8,12) {};
  	\node (I) at (12,1) {};
  	\node (J) at (4,1) {};
  	\node (H) at (14,9) {};
  	\node (F) at (2,9) {};

  \foreach \from/\to in {C/D, E/A, B/G, I/H, F/J}
    	\draw [very thick, blue](\from) -- (\to);

  \foreach \from/\to in {A/B, D/E, C/H, J/I, G/F}
    	\draw [very thick, violet](\from) -- (\to);

 \foreach \from/\to in {B/C, A/F,  D/I, E/J, H/G}
    	\draw [red, dashed](\from) -- (\to);

\end{tikzpicture}
\caption{Colored Hamiltonian cycle}
\label{fct Colored Hamiltonian cycle}
\end{center}

\end{minipage}
\end{figure}

Now step 1 is to throw out the red edges to obtain a polygon consisting of a subset of faces of the graph. The whole polygon is colored blue, as shown in Figure~\ref{fct blue polygon}. Step 2 is to throw out the blue edges of  Figure~\ref{fct Colored Hamiltonian cycle} to obtain one or more polygons which are colored red, see  Figure~\ref{fct red polygons}. By overlaying the two colorings, we obtain a four-coloring, where a face is colored purple when it is colored blue in the first step and red in the second step. See Figure~\ref{fct coloring}.

\begin{figure}[h]
\centering
\begin{minipage}{.32\textwidth}

\begin{center}
 \begin{tikzpicture}
  [scale=.35,auto=left,every node/.style={circle,draw, fill=black!100,scale=0.7}]

\draw [fill=blue,opacity=0.2] (8,12)--(2,9)--(4,1)--(12,1)--(14,9)--(11,7)--(10,4)--(6,4)--(5,7)--(8,9)--(8,12);

  	\node (A) at (5,7) {};
  	\node (B) at (8,9) {};
  	\node (C) at (11,7) {};
  	\node (D) at (10,4) {};
  	\node (E) at (6,4) {};
  	\node (G) at (8,12) {};
  	\node (I) at (12,1) {};
  	\node (J) at (4,1) {};
  	\node (H) at (14,9) {};
  	\node (F) at (2,9) {};

  \foreach \from/\to in {C/D, E/A, B/G, I/H, F/J}
    	\draw [very thick, blue](\from) -- (\to);

  \foreach \from/\to in {A/B, D/E, C/H, J/I, G/F}
    	\draw [very thick, violet](\from) -- (\to);

\end{tikzpicture}
\caption{Blue polygon}
\label{fct blue polygon}
\end{center}

\end{minipage}
\begin{minipage}{.32\textwidth}

\begin{center}
 \begin{tikzpicture}
  [scale=.35,auto=left,every node/.style={circle,draw, fill=black!100,scale=0.7}]

\draw [fill=red,opacity=0.2] (5,7)--(8,9)--(11,7)--(14,9)--(8,12)--(2,9)--(5,7);

\draw [fill=red,opacity=0.2] (6,4)--(10,4)--(12,1)--(4,1)--(6,4);

  	\node (A) at (5,7) {};
  	\node (B) at (8,9) {};
  	\node (C) at (11,7) {};
  	\node (D) at (10,4) {};
  	\node (E) at (6,4) {};
  	\node (G) at (8,12) {};
  	\node (I) at (12,1) {};
  	\node (J) at (4,1) {};
  	\node (H) at (14,9) {};
  	\node (F) at (2,9) {};

  \foreach \from/\to in {A/B, D/E, C/H, J/I, G/F}
    	\draw [very thick, violet](\from) -- (\to);

 \foreach \from/\to in {B/C, A/F,  D/I, E/J, H/G}
    	\draw [red, dashed](\from) -- (\to);

\end{tikzpicture}
\caption{Red polygons}
\label{fct red polygons}
\end{center}

\end{minipage}
\begin{minipage}{.32\textwidth}

\begin{center}
 \begin{tikzpicture}
  [scale=.35,auto=left,every node/.style={circle,draw, fill=black!100,scale=0.7}]

\draw [fill=blue,opacity=0.2] (8,12)--(2,9)--(4,1)--(12,1)--(14,9)--(11,7)--(10,4)--(6,4)--(5,7)--(8,9)--(8,12);

\draw [fill=red,opacity=0.2] (5,7)--(8,9)--(11,7)--(14,9)--(8,12)--(2,9)--(5,7);

\draw [fill=red,opacity=0.2] (6,4)--(10,4)--(12,1)--(4,1)--(6,4);

  	\node (A) at (5,7) {};
  	\node (B) at (8,9) {};
  	\node (C) at (11,7) {};
  	\node (D) at (10,4) {};
  	\node (E) at (6,4) {};
  	\node (G) at (8,12) {};
  	\node (I) at (12,1) {};
  	\node (J) at (4,1) {};
  	\node (H) at (14,9) {};
  	\node (F) at (2,9) {};

  \foreach \from/\to in {A/B, D/E, C/H, J/I, G/F}
    	\draw [very thick, violet](\from) -- (\to);

 \foreach \from/\to in {B/C, A/F,  D/I, E/J, H/G}
    	\draw [red, dashed](\from) -- (\to);

  \foreach \from/\to in {C/D, E/A, B/G, I/H, F/J}
    	\draw [very thick, blue](\from) -- (\to);

\end{tikzpicture}
\caption{Coloring}
\label{fct coloring}
\end{center}

\end{minipage}
\end{figure}

To see that this really is a four-coloring, consider two adjacent faces. Suppose the edge between these faces is colored blue. Then this edge is a side of the blue polygon, so in the first step one face will be colored blue and the other is not colored. Whatever happens in step two, overlaying these faces with red or uncolored faces can never give the two faces the same color. 

In the other cases we get the same conclusion that adjacent faces can never have the same color, so this procedure really gives a four coloring.

\chapter{Tutte's Conjecture}
\label{Tutte}
None of the known counterexamples to Tait's conjecture is bipartite, so Tutte formulated his own conjecture in 1971 \cite{Tutte1971}. 

\begin{conjecture}[Tutte's conjecture]
Every cubic 3-connected bipartite graph contains a Hamiltonian cycle.
\end{conjecture}

The first counterexample was found in 1976 by Horton and has 96 vertices. This Horton-96 graph first appeared in \cite{BM1976}, where it is left as an exercise to show that the graph is non-Hamiltonian. 

Later smaller counterexamples were found:
\begin{itemize}
\item In 1982 Horton found the Horton-92 graph with 92 vertices \cite{Horton1982}.
\item In 1983 Owens found a graph on 78 vertices \cite{Owens1983}.
\item Ellingham and Horton published two counterexamples: The Ellingham-Horton graph with 78 vertices in 1981 \cite{EH1981} and the Ellingham-Horton graph with 54 vertices in 1983 \cite{EH1983}.
\end{itemize}

The smallest counterexample currently known is the Georges's graph on 50 vertices, found in 1989  \cite{Georges1989}. 

The existence of these counterexamples shows that planarity in Barnette's conjecture is a necessary condition. 

In the first section of this chapter the Horton fragment is considered, which plays an important role in the Horton-96 and the Horton-92 graph. The Ellingham fragment in the second section is the main ingredient of Georges's graph.

\noindent \begin{minipage}{0.7\textwidth}

\section{Horton fragment}

The Horton graphs consist of Horton fragments; one of these fragments is shown in Figure~\ref{Horton fragment}. Each fragment contains two Horton circles; one of these circles is shown in Figure~\ref{Horton circle}. 

\begin{lemma}
No Hamiltonian cycle in a Horton circle contains exactly one of the red edges ($e_1$ and $e_2$) in Figure~\ref{Horton circle}.
\label{Hortoncirle lemma1}
\end{lemma}

\end{minipage}
\begin{minipage}{0.3\textwidth}


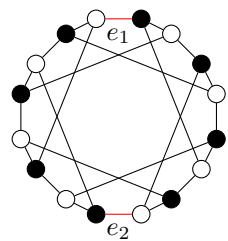
\begin{figure}[H]
\centering
\begin{tikzpicture}
  [scale=.2,auto=left,every node/.style={circle, draw, fill=white, scale=.7}]
  	\node [black](A) at (11,16) {};
  	\node (B) at (13,15) {};
  	\node [black](C) at (15,13) {};
  	\node (D) at (16,11) {};
  	\node [black](E) at (16,8) {};
  	\node (F) at (15,6) {};
  	\node [black](G) at (13,4) {};
  	\node (H) at (11,3) {};
  	\node [black](I) at (8,3) {};
  	\node (J) at (6,4) {};
  	\node [black](K) at (4,6) {};
  	\node (L) at (3,8) {};
  	\node [black](M) at (3,11) {};
  	\node (N) at (4,13) {};
  	\node [black](O) at (6,15) {};
  	\node (P) at (8,16) {};

  \foreach \from/\to in {A/B, B/C, C/D, D/E, E/F, F/G, G/H, I/J, J/K, K/L, L/M, M/N, N/O, O/P, 
			 A/F, B/M, C/H, D/O, E/J, G/L, I/N, K/P}
    	\draw (\from) -- (\to);

  \path[every node/.style={font=\sffamily}]
	(A) edge [red] node [black] {$e_1$} (P)
	(H) edge [red] node [black] {$e_2$} (I);

\end{tikzpicture}
\caption{Horton circle}
\label{Horton circle}
\end{figure}

\end{minipage}

\begin{proof}

Suppose there is a Hamiltonian cycle that uses exactly one of the red edges. Without loss of generality we can assume $e_1$ is used and $e_2$ is not used. In the following figures thick edges will represent edges that are on the Hamiltonian cycle, the dashed edges will be the ones that are not used.

\begin{itemize}
\item Suppose exactly one of the other edges between the right and the left half of the circle is used. There are two possiblities in the red vertex:

\begin{itemize}
\item Go up. Then we have the situation in Figure~\ref{Horton circle 2}. Because we used one edge to cross from right to left, the other three crossing-edges may not be used by assumption, so we get the situation given in Figure~\ref{Horton circle 3}. The vertices with one dashed edge now must use the other two edges, which is shown in Figure~\ref{Horton circle 4}. This gives two cycles, while we want one in the whole graph.

\begin{figure}[h]
\centering
\begin{minipage}{.32\textwidth}


\begin{center}
\begin{tikzpicture}
 [scale=.25,auto=left,every node/.style={circle,draw, fill=black!100, scale=.7}]
  	\node (A) at (11,16) {};
  	\node [fill=white](B) at (13,15) {};
  	\node (C) [fill=red] at (15,13) {};
  	\node [fill=white](D) at (16,11) {};
  	\node (E) at (16,8) {};
  	\node [fill=white](F) at (15,6) {};
  	\node (G) at (13,4) {};
  	\node [fill=white](H) at (11,3) {};
  	\node (I) at (8,3) {};
  	\node [fill=white](J) at (6,4) {};
  	\node (K) at (4,6) {};
  	\node [fill=white](L) at (3,8) {};
  	\node (M) at (3,11) {};
  	\node [fill=white](N) at (4,13) {};
  	\node (O) at (6,15) {};
  	\node [fill=white](P) at (8,16) {};

  \foreach \from/\to in {A/B, E/F, F/G, G/H, I/J, J/K, K/L, L/M, M/N, N/O, O/P, 
			A/F, C/H, I/N, K/P, B/M, F/A, G/L, E/J}
    	\draw (\from) -- (\to);

  \foreach \from/\to in {A/P, G/H, I/J, C/H, I/N, B/C, D/E, D/O} 
    	\draw [very thick](\from) -- (\to);

 \foreach \from/\to in {H/I, C/D}
    	\draw [dashed](\from) -- (\to);

\end{tikzpicture}
\caption{}
\label{Horton circle 2}
\end{center}

\end{minipage}
\begin{minipage}{.32\textwidth}


\begin{center}
\begin{tikzpicture}
 [scale=.25,auto=left,every node/.style={circle,draw, fill=black!100, scale=.7}]
  	\node (A) at (11,16) {};
  	\node [fill=white](B) at (13,15) {};
  	\node (C) at (15,13) {};
  	\node [fill=white](D) at (16,11) {};
  	\node (E) at (16,8) {};
  	\node [fill=white](F) at (15,6) {};
  	\node (G) at (13,4) {};
  	\node [fill=white](H) at (11,3) {};
  	\node (I) at (8,3) {};
  	\node [fill=white](J) at (6,4) {};
  	\node (K) at (4,6) {};
  	\node [fill=white](L) at (3,8) {};
  	\node (M) at (3,11) {};
  	\node [fill=white](N) at (4,13) {};
  	\node (O) at (6,15) {};
  	\node [fill=white](P) at (8,16) {};

  \foreach \from/\to in {A/B, E/F, F/G, G/H, I/J, J/K, K/L, L/M, M/N, N/O, O/P, 
			A/F, C/H, I/N, K/P, F/A}
    	\draw (\from) -- (\to);

  \foreach \from/\to in {A/P, G/H, I/J, C/H, I/N, B/C, D/E, D/O} 
    	\draw [very thick](\from) -- (\to);

 \foreach \from/\to in {H/I, C/D, B/M, E/J, G/L}
    	\draw [dashed](\from) -- (\to);

\end{tikzpicture}
\caption{}
\label{Horton circle 3}
\end{center}

\end{minipage}
\begin{minipage}{.32\textwidth}


\begin{center}
\begin{tikzpicture}
 [scale=.25,auto=left,every node/.style={circle,draw, fill=black!100, scale=.7}]
  	\node (A) at (11,16) {};
  	\node [fill=white](B) at (13,15) {};
  	\node (C) at (15,13) {};
  	\node [fill=white](D) at (16,11) {};
  	\node (E) at (16,8) {};
  	\node [fill=white](F) at (15,6) {};
  	\node (G) at (13,4) {};
  	\node [fill=white](H) at (11,3) {};
  	\node (I) at (8,3) {};
  	\node [fill=white](J) at (6,4) {};
  	\node (K) at (4,6) {};
  	\node [fill=white](L) at (3,8) {};
  	\node (M) at (3,11) {};
  	\node [fill=white](N) at (4,13) {};
  	\node (O) at (6,15) {};
  	\node [fill=white](P) at (8,16) {};

  \foreach \from/\to in {G/H, I/J, N/O, O/P, 
			A/F, C/H, I/N, K/P, F/A}
    	\draw (\from) -- (\to);

  \foreach \from/\to in {A/P, G/H, I/J, C/H, I/N, B/C, D/E, D/O, A/B, E/F, F/G, J/K, K/L,  L/M, M/N, O/P} 
    	\draw [very thick](\from) -- (\to);

 \foreach \from/\to in {H/I, C/D, B/M, E/J, G/L}
    	\draw [dashed](\from) -- (\to);

\end{tikzpicture}
\caption{}
\label{Horton circle 4}
\end{center}

\end{minipage}
\end{figure}

\item Go down.  Then we have the situation in Figure~\ref{Horton circle 5}. Because we used one edge to cross from right to left, the other three crossing-edges may not be used, so we get the situation given in Figure~\ref{Horton circle 6}. The vertices with one dashed edge now must use the other two edges, which is shown in Figure~\ref{Horton circle 7}. This gives two cycles, while we want one in the whole graph. 

\begin{figure}[h]
\centering
\begin{minipage}{.32\textwidth}


\begin{center}
\begin{tikzpicture}
 [scale=.25,auto=left,every node/.style={circle,draw, fill=black!100, scale=.7}]
  	\node (A) at (11,16) {};
  	\node [fill=white](B) at (13,15) {};
  	\node (C) [fill=red] at (15,13) {};
  	\node [fill=white](D) at (16,11) {};
  	\node (E) at (16,8) {};
  	\node [fill=white](F) at (15,6) {};
  	\node (G) at (13,4) {};
  	\node [fill=white](H) at (11,3) {};
  	\node (I) at (8,3) {};
  	\node [fill=white](J) at (6,4) {};
  	\node (K) at (4,6) {};
  	\node [fill=white](L) at (3,8) {};
  	\node (M) at (3,11) {};
  	\node [fill=white](N) at (4,13) {};
  	\node (O) at (6,15) {};
  	\node [fill=white](P) at (8,16) {};

  \foreach \from/\to in {E/F, F/G, G/H, I/J, J/K, K/L, L/M, M/N, N/O, O/P, 
			A/F, C/H, I/N, K/P, F/A, G/L, E/J, D/O, D/E, B/M}
    	\draw (\from) -- (\to);

  \foreach \from/\to in {A/P, G/H, I/J, C/H, I/N, C/D, A/B, M/B} 
    	\draw [very thick](\from) -- (\to);

 \foreach \from/\to in {H/I, B/C}
    	\draw [dashed](\from) -- (\to);

\end{tikzpicture}
\caption{}
\label{Horton circle 5}
\end{center}

\end{minipage}
\begin{minipage}{.32\textwidth}


\begin{center}
\begin{tikzpicture}
 [scale=.25,auto=left,every node/.style={circle,draw, fill=black!100, scale=.7}]
  	\node (A) at (11,16) {};
  	\node [fill=white](B) at (13,15) {};
  	\node (C) at (15,13) {};
  	\node [fill=white](D) at (16,11) {};
  	\node (E) at (16,8) {};
  	\node [fill=white](F) at (15,6) {};
  	\node (G) at (13,4) {};
  	\node [fill=white](H) at (11,3) {};
  	\node (I) at (8,3) {};
  	\node [fill=white](J) at (6,4) {};
  	\node (K) at (4,6) {};
  	\node [fill=white](L) at (3,8) {};
  	\node (M) at (3,11) {};
  	\node [fill=white](N) at (4,13) {};
  	\node (O) at (6,15) {};
  	\node [fill=white](P) at (8,16) {};

  \foreach \from/\to in {E/F, F/G, G/H, I/J, J/K, K/L, L/M, M/N, N/O, O/P, 
			A/F, C/H, I/N, K/P, F/A, D/E}
    	\draw (\from) -- (\to);

  \foreach \from/\to in {A/P, G/H, I/J, C/H, I/N, C/D, A/B, B/M} 
    	\draw [very thick](\from) -- (\to);

 \foreach \from/\to in {H/I, B/C, D/O,  E/J, G/L}
    	\draw [dashed](\from) -- (\to);

\end{tikzpicture}
\caption{}
\label{Horton circle 6}
\end{center}

\end{minipage}
\begin{minipage}{.32\textwidth}


\begin{center}
\begin{tikzpicture}
 [scale=.25,auto=left,every node/.style={circle,draw, fill=black!100, scale=.7}]
  	\node (A) at (11,16) {};
  	\node [fill=white](B) at (13,15) {};
  	\node (C) at (15,13) {};
  	\node [fill=white](D) at (16,11) {};
  	\node (E) at (16,8) {};
  	\node [fill=white](F) at (15,6) {};
  	\node (G) at (13,4) {};
  	\node [fill=white](H) at (11,3) {};
  	\node (I) at (8,3) {};
  	\node [fill=white](J) at (6,4) {};
  	\node (K) at (4,6) {};
  	\node [fill=white](L) at (3,8) {};
  	\node (M) at (3,11) {};
  	\node [fill=white](N) at (4,13) {};
  	\node (O) at (6,15) {};
  	\node [fill=white](P) at (8,16) {};

  \foreach \from/\to in {E/F, F/G, G/H, I/J, J/K, K/L, L/M, N/O, O/P, 
			A/F, C/H, I/N, K/P, F/A, D/E}
    	\draw (\from) -- (\to);

  \foreach \from/\to in {A/P, G/H, I/J, C/H, I/N, C/D, A/B, B/M, D/E, E/F, F/G, J/K, K/L, L/M, N/O, O/P} 
    	\draw [very thick](\from) -- (\to);

 \foreach \from/\to in {H/I, B/C, D/O,  E/J, G/L, M/N}
    	\draw [dashed](\from) -- (\to);

\end{tikzpicture}
\caption{}
\label{Horton circle 7}
\end{center}

\end{minipage}
\end{figure}

\end{itemize}

\item Suppose three of the edges between the right and the left part of the circle are used. Because of symmetry there are two possibilities:

\begin{itemize}
\item One of the upper crossing edges is not used. Without loss of generality we can assume the dashed edge in Figure~\ref{Horton circle 8} is not used and the other three are used. The vertices with two thick edges can not use the other edge and the vertices with a dashed edge must use the other two edges, so we get the situation of Figure~\ref{Horton circle 9}. Now consider the red edge in  Figure~\ref{Horton circle 10}. On one side it is connected to a vertex that already used two edges, so the red edge can not be on the Hamitonian cylce, while on the other end there is a dashed edge so the red edge must be used. This gives a contradiction.

\begin{figure}[h]
\centering
\begin{minipage}{.32\textwidth}


\begin{center}
\begin{tikzpicture}
 [scale=.25,auto=left,every node/.style={circle,draw, fill=black!100, scale=.7}]
  	\node (A) at (11,16) {};
  	\node [fill=white](B) at (13,15) {};
  	\node (C) at (15,13) {};
  	\node [fill=white](D) at (16,11) {};
  	\node (E) at (16,8) {};
  	\node [fill=white](F) at (15,6) {};
  	\node (G) at (13,4) {};
  	\node [fill=white](H) at (11,3) {};
  	\node (I) at (8,3) {};
  	\node [fill=white](J) at (6,4) {};
  	\node (K) at (4,6) {};
  	\node [fill=white](L) at (3,8) {};
  	\node (M) at (3,11) {};
  	\node [fill=white](N) at (4,13) {};
  	\node (O) at (6,15) {};
  	\node [fill=white](P) at (8,16) {};

  \foreach \from/\to in {A/B, B/C, C/D, D/E, E/F, F/G, G/H, I/J, J/K, K/L, L/M, M/N, N/O, O/P, 
			 A/F, C/H, I/N, K/P}
    	\draw (\from) -- (\to);

  \foreach \from/\to in {A/P, G/H, I/J, C/H, I/N, D/O, E/J, G/L} 
    	\draw [very thick](\from) -- (\to);

 \foreach \from/\to in {H/I, B/M}
    	\draw [dashed](\from) -- (\to);

\end{tikzpicture}
\caption{}
\label{Horton circle 8}
\end{center}

\end{minipage}
\begin{minipage}{.32\textwidth}


\begin{center}
\begin{tikzpicture}
 [scale=.25,auto=left,every node/.style={circle,draw, fill=black!100, scale=.7}]
  	\node (A) at (11,16) {};
  	\node [fill=white](B) at (13,15) {};
  	\node (C) at (15,13) {};
  	\node [fill=white](D) at (16,11) {};
  	\node (E) at (16,8) {};
  	\node [fill=white](F) at (15,6) {};
  	\node (G) at (13,4) {};
  	\node [fill=white](H) at (11,3) {};
  	\node (I) at (8,3) {};
  	\node [fill=white](J) at (6,4) {};
  	\node (K) at (4,6) {};
  	\node [fill=white](L) at (3,8) {};
  	\node (M) at (3,11) {};
  	\node [fill=white](N) at (4,13) {};
  	\node (O) at (6,15) {};
  	\node [fill=white](P) at (8,16) {};

  \foreach \from/\to in {A/B, B/C, C/D, D/E, E/F, G/H, I/J, K/L, L/M, M/N, N/O, O/P, 
			 A/F, C/H, I/N, K/P}
    	\draw (\from) -- (\to);

  \foreach \from/\to in {A/P, G/H, I/J, C/H, I/N, D/O, E/J, G/L, A/B, B/C, L/M, M/N} 
    	\draw [very thick](\from) -- (\to);

 \foreach \from/\to in {H/I, B/M, F/G, J/K}
    	\draw [dashed](\from) -- (\to);

\end{tikzpicture}
\caption{}
\label{Horton circle 9}
\end{center}

\end{minipage}
\begin{minipage}{.32\textwidth}


\begin{center}
\begin{tikzpicture}
 [scale=.25,auto=left,every node/.style={circle,draw, fill=black!100, scale=.7}]
  	\node (A) at (11,16) {};
  	\node [fill=white](B) at (13,15) {};
  	\node (C) at (15,13) {};
  	\node [fill=white](D) at (16,11) {};
  	\node (E) at (16,8) {};
  	\node [fill=white](F) at (15,6) {};
  	\node (G) at (13,4) {};
  	\node [fill=white](H) at (11,3) {};
  	\node (I) at (8,3) {};
  	\node [fill=white](J) at (6,4) {};
  	\node (K) at (4,6) {};
  	\node [fill=white](L) at (3,8) {};
  	\node (M) at (3,11) {};
  	\node [fill=white](N) at (4,13) {};
  	\node (O) at (6,15) {};
  	\node [fill=white](P) at (8,16) {};

  \foreach \from/\to in {A/B, B/C, C/D, D/E, E/F, G/H, I/J, K/L, L/M, M/N, N/O, O/P, 
			 C/H, I/N, K/P}
    	\draw (\from) -- (\to);

  \foreach \from/\to in {A/P, G/H, I/J, C/H, I/N, D/O, E/J, G/L, A/B, B/C, L/M, M/N} 
    	\draw [very thick](\from) -- (\to);

 \foreach \from/\to in {H/I, B/M, F/G, J/K}
    	\draw [dashed](\from) -- (\to);

  \foreach \from/\to in {A/F}
    	\draw (\from) -- (\to)[red];

\end{tikzpicture}
\caption{}
\label{Horton circle 10}
\end{center}

\end{minipage}
\end{figure}

\item One of the lower crossing edges is not used. Without loss of generality we can assume the dashed edge in Figure~\ref{Horton circle 11} is not used and the other three are used as is shown in Figure~\ref{Horton circle 12}.  The vertices with two thick edges can not use the other edge and the vertices with a dashed edge must use the other two edges. The red edge can not be used, because this would give a cycle which does not visit each vertex. So we get the situation in Figure~\ref{Horton circle 13} were we have two cycles, which is not what we wanted. 
\end{itemize}

\begin{figure}[h]
\centering
\begin{minipage}{.32\textwidth}


\begin{center}
\begin{tikzpicture}
 [scale=.25,auto=left,every node/.style={circle,draw, fill=black!100, scale=.7}]
  	\node (A) at (11,16) {};
  	\node [fill=white](B) at (13,15) {};
  	\node (C) at (15,13) {};
  	\node [fill=white](D) at (16,11) {};
  	\node (E) at (16,8) {};
  	\node [fill=white](F) at (15,6) {};
  	\node (G) at (13,4) {};
  	\node [fill=white](H) at (11,3) {};
  	\node (I) at (8,3) {};
  	\node [fill=white](J) at (6,4) {};
  	\node (K) at (4,6) {};
  	\node [fill=white](L) at (3,8) {};
  	\node (M) at (3,11) {};
  	\node [fill=white](N) at (4,13) {};
  	\node (O) at (6,15) {};
  	\node [fill=white](P) at (8,16) {};

  \foreach \from/\to in {A/B, B/C, C/D, D/E, E/F, F/G, G/H, I/J, J/K, K/L, L/M, M/N, N/O, O/P, 
			 A/F, B/M, C/H, D/O, E/J, I/N, K/P}
    	\draw (\from) -- (\to);

  \foreach \from/\to in {A/P, G/H, I/J, C/H, I/N, B/M, D/O, E/J} 
    	\draw [very thick](\from) -- (\to);

 \foreach \from/\to in {H/I, G/L}
    	\draw [dashed](\from) -- (\to);

\end{tikzpicture}
\caption{}
\label{Horton circle 11}
\end{center}

\end{minipage}
\begin{minipage}{.32\textwidth}


\begin{center}
\begin{tikzpicture}
 [scale=.25,auto=left,every node/.style={circle,draw, fill=black!100, scale=.7}]
  	\node (A) at (11,16) {};
  	\node [fill=white](B) at (13,15) {};
  	\node (C) at (15,13) {};
  	\node [fill=white](D) at (16,11) {};
  	\node (E) at (16,8) {};
  	\node [fill=white](F) at (15,6) {};
  	\node (G) at (13,4) {};
  	\node [fill=white](H) at (11,3) {};
  	\node (I) at (8,3) {};
  	\node [fill=white](J) at (6,4) {};
  	\node (K) at (4,6) {};
  	\node [fill=white](L) at (3,8) {};
  	\node (M) at (3,11) {};
  	\node [fill=white](N) at (4,13) {};
  	\node (O) at (6,15) {};
  	\node [fill=white](P) at (8,16) {};

  \foreach \from/\to in {B/C, C/D, D/E, E/F, F/G, G/H, I/J, K/L, L/M, N/O, 
			 A/F, B/M, C/H, D/O, E/J, I/N}
    	\draw (\from) -- (\to);

  \foreach \from/\to in {A/P, G/H, I/J, C/H, I/N, B/M, D/O, E/J, K/L, L/M, F/G, N/O, K/P} 
    	\draw [very thick](\from) -- (\to);

 \foreach \from/\to in {H/I, G/L, J/K, M/N, O/P}
    	\draw [dashed](\from) -- (\to);

  \foreach \from/\to in {A/B}
    	\draw (\from) -- (\to)[red];

\end{tikzpicture}
\caption{}
\label{Horton circle 12}
\end{center}

\end{minipage}
\begin{minipage}{.32\textwidth}


\begin{center}
\begin{tikzpicture}
 [scale=.25,auto=left,every node/.style={circle,draw, fill=black!100, scale=.7}]
  	\node (A) at (11,16) {};
  	\node [fill=white](B) at (13,15) {};
  	\node (C) at (15,13) {};
  	\node [fill=white](D) at (16,11) {};
  	\node (E) at (16,8) {};
  	\node [fill=white](F) at (15,6) {};
  	\node (G) at (13,4) {};
  	\node [fill=white](H) at (11,3) {};
  	\node (I) at (8,3) {};
  	\node [fill=white](J) at (6,4) {};
  	\node (K) at (4,6) {};
  	\node [fill=white](L) at (3,8) {};
  	\node (M) at (3,11) {};
  	\node [fill=white](N) at (4,13) {};
  	\node (O) at (6,15) {};
  	\node [fill=white](P) at (8,16) {};

  \foreach \from/\to in {B/C, C/D, D/E, E/F, F/G, G/H, I/J, K/L, L/M, N/O, 
			 A/F, B/M, C/H, D/O, E/J, I/N}
    	\draw (\from) -- (\to);

  \foreach \from/\to in {A/P, G/H, I/J, C/H, I/N, B/M, D/O, E/J, K/L, L/M, F/G, N/O, K/P, A/F, B/C} 
    	\draw [very thick](\from) -- (\to);

 \foreach \from/\to in {H/I, G/L, J/K, M/N, O/P, A/B}
    	\draw [dashed](\from) -- (\to);

  \foreach \from/\to in {}
    	\draw (\from) -- (\to)[red];

\end{tikzpicture}
\caption{}
\label{Horton circle 13}
\end{center}

\end{minipage}
\end{figure}

These are all the possibilities, because these are all situations in which an even number of edges is used to cross from right to left (note that the upper edge is also a crossing edge). Thus no Hamiltionian cycle can use exactly one of the red edges. 

\end{itemize}
\end{proof}

\noindent \begin{minipage}{0.5\textwidth}
Now we connect two Horton circles as shown in Figure~\ref{Horton fragment}.

\begin{lemma}
Every Hamiltonian cycle in a Horton fragment contains the red edge $e$ (depicted in Figure~\ref{Horton fragment}).
\label{Hortonfragment lemma2}
\end{lemma}

\end{minipage}
\begin{minipage}{0.6\textwidth}


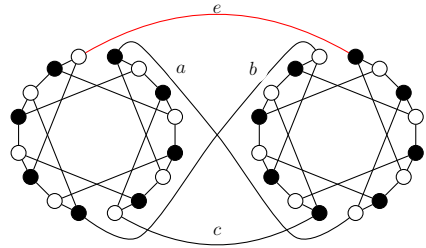
\begin{figure}[H]
\centering
\begin{tikzpicture}
  [scale=.16,auto=left,every node/.style={circle,draw, fill=white, scale=.6}]
\draw plot [smooth] coordinates { (8,3) (13,1) (19.5,9.5) (26,17) (28,16)};
\draw plot [smooth] coordinates { (11,16) (13,17) (19.5,9.5) (26,1) (31,3)};

\node [draw=none, font=\large] (a) at (16.5,15) {$a$};
\node [draw=none, font=\large] (b) at (22.5,15) {$b$};
\node [draw=none, font=\large] (c) at (19.5,1.5) {$c$};
\node [draw=none, font=\large] (c) at (19.5,20) {$e$};

  	\node [black](A1) at (11,16) {};
  	\node (B1) at (13,15) {};
  	\node [black](C1) at (15,13) {};
  	\node (D1) at (16,11) {};
  	\node [black](E1) at (16,8) {};
  	\node (F1) at (15,6) {};
  	\node [black](G1) at (13,4) {};
  	\node (H1) at (11,3) {};
  	\node [black](I1) at (8,3) {};
  	\node (J1) at (6,4) {};
  	\node [black](K1) at (4,6) {};
  	\node (L1) at (3,8) {};
  	\node [black](M1) at (3,11) {};
  	\node (N1) at (4,13) {};
  	\node [black](O1) at (6,15) {};
  	\node (P1) at (8,16) {};

  \foreach \from/\to in {A1/B1, B1/C1, C1/D1, D1/E1, E1/F1, F1/G1, G1/H1, I1/J1, J1/K1, K1/L1, L1/M1, M1/N1, N1/O1, O1/P1, 
				 A1/F1, B1/M1, C1/H1, D1/O1, E1/J1, G1/L1, I1/N1, K1/P1}
    	\draw (\from) -- (\to);

  	\node [black](A2) at (31,16) {};
  	\node (B2) at (33,15) {};
  	\node [black](C2) at (35,13) {};
  	\node (D2) at (36,11) {};
  	\node [black](E2) at (36,8) {};
  	\node (F2) at (35,6) {};
  	\node [black](G2) at (33,4) {};
  	\node (H2) at (31,3) {};
  	\node [black](I2) at (28,3) {};
  	\node (J2) at (26,4) {};
  	\node [black](K2) at (24,6) {};
  	\node (L2) at (23,8) {};
  	\node [black](M2) at (23,11) {};
  	\node (N2) at (24,13) {};
  	\node [black](O2) at (26,15) {};
  	\node (P2) at (28,16) {};

  \foreach \from/\to in {A2/B2, B2/C2, C2/D2, D2/E2, E2/F2, F2/G2, G2/H2, I2/J2, J2/K2, K2/L2, L2/M2, M2/N2, N2/O2, O2/P2, 
			 	A2/F2, B2/M2, C2/H2, D2/O2, E2/J2, G2/L2, I2/N2, K2/P2}
    	\draw (\from) -- (\to);

 	\draw (H1) to [bend right=30] (I2);

  \path[every node/.style={font=\sffamily}]
	(P1) edge [red, bend left=30] node [black] {} (A2);

\end{tikzpicture}
\caption{Horton fragment}
\label{Horton fragment}
\end{figure}

\end{minipage}

\begin{proof}
Suppose there is a Hamiltonian cycle that does not use edge $e$. Then there must be a pair of edges from $a, b, c$ in the Hamiltonian cycle. There are three possible pairs: 
\begin{itemize}
\item $a, b$. Then there would be a Hamiltonian path in the right Horton circle starting in white and ending in white, which is not possible because 16 vertices must be visited in alternating colors. 
\item $b,c$. Then there would be a Hamiltonian cycle in the left Horton circle using edge $e_2$ and not $e_1$.
\item $a, c$. Then there would be a Hamiltonian cycle in the right Horton circle using only edge $e_2$ and not $e_1$.
\end{itemize}

None of these cases is possible because of Lemma~\ref{Hortoncirle lemma1}, thus every Hamiltonian cycle must use edge $e$.  
\end{proof}

We have even proved that a Hamiltonian cycle must use all four edges connecting the two circles, but it is enough to know that edge $e$ must be used. 

\subsection{Horton-96}

In the Horton-96 graph three Horton fragments are connected as shown in  Figure~\ref{Horton-96 graph}.

\begin{thrm}
The Horton-96 graph is non-Hamiltonian. 

\end{thrm}

\begin{proof}
Suppose the red edge in Figure~\ref{Horton-96 graph} is not used. Then we can modify the bottom fragment such that we get the graph shown in Figure~\ref{Horton fragment}. Then we obtain a Hamiltonian cycle in this figure not using edge $e$, which is a contradiction to  Lemma~\ref{Hortonfragment lemma2}.

This is true for all the three fragments of the graph. So if there is a Hamiltonian cycle it must use all three edges incident to the red vertex which is not possible. 
\end{proof}


\begin{figure}[H]
\begin{center}
\begin{tikzpicture}
  [scale=.27,auto=left,every node/.style={circle, draw ,fill=white, scale=0.7}]

\draw [green] plot [smooth] coordinates { (45,6) (41,16) (29.5,25)}; 
\draw plot [smooth] coordinates { (5.3,21.3) (15.9,19.7) (29.5,25)};
\draw plot [smooth] coordinates { (38.4,47.8) (31.7,39.4) (29.5,25)};


\draw plot [smooth] coordinates { (18,3) (23,1) (29.5,9.5) (36,17) (38,16)};
\draw plot [smooth] coordinates { (21,16) (23,17) (29.5,9.5) (36,1) (41,3)};

  	\node [black](aA1) at (21,16) {};
  	\node (aB1) at (23,15) {};
  	\node [black](aC1) at (25,13) {};
  	\node (aD1) at (26,11) {};
  	\node [black](aE1) at (26,8) {};
  	\node (aF1) at (25,6) {};
  	\node [black](aG1) at (23,4) {};
  	\node (aH1) at (21,3) {};
  	\node [black](aI1) at (18,3) {};
  	\node (aJ1) at (16,4) {};
  	\node [black](aK1) at (14,6) {};
  	\node (aL1) at (13,8) {};
  	\node [black](aM1) at (13,11) {};
  	\node (aN1) at (14,13) {};
  	\node [black](aO1) at (16,15) {};
  	\node (aP1) at (18,16) {};

  \foreach \from/\to in {aA1/aB1, aB1/aC1, aC1/aD1, aD1/aE1, aE1/aF1, aF1/aG1, aG1/aH1, aI1/aJ1, aJ1/aK1, aK1/aL1, aL1/aM1,   
  				aM1/aN1, aN1/aO1, aO1/aP1, 
				 aA1/aF1, aB1/aM1, aC1/aH1, aD1/aO1, aE1/aJ1, aG1/aL1, aI1/aN1, aK1/aP1}
    	\draw (\from) -- (\to);

  	\node (aB2) at (43,15) {};
  	\node [black](aC2) at (45,13) {};
  	\node (aD2) at (46,11) {};
  	\node [black](aE2) at (46,8) {};
  	\node (aF2) at (45,6) {};
  	\node [black](aG2) at (43,4) {};
  	\node (aH2) at (41,3) {};
  	\node [black](aI2) at (38,3) {};
  	\node (aJ2) at (36,4) {};
  	\node [black](aK2) at (34,6) {};
  	\node (aL2) at (33,8) {};
  	\node [black](aM2) at (33,11) {};
  	\node (aN2) at (34,13) {};
  	\node [black](aO2) at (36,15) {};
  	\node (aP2) at (38,16) {};

  \foreach \from/\to in {aB2/aC2, aC2/aD2, aD2/aE2, aE2/aF2, aF2/aG2, aG2/aH2, aI2/aJ2, aJ2/aK2, aK2/aL2, aL2/aM2,  
 				aM2/aN2, aN2/aO2, aO2/aP2, 
			 	aB2/aM2, aC2/aH2, aD2/aO2, aE2/aJ2, aG2/aL2, aI2/aN2, aK2/aP2}
    	\draw (\from) -- (\to);

 	\draw (aH1) to [bend right=30] (aI2);


\draw plot [smooth] coordinates { (15,43.5) (12.3,42.8) (16.2, 33) (19.4,23.4) (17.5,22.2)};
\draw plot [smooth] coordinates { (26.1,36.9) (26,34.7) (16.2,33) (5.8,31.7) (6.3,28.8)};

  	\node [black](bA1) at (26.1,36.9) {};
  	\node (bB1) at (24.3,35.6) {};
  	\node [black](bC1) at (21.5,34.9) {};
  	\node (bD1) at (19.3,35.1) {};
  	\node [black](bE1) at (16.7,36.6) {};
  	\node (bF1) at (15.6,38.5) {};
  	\node [black](bG1) at (14.8,41.2) {};
  	\node (bH1) at (15,43.5) {};
  	\node [black](bI1) at (16.5,46) {};
  	\node (bJ1) at (18.4,47.2) {};
  	\node [black](bK1) at (21.1,47.9) {};
  	\node (bL1) at (23.3,47.7) {};
  	\node [black](bM1) at (25.9,46.3) {};
  	\node (bN1) at (27.2,44.3) {};
  	\node [black](bO1) at (27.8,41.7) {};
  	\node (bP1) at (27.7,39.4) {};

  \foreach \from/\to in {bA1/bB1, bB1/bC1, bC1/bD1, bD1/bE1, bE1/bF1, bF1/bG1, bG1/bH1, bI1/bJ1, bJ1/bK1, bK1/bL1, bL1/bM1,   
  				bM1/bN1, bN1/bO1, bO1/bP1, 
				 bA1/bF1, bB1/bM1, bC1/bH1, bD1/bO1, bE1/bJ1, bG1/bL1, bI1/bN1, bK1/bP1}
    	\draw (\from) -- (\to);

  	\node (bB2) at (14.2,18.4) {};
  	\node [black](bC2) at (11.3,17.8) {};
  	\node (bD2) at (9.1,17.9) {};
  	\node [black](bE2) at (6.7,19.4) {};
  	\node (bF2) at (5.3,21.3) {};
  	\node [black](bG2) at (4.7,24) {};
  	\node (bH2) at (4.8,26.2) {};
  	\node [black](bI2) at (6.3,28.8) {};
  	\node (bJ2) at (8.2,30) {};
  	\node [black](bK2) at (11,30.8) {};
  	\node (bL2) at (13.2,30.6) {};
  	\node [black](bM2) at (15.7,29.1) {};
  	\node (bN2) at (16.9,27.1) {};
  	\node [black](bO2) at (17.6,24.5) {};
  	\node (bP2) at (17.5,22.2) {};

  \foreach \from/\to in {bB2/bC2, bC2/bD2, bD2/bE2, bE2/bF2, bF2/bG2, bG2/bH2, bI2/bJ2, bJ2/bK2, bK2/bL2, bL2/bM2,  
 				bM2/bN2, bN2/bO2, bO2/bP2, 
			 	bB2/bM2, bC2/bH2, bD2/bO2, bE2/bJ2, bG2/bL2, bI2/bN2, bK2/bP2}
    	\draw (\from) -- (\to);

 	\draw (bH1) to [bend right=30] (bI2);


\draw plot [smooth] coordinates { (41.6,22) (39.8,23.2) (43.1, 32.7) (47.2,42.4) (42.9,45.8)};
\draw plot [smooth] coordinates { (54.4,25.8) (53.6,31.1) (43.1, 32.7) (33.3,34.6) (33.2,36.8)};

  	\node [black](cA1) at (41.6,22) {};
  	\node (cB1) at (41.6,24.2) {};
  	\node [black](cC1) at (42.2,26.9) {};
  	\node (cD1) at (43.5,28.8) {};
  	\node [black](cE1) at (46.1,30.3) {};
  	\node (cF1) at (48.2,30.4) {};
  	\node [black](cG1) at (51,29.8) {};
  	\node (cH1) at (52.9,28.5) {};
  	\node [black](cI1) at (54.4,25.8) {};
  	\node (cJ1) at (54.5,23.7) {};
  	\node [black](cK1) at (53.8,20.9) {};
  	\node (cL1) at (52.6,19) {};
  	\node [black](cM1) at (49.9,17.5) {};
  	\node (cN1) at (47.7,17.4) {};
  	\node [black](cO1) at (45,18.2) {};
  	\node (cP1) at (43.2,19.4) {};

  \foreach \from/\to in {cA1/cB1, cB1/cC1, cC1/cD1, cD1/cE1, cE1/cF1, cF1/cG1, cG1/cH1, cI1/cJ1, cJ1/cK1, cK1/cL1, cL1/cM1,   
  				cM1/cN1, cN1/cO1, cO1/cP1, 
				 cA1/cF1, cB1/cM1, cC1/cH1, cD1/cO1, cE1/cJ1, cG1/cL1, cI1/cN1, cK1/cP1}
    	\draw (\from) -- (\to);

  	\node (cB2) at (31.6,41.7) {};
  	\node [black](cC2) at (32.2,44.2) {};
  	\node (cD2) at (33.6,46.2) {};
  	\node [black](cE2) at (36.1,47.7) {};
  	\node (cF2) at (38.4,47.8) {};
  	\node [black](cG2) at (41.1,47) {};
  	\node (cH2) at (42.9,45.8) {};
  	\node [black](cI2) at (44.4,43.2) {};
  	\node (cJ2) at (44.5,40.9) {};
  	\node [black](cK2) at (43.8,38.2) {};
  	\node (cL2) at (42.6,36.3) {};
  	\node [black](cM2) at (40,34.9) {};
  	\node (cN2) at (37.8,34.8) {};
  	\node [black](cO2) at (35.1,35.6) {};
  	\node (cP2) at (33.2,36.8) {};

  \foreach \from/\to in {cB2/cC2, cC2/cD2, cD2/cE2, cE2/cF2, cF2/cG2, cG2/cH2, cI2/cJ2, cJ2/cK2, cK2/cL2, cL2/cM2,  
 				cM2/cN2, cN2/cO2, cO2/cP2, 
			 	cB2/cM2, cC2/cH2, cD2/cO2, cE2/cJ2, cG2/cL2, cI2/cN2, cK2/cP2}
    	\draw (\from) -- (\to);

 	\draw (cH1) to [bend right=30] (cI2);


  	\node [red](X) at (27,27) {};
  	\node [black](Y) at (29.5,25) {};
 	\node [black](Z) at (32,27) {};

\foreach \from/\to in {X/bP1, Z/cB2}
    	\draw (\from) -- (\to);

\draw [red] (X) to (aP1);
\draw [green] (Z) to (aB2);

 \draw (X) to [bend right=45] (cP1); 
 \draw (Z) to [bend left=45] (bB2); 

\end{tikzpicture}
\caption{Horton-96 graph}
\label{Horton-96 graph}
\end{center}
\end{figure}

\subsection{Horton-92}

The Horton fragment in  Figure~\ref{Horton fragment part} can be represented by a symbol, shown in  Figure~\ref{Horton fragment symbol}. The upper edge corresponds to the red edge in Figure~\ref{Horton-96 graph}, the other two are the ones indicated in green. Note that the three other lines do not correspond to edges. 

\begin{figure}[H]

\centering
\begin{minipage}{.45\textwidth}


\begin{center}
\begin{tikzpicture}
  [scale=.15,auto=left,every node/.style={circle, draw, fill=white, scale=0.5}]

\draw [green] plot [smooth] coordinates { (45,6) (41,16) (29.5,25)}; 
     
\draw plot [smooth] coordinates { (18,3) (23,1) (29.5,9.5) (36,17) (38,16)};
\draw plot [smooth] coordinates { (21,16) (23,17) (29.5,9.5) (36,1) (41,3)};

  	\node [black](aA1) at (21,16) {};
  	\node (aB1) at (23,15) {};
  	\node [black](aC1) at (25,13) {};
  	\node (aD1) at (26,11) {};
  	\node [black](aE1) at (26,8) {};
  	\node (aF1) at (25,6) {};
  	\node [black](aG1) at (23,4) {};
  	\node (aH1) at (21,3) {};
  	\node [black](aI1) at (18,3) {};
  	\node (aJ1) at (16,4) {};
  	\node [black](aK1) at (14,6) {};
  	\node (aL1) at (13,8) {};
  	\node [black](aM1) at (13,11) {};
  	\node (aN1) at (14,13) {};
  	\node [black](aO1) at (16,15) {};
  	\node (aP1) at (18,16) {};

  \foreach \from/\to in {aA1/aB1, aB1/aC1, aC1/aD1, aD1/aE1, aE1/aF1, aF1/aG1, aG1/aH1, aI1/aJ1, aJ1/aK1, aK1/aL1, aL1/aM1,   
  				aM1/aN1, aN1/aO1, aO1/aP1, 
				 aA1/aF1, aB1/aM1, aC1/aH1, aD1/aO1, aE1/aJ1, aG1/aL1, aI1/aN1, aK1/aP1}
    	\draw (\from) -- (\to);

  	\node (aB2) at (43,15) {};
  	\node [black](aC2) at (45,13) {};
  	\node (aD2) at (46,11) {};
  	\node [black](aE2) at (46,8) {};
  	\node (aF2) at (45,6) {};
  	\node [black](aG2) at (43,4) {};
  	\node (aH2) at (41,3) {};
  	\node [black](aI2) at (38,3) {};
  	\node (aJ2) at (36,4) {};
  	\node [black](aK2) at (34,6) {};
  	\node (aL2) at (33,8) {};
  	\node [black](aM2) at (33,11) {};
  	\node (aN2) at (34,13) {};
  	\node [black](aO2) at (36,15) {};
  	\node (aP2) at (38,16) {};

  \foreach \from/\to in {aB2/aC2, aC2/aD2, aD2/aE2, aE2/aF2, aF2/aG2, aG2/aH2, aI2/aJ2, aJ2/aK2, aK2/aL2, aL2/aM2,  
 				aM2/aN2, aN2/aO2, aO2/aP2, 
			 	aB2/aM2, aC2/aH2, aD2/aO2, aE2/aJ2, aG2/aL2, aI2/aN2, aK2/aP2}
    	\draw (\from) -- (\to);

 	\draw (aH1) to [bend right=30] (aI2);

  	\node [white](X) at (27,27) {};
  	\node [white](Y) at (29.5,25) {};
 	\node [white](Z) at (32,27) {};

\draw [red] (X) to (aP1);
\draw [green] (Z) to (aB2);

\end{tikzpicture}
\caption{Horton fragment part}
\label{Horton fragment part}
\end{center}

\end{minipage}
\begin{minipage}{.45\textwidth}


\begin{center}
\begin{tikzpicture}
  [scale=1.08,auto=left,every node/.style={circle, draw, fill=white,  scale=.7}]

\draw [red] plot coordinates { (0,2) (0,1)};
\draw [green] plot coordinates { (-2,-2) (-1,-1)};
\draw [green] plot coordinates { (2,-2) (1,-1)};

  	\node (P1) at (0,1) {};
  	\node (F2) at (-1,-1) {};
  	\node (B2) at (1,-1) {};

 	\draw (P1) to [bend right=50] (F2); 
 	\draw (F2) to [bend right=30] (B2); 
 	\draw (B2) to [bend right=50] (P1);

\end{tikzpicture}
\caption{Horton fragment symbol}
\label{Horton fragment symbol}
\end{center}

\end{minipage}
\end{figure}

We can now draw the Horton graph more clearly, as is shown in Figure~\ref{Horton-96 graph symbol}. It is easy to see that this is a counterexample to Tutte's conjecture, because we know that a Hamiltonian cycle through the graph must use the top edge of each fragment. But then all three edges of the top vertex must be used, which is not possible for a Hamiltonian cycle.

Now it is also easy to show the smaller counterexample Horton found and to imagine how he could have found it. The Horton-96 graph uses three fragments and three extra vertices. In the Horton-92 graph these three extra vertices are not necessary. One Horton fragment is colored the other way around. Then we can connect the vertices in the way shown in  Figure~\ref{Horton-92 graph symbol}. The cross means that we delete this white vertex. In the fragment this vertex is connected to two black vertices, which are the two black vertices shown in the bottom-right fragment. Note that the third edge of the deleted vertex is the one leaving the fragment. This fragment now has 30 vertices, the others have 31, which gives 92 vertices in total. 

\begin{figure}[H]

\centering
\begin{minipage}{.45\textwidth}


\begin{center}
\begin{tikzpicture}
  [scale=0.5,auto=left,every node/.style={circle, draw, fill=white,  scale=.7}]

  	\node (aP1) at (2,6) {};
  	\node (aF2) at (1,4) {};
  	\node (aB2) at (3,4) {};

  	\node (bP1) at (7,6) {};
  	\node (bF2) at (6,4) {};
  	\node (bB2) at (8,4) {};

  	\node (cP1) at (12,6) {};
  	\node (cF2) at (11,4) {};
  	\node (cB2) at (13,4) {};

	\node (X) [black] at (7,9) {};
  	\node (Y) [black] at (4.5,1) {};
  	\node (Z) [black] at (9.5,1) {};

 	\draw (aP1) to [bend right=50] (aF2); 
 	\draw (aF2) to [bend right=30] (aB2); 
 	\draw (aB2) to [bend right=50] (aP1); 

	\draw (bP1) to [bend right=50] (bF2); 
 	\draw (bF2) to [bend right=30] (bB2); 
 	\draw (bB2) to [bend right=50] (bP1); 

	\draw (cP1) to [bend right=50] (cF2); 
 	\draw (cF2) to [bend right=30] (cB2); 
 	\draw (cB2) to [bend right=50] (cP1); 

 	\draw (aP1) to (X); 
 	\draw (bP1) to (X); 
 	\draw (cP1) to (X); 

 	\draw (aF2) to (Y); 
 	\draw (bF2) to (Y);
  	\draw (cF2) to (Y); 

 	\draw (aB2) to (Z); 
 	\draw (bB2) to (Z);
  	\draw (cB2) to (Z); 

\end{tikzpicture}
\caption{Horton-96 graph symbol}
\label{Horton-96 graph symbol}
\end{center}

\end{minipage}
\begin{minipage}{.45\textwidth}


\begin{center}
\begin{tikzpicture}
  [scale=0.62,auto=left,every node/.style={circle, draw, fill=white,  scale=.7}]

\draw plot [smooth] coordinates { (4,6) (4.5,6.5) (7,4) (8,1) (7,1)};
\draw plot [smooth] coordinates { (3,1) (3.3,1) (4,1.5) (4.6,2) (4.9,2)};
\draw plot [smooth] coordinates { (1,1)(0.7,0.5) (3.5,0) (6,0.2) (6,0.7)};

  	\node (aP1) [black] at (4,6) {};
  	\node (aF2) [black] at (3,4) {};
  	\node (aB2) [black] at (5,4) {};

  	\node (bP1) at (2,3) {};
  	\node (bF2) at (1,1) {};
  	\node (bB2) at (3,1) {};

  	\node (cP1) at (6,3) {};
  	\node (cF2) at (5,1) {};
  	\node (cB2) at (7,1) {};

  	\node (cQ) [black] at (4.9,2) {};
  	\node (cR) [black] at (6,0.7) {};

 	\draw (aP1) to [bend right=50] (aF2); 
 	\draw (aF2) to [bend right=30] (aB2); 
 	\draw (aB2) to [bend right=50] (aP1); 

	\draw (bP1) to [bend right=50] (bF2); 
 	\draw (bF2) to [bend right=30] (bB2); 
 	\draw (bB2) to [bend right=50] (bP1); 

	\draw (cB2) to [bend right=50] (cP1); 

 	\draw (cQ) to [bend right=15] (cF2); 
	\draw (cF2) to [bend right=15] (cR); 

 	\draw (cQ) to [bend left=25] (cP1); 
	\draw (cB2) to [bend left=15] (cR); 

 	\draw (aF2) to (bP1); 
 	\draw (aB2) to (cP1); 

\draw plot [smooth] coordinates {(4.5,1.5)(5.5,0.5)};
\draw plot [smooth] coordinates {(4.5,0.5)(5.5,1.5)};

\end{tikzpicture}
\caption{Horton-92 graph symbol}
\label{Horton-92 graph symbol}
\end{center}

\end{minipage}
\end{figure}

As Horton mentioned in \cite{Horton1982} the advantage of the construction of the Horton-92 graph rather than the construction used in the Horton-96 graph, is that if the Horton fragment would be replaced by a planar fragment, then the whole graph would be planar. Then it would be a counterexample to Barnette's conjecture. 

Note that the Horton fragment is not planar, because $K_{3,3}$ is a minor of the graph. See Figure~\ref{K3,3}. 


\begin{figure}
\begin{center}
\begin{tikzpicture}
  [scale=.2,auto=left,every node/.style={circle, draw=black, fill=white, scale=0.7}]

\draw [black, dashed] plot [smooth] coordinates { (45,6) (41,16) (29.5,25)}; 
     
\draw [red, very thick] plot [smooth] coordinates { (18,3) (23,1) (29.5,9.5) (36,17) (38,16)};
\draw [red, very thick] plot [smooth] coordinates { (21,16) (23,17) (29.5,9.5) (36,1) (41,3)};

  	\node [black](aA1) at (21,16) {};
  	\node (aB1) at (23,15) {};
  	\node [black](aC1) at (25,13) {};
  	\node (aD1) at (26,11) {};
  	\node [black](aE1) at (26,8) {};
  	\node (aF1) at (25,6) {};
  	\node [black](aG1) at (23,4) {};
  	\node [draw=black, fill=green] (aH1) at (21,3) {};
  	\node [black](aI1) at (18,3) {};
  	\node (aJ1) at (16,4) {};
  	\node [black](aK1) at (14,6) {};
  	\node (aL1) at (13,8) {};
  	\node [black](aM1) at (13,11) {};
  	\node (aN1) at (14,13) {};
  	\node [black](aO1) at (16,15) {};
  	\node (aP1) at (18,16) {};

  \foreach \from/\to in {aA1/aB1, aB1/aC1, aC1/aD1, aD1/aE1, aE1/aF1, aF1/aG1, aG1/aH1, aI1/aJ1, aJ1/aK1, aK1/aL1, aL1/aM1,   
  				aM1/aN1, aN1/aO1, aO1/aP1, 
				 aA1/aF1, aB1/aM1, aC1/aH1, aD1/aO1, aE1/aJ1, aG1/aL1, aI1/aN1, aK1/aP1}
    	\draw [dashed] (\from) -- (\to);

  	\node (aB2) at (43,15) {};
  	\node [black](aC2) at (45,13) {};
  	\node (aD2) at (46,11) {};
  	\node [draw=black, fill=blue](aE2) at (46,8) {};
  	\node (aF2) at (45,6) {};
  	\node [black](aG2) at (43,4) {};
  	\node (aH2) at (41,3) {};
  	\node [draw=black, fill=blue](aI2) at (38,3) {};
  	\node [draw=black, fill=green](aJ2) at (36,4) {};
  	\node [draw=black, fill=blue](aK2) at (34,6) {};
  	\node (aL2) at (33,8) {};
  	\node [black](aM2) at (33,11) {};
  	\node [draw=black, fill=green](aN2) at (34,13) {};
  	\node [black](aO2) at (36,15) {};
  	\node (aP2) at (38,16) {};

  \foreach \from/\to in {aB2/aC2, aC2/aD2, aD2/aE2, aE2/aF2, aF2/aG2, aG2/aH2, aI2/aJ2, aJ2/aK2, aK2/aL2, aL2/aM2,  
 				aM2/aN2, aN2/aO2, aO2/aP2, 
			 	aB2/aM2, aC2/aH2, aD2/aO2, aE2/aJ2, aG2/aL2, aI2/aN2, aK2/aP2}
    	\draw [dashed] (\from) -- (\to);

 	\draw [red, very thick] (aH1) to [bend right=30] (aI2); 

  \foreach \from/\to in {aA1/aB1, aB1/aC1, aE1/aF1, aF1/aG1, aG1/aH1, aI1/aJ1, aC1/aH1, aE1/aJ1}
    	\draw  [red, very thick] (\from) -- (\to);

  \foreach \from/\to in {aD2/aE2, aE2/aF2, aF2/aG2, aG2/aH2, aI2/aJ2, aJ2/aK2, aK2/aL2, aL2/aM2,  
 				aM2/aN2, aN2/aO2, 
			 	aD2/aO2, aI2/aN2, aK2/aP2, aE2/aJ2}
    	\draw  [red, very thick] (\from) -- (\to);

  	\node [white](X) at (27,27) {};
  	\node [white](Y) at (29.5,25) {};
 	\node [white](Z) at (32,27) {};

\draw [black, dashed] (X) to (aP1);
\draw [black, dashed] (Z) to (aB2);

\end{tikzpicture}
\caption{$K_{3,3}$}
\label{K3,3}
\end{center}
\end{figure}

\section{Ellingham Fragment}

The Ellingham fragment shown in Figure~\ref{Ellingham fragment} is used in several counterexamples.

\begin{lemma}
No Hamiltonian cycle in an Ellingham fragment contains both red edges.
\label{Ellingham}
\end{lemma}

\noindent \begin{minipage}{0.6\textwidth}

{\bf Sketch of the proof:}
Suppose there is a Hamiltonian cycle in the Ellingham fragment that uses both red edges. Consider the four endpoints of these red edges. For each of these vertices there are two possibilities: go up or go down. Because of symmetry this gives 10 different situations of subpaths containing the red edges and one other edge for each of the four vertices on the red edges. In each case we can use the same arguments as in Horton's graph to see that there would be two cycles, or that all three edges of a vertex are contained in the Hamiltonian cycle. This gives a contradiction to the assumption, so there can not be a Hamiltonion cycle using both red edges.

\end{minipage}
\begin{minipage}{0.4\textwidth}


\begin{figure}[H]
\centering

\begin{tikzpicture}
  [scale=.27,auto=left,every node/.style={circle, draw, fill=white, scale=.7}]
  	\node [black](A) at (5,14) {};
  	\node (B) at (12,14) {};
  	\node [black](C) at (15,8) {};
  	\node (D) at (12,2) {};
  	\node [black](E) at (5,2) {};
  	\node (F) at (2,8) {};
  	\node (G) at (6,12) {};
  	\node [black](H) at (11,12) {};
  	\node (I) at (13,8) {};
  	\node [black](J) at (11,4) {};
  	\node (K) at (6,4) {};
  	\node [black](L) at (4,8) {};
  	\node [black](M) at (7,10) {};
  	\node (N) at (10,10) {};
  	\node [black](O) at (11,8) {};
  	\node (P) at (10,6) {};
  	\node [black](Q) at (7,6) {};
  	\node (R) at (6,8) {};

  \foreach \from/\to in {B/C, C/D, D/E, E/F, F/A, 
				N/O, O/P, P/Q, Q/R, R/M,
				A/G, G/M, B/H, H/N, C/I, I/O, D/J, J/P, E/K, K/Q, F/L, L/R,
				K/J, A/N, B/M,
				G/H}
    	\draw (\from) -- (\to);

 	\draw (L) to [bend right=100] (I); 

	\draw [red] (G) to (H);
	\draw [red] (K) to (J);

\end{tikzpicture}
\caption{Ellingham fragment}
\label{Ellingham fragment}

\end{figure}
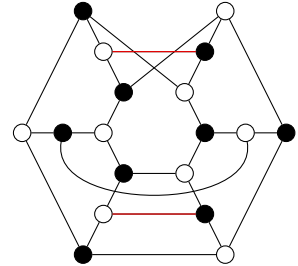

\end{minipage}\\

\begin{figure}[H]
\begin{center}
\begin{tikzpicture}
  [scale=.2,auto=left,every node/.style={circle, draw, fill=white, scale=.5}]

\draw plot [smooth] coordinates { (19,14) (24,17) (31,19) (38,17) (44,8) (41,0)(38,0)(36,4)};
\draw plot [smooth] coordinates { (24,2)   (19,-1)  (12,-3)  (5, -1)  (-1,8)  (2,16) (5,16)(7,12)};

\draw plot [smooth] coordinates { (10,4)   (5,0)       (1,0)       (-4,8)   (3,19)    (15, 19) (24, 14)};
\draw plot [smooth] coordinates { (33,12) (38, 16)  (42, 16)   (47, 8)(41,-3)    (31, -3)  (24,2)};


  	\node [black](A1) at (5,14) {};
  	\node (B1) at (12,14) {};
  	\node [black](C1) at (15,8) {};
  	\node (D1) at (12,2) {};
  	\node [black](E1) at (5,2) {};
  	\node (F1) at (2,8) {};
  	\node (G1) at (6,12) {};
  	\node [black](H1) at (11,12) {};
  	\node (I1) at (13,8) {};
  	\node [black](J1) at (11,4) {};
  	\node (K1) at (6,4) {};
  	\node [black](L1) at (4,8) {};
  	\node [black](M1) at (7,10) {};
  	\node (N1) at (10,10) {};
  	\node [black](O1) at (11,8) {};
  	\node (P1) at (10,6) {};
  	\node [black](Q1) at (7,6) {};
  	\node (R1) at (6,8) {};

  \foreach \from/\to in {B1/C1, C1/D1, D1/E1, E1/F1, F1/A1, 
				N1/O1, O1/P1, P1/Q1, Q1/R1, R1/M1,
				A1/G1, G1/M1, B1/H1, H1/N1, C1/I1, I1/O1, D1/J1, J1/P1, E1/K1, K1/Q1, F1/L1, L1/R1,
				K1/J1, A1/N1, B1/M1,
				G1/H1, K1/J1}
    	\draw (\from) -- (\to);

 	\draw (L1) to [bend right=100] (I1);

 	\node [black](S1) at (7,12) {};
  	\node (T1) at (10,12) {};
  	\node [black](U1) at (7,4) {};
  	\node (V1) at (10,4) {};


\node [black](A2) at (31,14) {};
  	\node (B2) at (38,14) {};
  	\node [black](C2) at (41,8) {};
  	\node (D2) at (38,2) {};
  	\node [black](E2) at (31,2) {};
  	\node (F2) at (28,8) {};
  	\node (G2) at (32,12) {};
  	\node [black](H2) at (37,12) {};
  	\node (I2) at (39,8) {};
  	\node [black](J2) at (37,4) {};
  	\node (K2) at (32,4) {};
  	\node [black](L2) at (30,8) {};
  	\node [black](M2) at (33,10) {};
  	\node (N2) at (36,10) {};
  	\node [black](O2) at (37,8) {};
  	\node (P2) at (36,6) {};
  	\node [black](Q2) at (33,6) {};
  	\node (R2) at (32,8) {};

  \foreach \from/\to in {B2/C2, C2/D2, D2/E2, E2/F2, F2/A2, 
				N2/O2, O2/P2, P2/Q2, Q2/R2, R2/M2,
				A2/G2, G2/M2, B2/H2, H2/N2, C2/I2, I2/O2, D2/J2, J2/P2, E2/K2, K2/Q2, F2/L2, L2/R2,
				K2/J2, A2/N2, B2/M2,
				G2/H2, K2/J2}
    	\draw (\from) -- (\to);

 	\draw (L2) to [bend right=100] (I2); 

 	\node [black](S2) at (33,12) {};
  	\node (T2) at (36,12) {};
  	\node [black](U2) at (33,4) {};
  	\node (V2) at (36,4) {};

  	\node [black](P) at (19,14) {};
  	\node [black](Q) at (24,14) {};
	\node (X) at (21.5,10) {};
  	\node [black](Y) at (21.5,6) {};
  	\node (R) at (19,2) {};
  	\node (S) at (24,2) {};

 \foreach \from/\to in {P/X, Q/X, X/Y, Y/R, Y/S}
    	\draw (\from) -- (\to);

 	\draw (U1) to [bend right=80] (R); 
 	\draw (R) to [bend right=80] (U2); 

 	\draw (T1) to [bend left=80] (P); 
 	\draw (Q) to [bend left=80] (T2);

\end{tikzpicture}
\caption{Georges's graph}
\label{Georges Graph}
\end{center}
\end{figure}
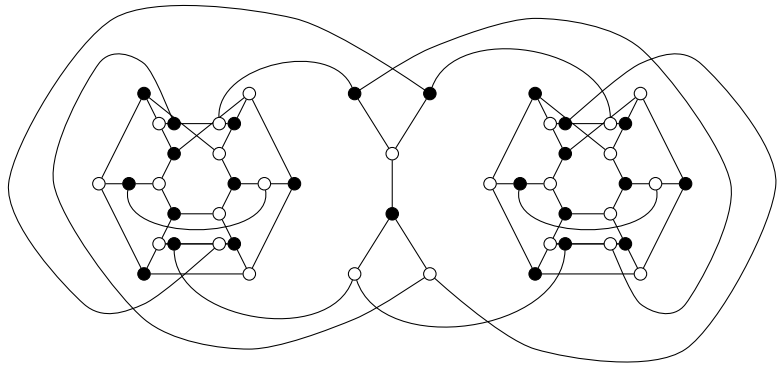

In 1989 John P. Georges found Georges's graph \cite{Georges1989}, which is shown in  Figure~\ref{Georges Graph}. It is clear that the graph is bipartite and cubic. We used the algorithm in Chapter~\ref{Programsinc} to check that the graph is 3-connected. 


The graph consists of two Ellingham fragments $B_1$ and $B_2$, with four extra vertices placed on the red edges of Lemma~\ref{Ellingham}. These extra vertices are connected to a tree $T$ of 6 vertices in the middle of the graph. 

To make life easier we only consider the edges connecting the Ellingham fragments to the tree. So we only take the red vertices of Figure~\ref{Georges Graph part colored} and depict Georges's graph as in Figure~\ref{Georges Graph part}.

\begin{figure}[H]

\centering
\begin{minipage}{.45\textwidth}


\begin{center}
\begin{tikzpicture}
  [scale=.14,auto=left,every node/.style={circle, draw, fill=white, scale=.4}]

\draw plot [smooth] coordinates { (19,14) (24,17) (31,19) (38,17) (44,8) (41,0)(38,0)(36,4)};
\draw plot [smooth] coordinates { (24,2)   (19,-1)  (12,-3)  (5, -1)  (-1,8)  (2,16) (5,16)(7,12)};

\draw plot [smooth] coordinates { (10,4)   (5,0)       (1,0)       (-4,8)   (3,19)    (15, 19) (24, 14)};
\draw plot [smooth] coordinates { (33,12) (38, 16)  (42, 16)   (47, 8)(41,-3)    (31, -3)  (24,2)};


  	\node [black](A1) at (5,14) {};
  	\node (B1) at (12,14) {};
  	\node [black](C1) at (15,8) {};
  	\node (D1) at (12,2) {};
  	\node [black](E1) at (5,2) {};
  	\node (F1) at (2,8) {};
  	\node (G1) at (6,12) {};
  	\node [black](H1) at (11,12) {};
  	\node (I1) at (13,8) {};
  	\node [black](J1) at (11,4) {};
  	\node (K1) at (6,4) {};
  	\node [black](L1) at (4,8) {};
  	\node [black](M1) at (7,10) {};
  	\node (N1) at (10,10) {};
  	\node [black](O1) at (11,8) {};
  	\node (P1) at (10,6) {};
  	\node [black](Q1) at (7,6) {};
  	\node (R1) at (6,8) {};

  \foreach \from/\to in {B1/C1, C1/D1, D1/E1, E1/F1, F1/A1, 
				N1/O1, O1/P1, P1/Q1, Q1/R1, R1/M1,
				A1/G1, G1/M1, B1/H1, H1/N1, C1/I1, I1/O1, D1/J1, J1/P1, E1/K1, K1/Q1, F1/L1, L1/R1,
				K1/J1, A1/N1, B1/M1,
				G1/H1, K1/J1}
    	\draw (\from) -- (\to);

 	\draw (L1) to [bend right=100] (I1);

 	\node [red](S1) at (7,12) {};
  	\node [red](T1) at (10,12) {};
  	\node [red](U1) at (7,4) {};
  	\node [red](V1) at (10,4) {};


\node [black](A2) at (31,14) {};
  	\node (B2) at (38,14) {};
  	\node [black](C2) at (41,8) {};
  	\node (D2) at (38,2) {};
  	\node [black](E2) at (31,2) {};
  	\node (F2) at (28,8) {};
  	\node (G2) at (32,12) {};
  	\node [black](H2) at (37,12) {};
  	\node (I2) at (39,8) {};
  	\node [black](J2) at (37,4) {};
  	\node (K2) at (32,4) {};
  	\node [black](L2) at (30,8) {};
  	\node [black](M2) at (33,10) {};
  	\node (N2) at (36,10) {};
  	\node [black](O2) at (37,8) {};
  	\node (P2) at (36,6) {};
  	\node [black](Q2) at (33,6) {};
  	\node (R2) at (32,8) {};

  \foreach \from/\to in {B2/C2, C2/D2, D2/E2, E2/F2, F2/A2, 
				N2/O2, O2/P2, P2/Q2, Q2/R2, R2/M2,
				A2/G2, G2/M2, B2/H2, H2/N2, C2/I2, I2/O2, D2/J2, J2/P2, E2/K2, K2/Q2, F2/L2, L2/R2,
				K2/J2, A2/N2, B2/M2,
				G2/H2, K2/J2}
    	\draw (\from) -- (\to);

 	\draw (L2) to [bend right=100] (I2); 

 	\node [red](S2) at (33,12) {};
  	\node [red] (T2) at (36,12) {};
  	\node [red](U2) at (33,4) {};
  	\node [red](V2) at (36,4) {};

  	\node [red](P) at (19,14) {};
  	\node [red](Q) at (24,14) {};
	\node [red] (X) at (21.5,10) {};
  	\node [red](Y) at (21.5,6) {};
  	\node [red](R) at (19,2) {};
  	\node [red](S) at (24,2) {};

 \foreach \from/\to in {P/X, Q/X, X/Y, Y/R, Y/S}
    	\draw (\from) -- (\to);

 	\draw (U1) to [bend right=80] (R); 
 	\draw (R) to [bend right=80] (U2); 

 	\draw (T1) to [bend left=80] (P); 
 	\draw (Q) to [bend left=80] (T2);

\end{tikzpicture}
\caption{}
\label{Georges Graph part colored}
\end{center}

\end{minipage}
\begin{minipage}{.45\textwidth}


\begin{center}
\begin{tikzpicture}
  [scale=.14,auto=left,every node/.style={circle, draw, fill=white, scale=.5}]

\draw plot [smooth] coordinates { (19,14) (24,17) (31,19) (38,17) (44,8) (41,0)(38,0)(36,4)};
\draw plot [smooth] coordinates { (24,2)   (19,-1)  (12,-3)  (5, -1)  (-1,8)  (2,16) (5,16)(7,12)};

\draw plot [smooth] coordinates { (10,4)   (5,0)       (1,0)       (-4,8)   (3,19)    (15, 19) (24, 14)};
\draw plot [smooth] coordinates { (33,12) (38, 16)  (42, 16)   (47, 8)(41,-3)    (31, -3)  (24,2)};





 	\node [black, text=white](S1) at (7,12) {$E_1$};
  	\node [](T1) at (10,12) {$F_1$};
  	\node [black, text=white](U1) at (7,4) {$G_1$};
  	\node [](V1) at (10,4) {$H_1$};





 	\node [black, text=white](S2) at (33,12) {$E_2$};
  	\node [] (T2) at (36,12) {$F_2$};
  	\node [black, text=white](U2) at (33,4) {$G_2$};
  	\node [] (V2) at (36,4) {$H_2$};

  	\node [black, text=white](P) at (19,14) {P};
  	\node [black, text=white](Q) at (24,14) {Q};
	\node [](X) at (21.5,10) {X};
  	\node [black, ,text=white](Y) at (21.5,6) {Y};
  	\node [](R) at (19,2) {R};
  	\node [](S) at (24,2) {S};

 \foreach \from/\to in {P/X, Q/X, X/Y, Y/R, Y/S}
    	\draw (\from) -- (\to);

 	\draw (U1) to [bend right=80] (R); 
 	\draw (R) to [bend right=80] (U2); 

 	\draw (T1) to [bend left=80] (P); 
 	\draw (Q) to [bend left=80] (T2);

\end{tikzpicture}
\caption{}
\label{Georges Graph part}
\end{center}

\end{minipage}
\end{figure}

\begin{thrm}
There is no Hamiltonian cycle in Georges's graph.
\end{thrm}

\begin{proof}
Suppose there is a Hamiltonian cycle $H$ in the graph. We consider four cases:
\begin{itemize}
\item Case 1: $H$ uses 2 edges joining $B_1$ and $T$ and 2 edges joining $B_2$ and $T$. Then there must be a Hamiltonian path in both Ellingham fragments. Because the number of vertices of a fragment is even, this path must have a white and a black endpoint. It can not go from $E_i$ to $F_i$ or from $G_i$ to $H_i$ because of  Lemma~\ref{Ellingham}, so the endpoint of the Hamiltonian paths must be $E_i$ and $H_i$ or $F_i$ and $G_i$. 

Suppose we take $E_1, H_1, E_2$ and $H_2$ as endpoints, then edges $E_1S, H_1Q, E_2S, H_2P$ must be used. Because of the assumption that only two edges are used, the edges $F_1P, G_1R, F_2Q, G_2R$ can not be used. Now it is clear that vertex $R$ will not be on the cycle, which is not possible. Another combination of endpoints gives a similar conclusion.

\item Case 2: $H$ uses 2 edges joining $B_1$ and $T$ and 4 edges joining $B_2$ and $T$. We again suppose we take $E_1$ and $H_1$ as endpoints, the other case is similar. Note that in $B_2$ there must be a path $P_1$ from $E_2$ to $F_2$ and a path $P_2$ from $G_2$ to $H_2$, because otherwise the two paths together with edges $E_2F_2$ and $G_2H_2$ would form a Hamiltonian cycle which is a contradiction to Lemma~\ref{Ellingham}. But then $P_2$ together with edges $G_2R$, $RY$, $YX$, $XP$ and $PH_2$ would form a cycle that does not include all vertices.

\item Case 3: $H$ uses 4 edges joining $B_1$ and $T$ and 2 edges joining $B_2$ and $T$. This case is similar to Case 2. 

\item Case 4: $H$ uses 4 edges joining $B_1$ and $T$ and 4 edges joining $B_2$ and $T$. Then it is clear that $H$ does not contain $X$ and $Y$ which is a contradiction.

\end{itemize}
So there can not be a Hamiltonian cycle in Georges's graph.

\end{proof}

\chapter{Complexity}
\label{complexity}

We will now consider the complexity of some problems related to Barnette's conjecture. Two famous complexity classes are P and \NPP. We say that a decision problem is in P if there exists a polynomial time algorithm that solves all instances of the problem. We say that a decision problem is in NP if we can check for all instances in polynomial time whether a given solution is correct. A subset of NP is the class of \NPP-complete problems. A problem is \NPP-complete if it is in NP and if all problems in NP can be reduced to it in polynomial time. Problem $\Pi_{1}$ can be reduced to problem $\Pi_{2}$ if there is a polynomial time function from the instances of $\Pi_{1}$ to the instances of $\Pi_{2}$, such that yes-instances correspond.
To show that a problem $\Pi$ is \NPP-complete it suffices to prove that there is a reduction from a known \NPP-complete problem to $\Pi$. This implies that $\Pi$ is \NPP-complete, because all problems in NP are reducible to the known \NPP-complete problem and reductions are transitive, so all problems in NP are reducible to $\Pi$.

In this chapter we will show that the problems corresponding to Tait's and Tutte's conjecture are \NPP-complete. To show this, we will use the \NPP-complete problem $3$-SAT-CNF to reduce from. This problem is defined in Chapter 3 of \cite{Terwijn}. An instance of the problem $3$-SAT-CNF is a Boolean formula $F$ in conjunctive normal form with exactly $3$ literals per clause, so $F=C_1 \land C_2 \land \dots \land C_m$, where $C_j = l_{j1} \lor l_{j2} \lor l_{j3}$ for literals $l_{ji} \in \{x_1, \dots, x_n\} \cup \{\neg x_1, \dots, \neg x_n\}$. The problem is to decide whether or not there exists a truth assignment to the variables $\{x_1, \dots, x_n\}$ such that $F$ is true.\\

\section{Tait's conjecture and the necessity of bipartiteness}
First we consider a decision problem corresponding to Tait's conjecture. An instance of this problem is a cubic, $3$-connected, planar graph. The question is whether or not there exists a Hamiltonian cycle in this graph.

\begin{thrm}
\label{HAM*}
Finding a Hamiltonian cycle in a cubic, $3$-connected, planar graph is \NPP-complete.
\end{thrm}
\begin{proof}
We will prove this theorem by showing that finding a Hamiltonian cycle in a cubic, $3$-connected, planar graph is in NP and that there exists a reduction from $3$-SAT-CNF to this problem. This reduction can be found in \cite{GareyJohnsonTarjan}.

It is clear that this problem is in NP, because given a cycle in a cubic, $3$-connected, planar graph, we can check in polynomial time whether or not this cycle is Hamiltonian. This can be done by checking whether the cycle visits all vertices exactly once. In the remainder of the proof we will give a polynomial time transformation from an instance of $3$-SAT-CNF to a cubic, $3$-connected, planar graph and prove that yes-instances correspond under this transformation. An instance of $3$-SAT-CNF is given by a Boolean formula $F=C_1 \land C_2 \land \dots \land C_m$, where $C_j = l_{j1} \lor l_{j2} \lor l_{j3}$ for literals $l_{ji} \in \{x_1, \dots, x_n\} \cup \{\neg x_1, \dots, \neg x_n\}$. Yes-instances correspond when there exists a truth-assignment of $\{x_1, \dots, x_n\}$ that makes $F$ true if and only if there exists a Hamiltonian cycle in the corresponding cubic, $3$-connected, planar graph.

We start the construction of a cubic, $3$-connected, planar graph with an empty graph. For every variable in $F$ we add $4$ vertices which are connected as shown in Figure \ref{Variables in the graph}. For every literal one parallel edge in this construction represents {\it true} and one parallel edge represents {\it false}. For every clause in $F$ we add $6$ vertices which are connected as shown in Figure \ref{Clauses in the graph}. By connecting $v_{11}$ to $w_{11}$ and $v_{n4}$ to $w_{m6}$ we finish the basic part of the graph. To this basic part we add some subgraphs with special properties. The subgraphs we add are a XOR-module, a $2$-input-OR-module and a $3$-input-OR-module. The construction of these modules can be found in \cite{GareyJohnsonTarjan}, the properties of these modules will be explained here.

\noindent \begin{minipage}{0.6\textwidth}
\begin{minipage}{0.5\textwidth}
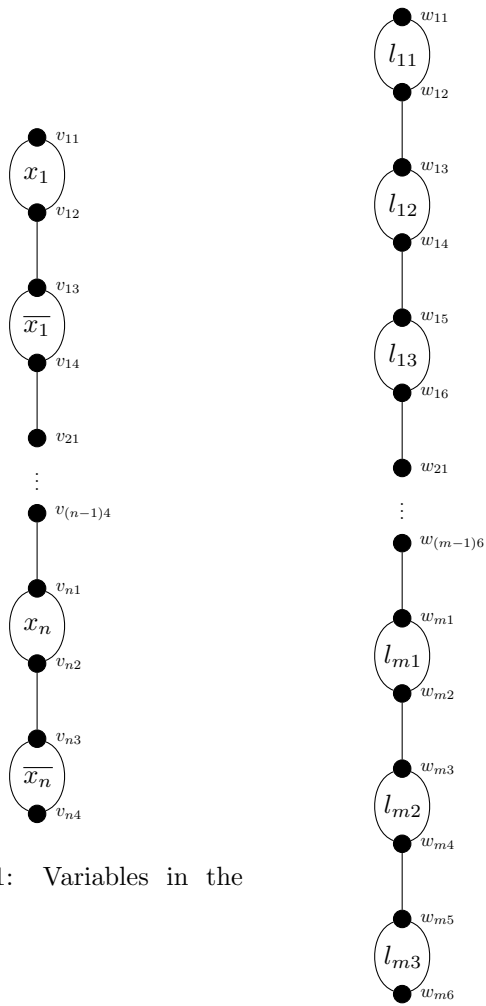
\begin{figure}[H]
\centering
\begin{tikzpicture}
  [scale=.5,auto=left,every node/.style={circle,draw, fill=black!100, scale=.7}]
  	\node [label=right:{$v_{11}$}] (A) at (0, 18) {};
  	\node [label=right:{$v_{12}$}] (B) at (0, 16) {};
  	\node [label=right:{$v_{13}$}] (C) at (0, 14) {};
  	\node [label=right:{$v_{14}$}] (D) at (0, 12) {};
  	\node [label=right:{$v_{21}$}] (E) at (0, 10) {};
	\node [label=right:{$v_{(n-1)4}$}] (F) at (0, 8) {};
  	\node [label=right:{$v_{n1}$}] (G) at (0, 6) {};
  	\node [label=right:{$v_{n2}$}] (H) at (0, 4) {};
  	\node [label=right:{$v_{n3}$}] (I) at (0, 2) {};
  	\node [label=right:{$v_{n4}$}] (J) at (0, 0) {};

  	\node  (X) at (0, 9) [draw=none, fill=white]{\vdots};

  \foreach \from/\to in {B/C,H/I, D/E, F/G}
    	\draw (\from) -- (\to);

  \path[every node/.style={}]
	(A) edge [bend right = 70] node [black] {\,$x_1$} (B)
	(C) edge [bend right = 70] node [black] {\,$\overline{x_1}$} (D)
	(G) edge [bend right = 70] node [black] {\,$x_n$} (H)
	(I) edge [bend right = 70] node [black] {\,$\overline{x_n}$} (J)

	(A) edge [bend left = 70] node [black] {} (B)
	(C) edge [bend left = 70] node [black] {} (D)
	(G) edge [bend left = 70] node [black] {} (H)
	(I) edge [bend left = 70] node [black] {} (J);
\end{tikzpicture}
\caption{Variables in the graph}
\label{Variables in the graph}
\end{figure}
\end{minipage}
\begin{minipage}{0.5\textwidth}
\vspace{4 mm}
\begin{figure}[H]
\centering
\begin{tikzpicture}
  [scale=.5,auto=left,every node/.style={circle,draw, fill=black!100, scale=.7}]
  	\node [label=right:{$w_{11}$}] (A) at (0, 26) {};
  	\node [label=right:{$w_{12}$}] (B) at (0, 24) {};
  	\node [label=right:{$w_{13}$}] (C) at (0, 22) {};
  	\node [label=right:{$w_{14}$}] (D) at (0, 20) {};
  	\node [label=right:{$w_{15}$}] (E) at (0, 18) {};
  	\node [label=right:{$w_{16}$}] (F) at (0, 16) {};
  	\node [label=right:{$w_{21}$}] (G) at (0, 14) {};
	\node [label=right:{$w_{(m-1)6}$}] (H) at (0, 12) {};
  	\node [label=right:{$w_{m1}$}] (I) at (0, 10) {};
  	\node [label=right:{$w_{m2}$}] (J) at (0, 8) {};
  	\node [label=right:{$w_{m3}$}] (K) at (0, 6) {};
  	\node [label=right:{$w_{m4}$}] (L) at (0, 4) {};
  	\node [label=right:{$w_{m5}$}] (M) at (0,2) {};
  	\node [label=right:{$w_{m6}$}] (N) at (0, 0) {};

  	\node  (X) at (0, 13) [draw=none, fill=white]{\vdots};

  \foreach \from/\to in {B/C, D/E, F/G, H/I, J/K, L/M}
    	\draw (\from) -- (\to);

  \path[every node/.style={}]
	(A) edge [bend right = 70] node [black] {\,$l_{11}$} (B)
	(C) edge [bend right = 70] node [black] {\,$l_{12}$} (D)
	(E) edge [bend right = 70] node [black] {\,$l_{13}$}(F)
	(I) edge [bend right = 70] node [black] {$l_{m1}$} (J)
	(K) edge [bend right = 70] node [black] {$l_{m2}$} (L)
	(M) edge [bend right = 70] node [black] {$l_{m3}$} (N)

	(A) edge [bend left = 70] node [black] {} (B)
	(C) edge [bend left = 70] node [black] {} (D)
	(E) edge [bend left = 70] node [black] {} (F)
	(I) edge [bend left = 70] node [black] {} (J)
	(K) edge [bend left = 70] node [black] {} (L)
	(M) edge [bend left = 70] node [black] {} (N);

\end{tikzpicture}
\caption{Clauses in the graph}
\label{Clauses in the graph}
\end{figure}
\end{minipage}
\end{minipage}
\begin{minipage}{0.4\textwidth}
\vspace{1 cm}
\begin{figure}[H]
\centering
\begin{tikzpicture}
  [scale=.5,auto=left,every node/.style={circle,draw, fill=black!100, scale=.7}]

  	\node (A1) at (0,6) {};
  	\node (B1) at (4,6) {};
  	\node (C1) at (5,6) {};
  	\node (D1) at (6,6) {};
 	\node (E1) at (7,6) {};
 	\node (F1) at (11,6) {};

  	\node (A2) at (0,0) {};
  	\node (B2) at (4,0) {};
  	\node (C2) at (5,0) {};
  	\node (D2) at (6,0) {};
 	\node (E2) at (7,0) {};
 	\node (F2) at (11,0) {};

  \foreach \from/\to in {A1/B1, B1/C1, C1/D1, D1/E1, E1/F1, A2/B2, B2/C2, C2/D2, D2/E2, E2/F2, B1/B2, C1/C2, D1/D2, E1/E2}
    	\draw (\from) -- (\to);

  	\node (X) [fill=white, text opacity=1] at (4,0)  {$\uparrow$};

  	\node (Y) [fill=white, text opacity=1] at (7,0)  {$\uparrow$};

\end{tikzpicture}
\caption{XOR-module}
\label{XOR-module}
\end{figure}
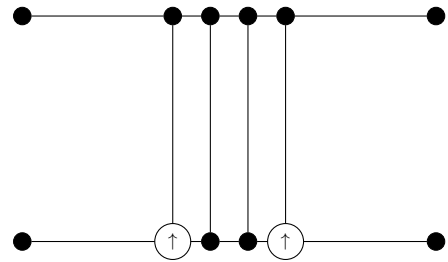
\vspace{2 cm}
\begin{figure}[H]
\centering
\begin{tikzpicture}
  [scale=.5,auto=left,every node/.style={circle,draw, fill=black!100, scale=.7}]

 \draw plot [smooth] coordinates {(2,2)(2,0)};
  	\node (A) at (0,2) {};
  	\node (B) at (4,2) {};
  	\node (C) at (0,0) {};
  	\node (D) at (4,0) {};

  	\node (X) [fill=white, text opacity=0] at (2,1)  {$X$};

  \foreach \from/\to in {A/B, C/D}
    	\draw (\from) -- (\to);

 \draw plot [smooth] coordinates {(1.65,1.4)(2.35,0.6)};
 \draw plot [smooth] coordinates {(1.65,0.6)(2.35,1.4)};

\end{tikzpicture}
\caption{Notation XOR-module}
\label{notation XOR-module}
\end{figure}
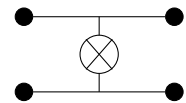
\end{minipage}

A XOR-module is a particular subgraph that can be implemented in a graph by replacing two disjoint edges of the graph by a XOR-module. This subgraph has the property that exactly one of those two edges has to be used in any Hamiltonian cycle. In this module two Tutte fragments (see Figure \ref{Tutte's fragment}) are used to obtain this property. The Tutte fragments are used because they have the property that a Hamiltonian cycle has to pass through one specified edge and through one of the other two edges. A XOR-module can be found in Figure \ref{XOR-module}. The Tutte fragments are indicated by the vertices that have an arrow inside. This arrow points to the top of a Tutte's fragment, where the edge is located that has to be used in any Hamiltonian path. The abbreviated notation for this subgraph is given in Figure \ref{notation XOR-module}. In Figure \ref{XOR-module} it can easily be checked that any Hamiltonian cycle has to use the upper edge or the lower edge in Figure \ref{notation XOR-module} and cannot use both. 

\begin{minipage}{0.6\textwidth}
A $2$-input-OR-module is a particular subgraph that can be implemented in a graph by replacing two disjoint edges by the two edges of a $2$-input-OR-module drawn in Figure \ref{2-input-or-module}. This subgraph has the property that at least one of those two edges has to be used in any Hamiltonian cycle. How this property is obtained can be found in \cite{GareyJohnsonTarjan}. The abbreviated notation for this subgraph is given in Figure \ref{2-input-or-module}.
\end{minipage}
\begin{minipage}{0.4\textwidth}
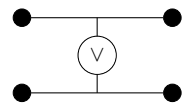
\begin{figure}[H]
\centering
\begin{tikzpicture}
  [scale=.5,auto=left,every node/.style={circle,draw, fill=black!100, scale=.7}]
 \draw plot [smooth] coordinates {(2,2)(2,0)};
  	\node (A) at (0,2) {};
  	\node (B) at (4,2) {};
  	\node (C) at (0,0) {};
  	\node (D) at (4,0) {};

  	\node (X) [fill=white, text opacity=0] at (2,1)  {$X$};
  	\node (X) [draw=none, fill=white, text opacity=1, font=\large] at (2,1)  {$\vee$};

  \foreach \from/\to in {A/B, C/D}
    	\draw (\from) -- (\to);
\end{tikzpicture}
\caption{2-input-OR-module}
\label{2-input-or-module}
\end{figure}
\end{minipage}

\begin{minipage}{0.6\textwidth}
A $3$-input-OR-module is a particular subgraph that can be implemented in a graph by replacing three disjoint edges by the three edges of a $3$-input-OR-module drawn in Figure \ref{3-input-or-module}. This subgraph has the property that at least one of those three edges has to be used in any Hamiltonian cycle. How this property is obtained can be found in \cite{GareyJohnsonTarjan}. The abbreviated notation for this subgraph is given in Figure \ref{3-input-or-module}.
\end{minipage}
\begin{minipage}{0.4\textwidth}
\begin{figure}[H]
\centering
\begin{tikzpicture}
  [scale=.5,auto=left,every node/.style={circle,draw, fill=black!100, scale=.7}]
 \draw plot [smooth] coordinates {(3,2)(3,1)};
 \draw plot [smooth] coordinates {(0,1)(6,1)};
  	\node (A) at (0,2) {};
  	\node (B) at (0,0) {};
  	\node (C) at (1,2) {};
  	\node (D) at (5,2) {};
	\node (E) at (6,2) {};
	\node (F) at (6,0) {};

  	\node (X) [fill=white, text opacity=0] at (3,1)  {$X$};
  	\node (X) [draw=none, fill=white, text opacity=1, font=\large] at (3,1)  {$\vee$};

  \foreach \from/\to in {A/B, C/D, E/F}
    	\draw (\from) -- (\to);
\end{tikzpicture}
\caption{3-input-OR-module}
\label{3-input-or-module}
\end{figure}
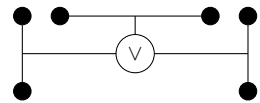
\end{minipage}

All of these subgraphs are planar. Note that Tutte's fragment, used in the XOR-module, is also planar. Now we will explain how to add the modules to the basic part of the graph. To understand this it might help to take a look at Figure \ref{Example}. The graph in this figure corresponds to the Boolean formula $F = (\neg x \lor y \lor \neg z) \land (x \lor \neg y \lor \neg z)$.

To eventually obtain the property that every variable is true or false and not both, we replace $v_{i1}v_{i2}$ and $v_{i3}v_{i4}$ for $i \in \{1, \dots, n\}$ by a XOR-module. To obtain the property that in every clause at least one literal is true, we replace $w_{j1}w_{j2}$, $w_{j3}w_{j4}$ and $w_{j5}w_{j6}$ for $j \in \{1, \dots, m\}$ by a $3$-input-OR-module. Furthermore we want to make sure that yes-instances correspond. To obtain this we add some more XOR-modules. If in $C_j$ literal $l_{jk} = x_i$, replace $w_{j(2k-1)}w_{j(2k)}$ and $v_{i1}v_{i2}$ by a XOR-module. If in $C_j$ literal $l_{jk}=\neg x_i$, replace $w_{j(2k-1)}w_{j(2k)}$ and $v_{i3}v_{i4}$ by a XOR-module. These XOR-modules use the parallel edges that are not yet used by other modules. Finally replace $v_{11}w_{11}$ and $v_{n4}w_{m6}$ by a $2$-input-OR-module to ensure $3$-connectedness. It is clear that in a Hamiltonian cycle both edges are used.

In the construction of this graph some XOR-modules will intersect each other, which makes the graph non-planar. This can be solved by making a new planar construction at every intersection. This planar construction can be found in \cite{GareyJohnsonTarjan}. This construction contains five XOR-modules and two Tutte fragments. For every intersection one such construction is added, which consists of 152 vertices.
This completes the construction of a cubic, $3$-connected, planar graph for a given Boolean formula $F=C_1 \land C_2 \land \dots \land C_m$. It is clear that this construction can be done in polynomial time and that the constructed graph is cubic, $3$-connected and planar.

\begin{claim} There exists a truth-assignment of $\{x_1, \dots, x_n\}$ that makes $F$ true if and only if there exists a Hamiltonian cycle in the constructed cubic, $3$-connected, planar graph.\end{claim}
\begin{proof}
First we will show that a truth assignment of $F$ corresponds to a Hamiltonian cycle in the constructed graph. If in the truth assignment of $F$ variable $x_i$ is true, use edge $v_{i1}v_{i2}$ that is connected by a XOR-module to the clauses. The XOR-module between $v_{i1}v_{i2}$ and $v_{i3}v_{i4}$ now ensures that the parallel edge $v_{i3}v_{i4}$ that is connected by a XOR-module to the clauses is not used. If in the truth assignment variable $x_i$ is not true, use the other of the two parallel edges. Furthermore, use all edges $v_{i2}v_{i3}$ and $v_{i4}v_{(i+1)1}$ to connect the vertices in the part of the graph corresponding to the variables. It is clear that the XOR-modules between $v_{i1}v_{i2}$ and $v_{i3}v_{i4}$ are satisfied for all  $i \in \{1, \dots, n\}$.

Use edge  $v_{n4}w_{m6}$ to go to the part of the graph corresponding to the clauses. The XOR-modules between the variables and the clauses leave only one option to pass through the part of the clauses. For example if $x_i$ is true, there is a XOR-module between $w_{j(2k-1)}w_{j(2k)}$ and $v_{i1}v_{i2}$. Because $x_i$ is true, in the part of the variables we have used the parallel edge $v_{i1}v_{i2}$ that is connected by a XOR-module to the clauses. So in the part corresponding to the clauses we can only use the parallel edge $w_{j(2k-1)}w_{j(2k)}$ that is connected by a $3$-input-OR-module to the rest of clause $C_j$. If we go through the part of the clauses like this, we make sure that every XOR-module between a variable and a clause is satisfied. When we reach $w_{11}$, we use edge $v_{11}w_{11}$ to complete the Hamiltonian cycle.

It is clear that the $2$-input-OR-module between $v_{11}w_{11}$ and $v_{n4}w_{m6}$ is satisfied, so we only have to check that for every clause the $3$-input-OR-module is satisfied. Because we started with a truth assignment of $F$, in every clause $C_j$ at least one literal $l_{jk}$ is true. If the literal $l_{jk}$ is true, the XOR-module between this literal and this clause is already satisfied in the part of the variables. The property of the XOR-module now ensures that for this literal in the part of the clauses the parallel edge connected to the $3$-input-OR-module has to be used. Now we can conclude that every $3$-input-OR-module is satisfied.

An example of a truth assignment and a corresponding Hamiltonian cycle is given in Figure \ref{Example}. This graph corresponds to the Boolean formula $F = (\neg x \lor y \lor \neg z) \land (x \lor \neg y \lor \neg z)$. The Hamiltonian path marked by the red edges corresponds to the truth assignment $(x, y, z) = (0, 0, 1)$.

\begin{minipage}{0.99\textwidth}
\begin{figure}[H]
\centering
\begin{tikzpicture}
  [scale=.5,auto=left,every node/.style={circle,draw, fill=black!100, scale=.7}]
	\node  (A) at (0, 18) {};
  	\node  (B) at (0, 16) {};
  	\node (C) at (0, 14) {};
  	\node  (D) at (0, 12) {};
  	\node  (E) at (0, 10) {};
  	\node (F) at (0, 8) {};
  	\node (G) at (0, 6) {};
  	\node  (H) at (0, 4) {};
  	\node  (I) at (0, 2) {};
  	\node (J) at (0, 0) {};
  	\node  (K) at (0, -2) {};
  	\node  (L) at (0, -4) {};

  \foreach \from/\to in {B/C, D/E, F/G, H/I, J/K}
    	\draw (\from) -- (\to)[red];

  \path[every node/.style={}]
	(A) edge [red, bend right = 70] node [black] {\,\,$x$} (B)
	(C) edge [bend right = 70] node [black] {\,\,$\overline{x}$} (D)
	(E) edge [red,bend right = 70] node [black] {\,\,$y$} (F)
	(G) edge [bend right = 70] node [black] {\,\,$\overline{y}$} (H)
	(I) edge [bend right = 70] node [black] {\,\,$z$} (J)
	(K) edge [red,bend right = 70] node [black] {\,\,$\overline{z}$} (L)

	(A) edge [bend left = 70] node [black] {} (B)
	(C) edge [red,bend left = 70] node [black] {} (D)
	(E) edge [bend left = 70] node [black] {} (F)
	(G) edge [red,bend left = 70] node [black] {} (H)
	(I) edge [red,bend left = 70] node [black] {} (J)
	(K) edge [bend left = 70] node [black] {} (L);

	\node  (A2) at (20, 18) {};
  	\node  (B2) at (20, 16) {};
  	\node (C2) at (20, 14) {};
  	\node  (D2) at (20, 12) {};
  	\node  (E2) at (20, 10) {};
  	\node (F2) at (20, 8) {};
  	\node (G2) at (20, 6) {};
  	\node  (H2) at (20, 4) {};
  	\node  (I2) at (20, 2) {};
  	\node (J2) at (20, 0) {};
  	\node  (K2) at (20, -2) {};
  	\node  (L2) at (20, -4) {};

  \foreach \from/\to in {B2/C2, D2/E2, F2/G2, H2/I2, J2/K2}
    	\draw (\from) -- (\to) [red];

  \path[every node/.style={}]
	(A2) edge [bend right = 70] node [black] {\,\,$\overline{x}$} (B2)
	(C2) edge [red,bend right = 70] node [black] {\,\,$y$} (D2)
	(E2) edge [red,bend right = 70] node [black] {\,\,$\overline{z}$} (F2)
	(G2) edge [red,bend right = 70] node [black] {\,\,$x$} (H2)
	(I2) edge [bend right = 70] node [black] {\,\,$\overline{y}$} (J2)
	(K2) edge [red,bend right = 70] node [black] {\,\,$\overline{z}$} (L2)

	(A2) edge [red,bend left = 70] node [black] {} (B2)
	(C2) edge [bend left = 70] node [black] {} (D2)
	(E2) edge [bend left = 70] node [black] {} (F2)
	(G2) edge [bend left = 70] node [black] {} (H2)
	(I2) edge [red,bend left = 70] node [black] {} (J2)
	(K2) edge [bend left = 70] node [black] {} (L2);
  	\node (X1) [fill=white, text opacity=0] at (-2,15)  {$X$};
 	\node (Y1) [fill=white, text opacity=0] at (-2,7)  {$X$};
 	\node (Z1) [fill=white, text opacity=0] at (-2,-1)  {$X$};

 	\draw (X1) to [bend left=30] (-0.7,17); 
 	\draw (X1) to [bend right=30] (-0.7,13); 
 	\draw (Y1) to [bend left=30] (-0.7,9); 
 	\draw (Y1) to [bend right=30] (-0.7,5); 
 	\draw (Z1) to [bend left=30] (-0.7,1); 
 	\draw (Z1) to [bend right=30] (-0.7,-3); 

 \draw plot [smooth] coordinates {(-2.35,15.4)(-1.65,14.6)};
 \draw plot [smooth] coordinates {(-2.35,14.6)(-1.65,15.4)};

 \draw plot [smooth] coordinates {(-2.35,7.4)(-1.65,6.6)};
 \draw plot [smooth] coordinates {(-2.35,6.6)(-1.65,7.4)};

 \draw plot [smooth] coordinates {(-2.35,-0.6)(-1.65,-1.4)};
 \draw plot [smooth] coordinates {(-2.35,-1.4)(-1.65,-0.6)};
  	\node (X2) [fill=white, text opacity=0] at (22,13)  {$X$};
 	\node (Y2) [fill=white, text opacity=0] at (22,1)  {$X$};

 	\draw (X2) to [bend right=15] (20.7,17); 
 	\draw (X2) to (20.7,13); 
 	\draw (X2) to [bend left=15] (20.7,9); 
 	\draw (Y2) to [bend right=15] (20.7,5); 
 	\draw (Y2) to (20.7,1); 
 	\draw (Y2) to [bend left=15] (20.7,-3); 

  	\node (X) [draw=none, fill=white, text opacity=1, font=\large] at (22,13)  {$\vee$};
  	\node (Y) [draw=none, fill=white, text opacity=1, font=\large] at (22,1)  {$\vee$};
 	\draw (0.7,17) to (19.3,5); 
 	\draw (0.7,13) to (19.3,17); 
 	\draw (0.7,9) to (19.3,13); 
 	\draw (0.7,5) to (19.3,1); 
 	\draw (0.7,-2.7) to (19.3,9); 
 	\draw (0.7,-3) to (19.3,-3); 

 	\draw [red] (A) to (A2); 
 	\draw [red] (L) to (L2); 

  	\node (P1) [fill=white, text opacity=0] at (13,15.7)  {$X$};
  	\node (P2) [fill=white, text opacity=0] at (13,11.7)  {$X$};
  	\node (P3) [fill=white, text opacity=0] at (13,9)  {$X$};
  	\node (P4) [fill=white, text opacity=0] at (13,5)  {$X$};
  	\node (P5) [fill=white, text opacity=0] at (13,2.4)  {$X$};
  	\node (P6) [fill=white, text opacity=0] at (13,-3)  {$X$};

 \draw plot [smooth] coordinates {(12.65,16.1)(13.35,15.3)};
 \draw plot [smooth] coordinates {(12.65,15.3)(13.35,16.1)};

 \draw plot [smooth] coordinates {(12.65,12.1)(13.35,11.3)};
 \draw plot [smooth] coordinates {(12.65,11.3)(13.35,12.1)};

 \draw plot [smooth] coordinates {(12.65,9.4)(13.35,8.6)};
 \draw plot [smooth] coordinates {(12.65,8.6)(13.35,9.4)};

 \draw plot [smooth] coordinates {(12.65,5.4)(13.35,4.6)};
 \draw plot [smooth] coordinates {(12.65,4.6)(13.35,5.4)};

 \draw plot [smooth] coordinates {(12.65,2.8)(13.35,2)};
 \draw plot [smooth] coordinates {(12.65,2)(13.35,2.8)};

 \draw plot [smooth] coordinates {(12.65,-2.6)(13.35,-3.4)};
 \draw plot [smooth] coordinates {(12.65,-3.4)(13.35,-2.6)};
 \draw plot [smooth] coordinates {(26,7)(25,16)(20,20)(14, 19)(13, 18)};
 \draw plot [smooth] coordinates {(26,7)(25,-2)(20,-6)(14, -5)(13, -4)};

	\node (X4) [fill=white, text opacity=0] at (26,7)  {$X$};

  	\node (X5) [draw=none, fill=white, text opacity=1, font=\large] at (26,7)  {$\vee$};
\end{tikzpicture}
\caption{Transformation from $F = (\neg x \lor y \lor \neg z) \land (x \lor \neg y \lor \neg z)$ to a planar, cubic, $3$-connected graph}
\label{Example}
\end{figure}
\end{minipage}

In a similar way you can check that every Hamiltonian cycle corresponds to a truth assignment of $F$. Every Hamiltonian cycle has to use exactly one of the edges $v_{i1}v_{i2}$ and $v_{i3}v_{i4}$ that is connected to the clauses, because of the XOR-module connecting these two edges. This defines a truth assignment of $F$, by making a variable true if the parallel edge $v_{i1}v_{i2}$ connected to the clauses is used in the Hamiltonian cycle and false otherwise.
\end{proof}
This completes the proof of the reduction from $3$-SAT-CNF to this decision problem corresponding to Tait's conjecture. So finding a Hamiltonian cycle in a cubic, $3$-connected, planar graph is \NPP-complete.
\end{proof}

\begin{corollary} \label{pad} Finding a Hamiltonian \emph{path} in a planar, cubic, $3$-connected graph is \NPP-complete.
\end{corollary}
\begin{proof}
It is clear that this problem is in NP, so we only have to give a reduction from an \NPP-complete problem. We use the same reduction from $3$-SAT-CNF as used in the proof of the theorem above, except that we replace the $2$-input-OR-module between $v_{11}w_{11}$ and $v_{n4}w_{m6}$ by a XOR-module. A truth assignment to the variables can be transformed into a Hamiltonian \emph{path} by using the same method used in the proof of Theorem \ref{HAM*}. Only the last step, were we used  $v_{11}w_{11}$, can be omitted here, because only a Hamiltonian path is needed. This satisfies all modules, because the only new module is the XOR-module between $v_{11}w_{11}$ and $v_{n4}w_{m6}$, but this one is satisfied because we do use $v_{n4}w_{m6}$.

It remains to show that a given Hamiltonian path in the graph corresponds to a truth assignment of the variables which makes $F$ true. A Hamiltonian path connects all vertices in the graph, so it also satisfies all modules. In particular it satisfies the XOR-module between $v_{11}w_{11}$ and $v_{n4}w_{m6}$. Without loss of generality suppose that the Hamiltonian path passes trough $v_{11}w_{11}$ and not through $v_{n4}w_{m6}$. The part of the variables is only connected to the part of the clauses by XOR-modules, so all vertices in the part of the variables have to be visited in a row by the Hamiltonian path. Stated otherwise: the Hamiltonian path can only change once from a vertex $v_{ab}$ to a vertex $w_{cd}$, which is done by using edge $v_{11}w_{11}$. So the Hamiltonian path has to start in $v_{n4}$, then pass through all vertices of the part of the variables, use edge $v_{11}w_{11}$, then pass through all vertices of the part of the clauses and end in $w_{m6}$. This gives a truth assignment of the variables which makes $F$ true in the same way as we have seen in the proof of Theorem \ref{HAM*}.
\end{proof}

\section{Tutte's conjecture and the necessity of planarity}
We will now consider a decision problem corresponding to Tutte's conjecture. An instance of this problem is a cubic, $3$-connected, bipartite graph. The question is whether or not there exists a Hamiltonian cycle in this graph.

\begin{thrm} \label{HAM**}
Finding a Hamiltonian cycle in a cubic, $3$-connected, bipartite graph is \NPP-complete.
\end{thrm}
\begin{proof}
This proof is similar to the proof of Theorem \ref{HAM*} and can be found in \cite{NPcompletenessbipartite}. The only difference is that the non-bipartite, planar Tutte fragments \ref{Tutte's fragment} in the XOR-module are replaced by bipartite, non-planar Horton fragments \ref{Horton fragment}.
\end{proof}

\section{Complexity of Barnette's conjecture}
We will now consider a decision problem corresponding to Barnette's conjecture. An instance of this problem is a cubic, $3$-connected, bipartite, planar graph. The question is whether or not there exists a Hamiltonian cycle in this graph. For the complexity of this problem we use the concept of a required edge fragment, see Definition \ref{requirededgefragment}. A required edge fragment has three half edges of which one is the \emph{required edge}. This required edge has to be used in any Hamiltonian cycle in the graph that contains the fragment. This is the same property that Tutte's fragment and the Horton fragment have.

\begin{thrm} \label{HAM***}
If there exists a cubic, $3$-connected, bipartite, planar required edge fragment, then finding a Hamiltonian cycle in a cubic, $3$-connected, bipartite, planar graph is \NPP-complete.
\end{thrm}
\begin{proof}
This proof is similar to the proof of Theorem \ref{HAM*}. The only difference is that the non-bipartite, planar Tutte fragments in the XOR-module are replaced by cubic, $3$-connected, bipartite, planar required edge fragments.
\end{proof}

The reduction used for these proofs is very similar to the reduction from $3$-SAT-CNF to the common Hamiltonian cycle problem. For the common Hamiltonian cycle no Tutte or Horton fragments are used, but a slightly different XOR-module is used to obtain the same property of this module. This reduction can be found in \cite{Terwijn}.

\section{Construction of counterexamples}
Every counterexample for Tait's and Tutte's conjecture that we discussed did contain a Hamiltonian {\it path}. This gave rise to the conjectures that every Tait graph and every Tutte graph contains a Hamiltonian path. We can disprove this conjecture immediately with the reduction given above and Corollary \ref{pad}. By making up a non-satisfiable Boolean formula $F$, we can construct the corresponding graph using Corollary \ref{pad} and Theorem \ref{HAM*} for a Tait graph and Corollary \ref{pad} and Theorem \ref{HAM**} for a Tutte graph. This graph will not contain a Hamiltonian path, because there is no possible truth assignment to $F$. An easy example of a non-satisfiable formula in $3$-CNF is $$F = (x \lor y \lor z) \land (x \lor y \lor \neg z) \land (x \lor \neg y \lor z) \land (x \lor \neg y \lor \neg z) \land (\neg x \lor y \lor z) \land (\neg x \lor y \lor \neg z) \land (\neg x \lor \neg y \lor z) \land (\neg x \lor \neg y \lor \neg z).$$
The corresponding graph has $12$ vertices for its $3$ variables and $48$ vertices for its $8$ clauses. For every variable there is a XOR-module, which itself has $36$ vertices. For every clause there is a $3$-input-OR-module, which consists of $510$ vertices. Furthermore we have a XOR-module for every literal used in the clauses, which are in total $24$. Finally there is an extra XOR-module connecting $v_{11}w_{11}$ and $v_{n4}w_{m6}$. This brings the total amount of vertices to $12 + 48 + 3 \cdot 36 + 8 \cdot 510 + 24 \cdot 36 + 36 = 5148$. Several XOR-modules intersect each other, which makes the graph non-planar. This problem can be solved by replacing every intersection by a new module, which can be found in \cite{GareyJohnsonTarjan}. This adds another $152$ vertices for every intersection. Even for this easy non-satisfiable formula in $3$-CNF, the corresponding counterexample is already very big. So a different method has to be used to find small counterexamples for Barnette's conjecture.

\chapter{Planar 4-connected graphs}
\label{4connected}

In this chapter we will give a sketch of the proof that every $4$-connected, planar graph $G$ is Hamiltonian. This is stated as Theorem \ref{hamilton} and the proof we sketch is based on a proof in \cite{Tutte1956} by Tutte. To understand the proof we will first introduce some definitions and notations used in \cite{Tutte1956}. In this chapter a graph is defined by a set of edges. Vertices of a graph are the endpoints of its edges. In this case a subgraph $H$ of a graph $G$ is a subset of the edges of $G$. The vertices of $H$ are the endpoints of the edges in this subset. By $\alpha_0(G)$ we denote the number of vertices of $G$ and by $\alpha_1(G)$ we denote the number of edges of $G$. We say that two subgraphs meet, when they have a common edge. We call two graphs disjoint if they have no edge in common; they can have a vertex in common.
\begin{definition}
An {\bf isthmus} of graph $G$ is an edge which is not part of any cycle of $G$.
\end{definition}
\begin{definition}
A {\bf vertex of attachment} is a common vertex of $H \subseteq G$ and $G-H$. We denote the total number of vertices of attachment by the {\bf attachment number $w(H)=w(G-H)$}.
\end{definition}
\begin{definition}
Let $J$ be a cycle of a graph $G$. A subset $H$ of $G-J$ is called {\bf $J$-bounded} if all its points of attachment are vertices of $J$.
\end{definition}
\begin{definition}
Let $J$ be a cycle of a graph $G$ and let $G-J$ be non-empty, then it has a unique expression as a union of disjoint minimal non-empty $J$-bounded subsets of $G-J$. These subsets are the {\bf bridges} of $J$ in $G$.
\end{definition}
\begin{definition}
Cycle $J$ is called the {\bf bounding cycle} of two {\bf residual domains} in which it divides the plane. In  a planar embedding the residual domains of bounding cycle $J$ are the 'inside' and the 'outside' of cycle $J$.
\end{definition}
\begin{definition}
\label{terminal}
Let edge $e$ not be an isthmus of a planar graph $G$ and let $D$ be a residual domain of cycle $J$ with $e \in J$. Then $D$ is a {\bf terminal domain of} $e$ if no cycle of $G$ containing $e$ has a residual domain that is a proper subset of $D$. In this case cycle $J$ is called a {\bf terminal cycle}.\end{definition}

Note that the definitions of a terminal domain and a terminal cycle are dependent on the planar embedding of the graph $G$. Similar to definition \ref{terminal} we can define a terminal domain of a path $L$, which is a subset of cycle $J$.

\begin{definition}
Let path $L$ be a subset of cycle $J$ of a planar graph $G$ and let $D$ be a residual domain of cycle $J$. Then $D$ is a {\bf terminal domain of} $L$ if no cycle of $G$ containing $L$ has a residual domain that is a proper subset of $D$. In this case cycle $J$ is called a {\bf terminal cycle}.
\end{definition}

Besides these definitions we also use another theorem and some lemma's to sketch the proof of Theorem \ref{hamilton}. The proof of Theorem \ref{stelling} can be found in \cite{Tutte1956} as the proof of Theorem I on page 100.

\begin{thrm}
\label{stelling}
Let $G$ be a planar graph. Let $e$ be an edge of $G$ which is not an isthmus and let $e'$ be an edge distinct from $e$ of a terminal cycle of $e$. Then there exists a cycle $J$ of $G$ having the following properties:\\
$(i)$ $e, e' \in J$.\\
$(ii)$ If $B$ is a bridge of $J$ in $G$, then $w(B) \leq 3$.\\
$(iii)$ If $B$ is a bridge of $J$ in $G$ which meets a terminal cycle of $e$, then $w(B)=2$.
\end{thrm}

\begin{lemma}
\label{connected}
Let $B$ be a bridge of a cycle $J$ of a $k$-connected planar graph $G$, for $k \geq 2$. Then if $w(B) < k$ one of the following alternatives is true:\\
$(i)$ $B$ has just one edge and both ends of this edge are vertices of $J$.\\
$(ii)$ $\alpha_1(J)=w(B)$ and $B$ is the only bridge of $J$ in $G$ having more than one edge.
\end{lemma}
\begin{proof}
Let $B$ be a bridge of a cycle $J$, with $w(B)<k$ and $B$ has more than one edge. Now we have to show that $\alpha_1(J)=w(B)$ and $B$ is the only bridge of $J$ in $G$ having more than one edge. Because $B$ has more than one edge, there is a vertex in $B$ that is not a vertex of attachment. We have assumed that $w(B)<k$, so after removing all vertices of attachment, graph $G$ is still connected. This implies that the only vertices left, after removing the vertices of attachment of $B$, are vertices of $B$. Otherwise the removal results in a disconnected graph. This implies that $J$ consists exactly of the $w(B)$ vertices and thus $\alpha_0(J)=\alpha_1(J)=w(B)$. Furthermore this implies that there are no vertices in the graph other than vertices of $J$ or vertices of $B$. So $B$ is the only bridge having more than one edge, because all other bridges only consist of vertices of $J$ and thus can have at most one edge.
\end{proof}

\begin{lemma}
\label{bridge}
Let $B$ be a bridge of a cycle $J$ and $D$ a residual domain of $J$ in graph $G$. Let $x$ and $y$ be distinct vertices of attachment of $B$ and let $L$ a path from $x$ to $y$ only using edges  of $J$. Let $T$ be the terminal domain of $L$ contained in $D$. Then there is a path $M$ from $x$ to $y$ in $D$, not using any points of $J$ other than $x$ and $y$, such that $L \cup M$ is the bounding cycle of $T$.
\end{lemma}
The proof of Lemma \ref{bridge} can be found in \cite{Tutte1956} as the proof of statement 4.3 on page 104.

\begin{lemma}
\label{isthmus}
A 2-connected planar graph $G$ such that $\alpha_1(G) \geq 2$ has no isthmus.
\end{lemma}
\begin{proof}
Suppose $G$ has an isthmus $e=v_1v_2$ and $G$ is 2-connected. Define $$A_1 = \{v \in V(G) | \mbox{there exists a path from} \ v_1 \ \mbox{to} \ v \ \mbox{not using} \ v_2\}$$ and $$A_2 = \{v \in V(G) | \mbox{there exists a path from} \ v_2 \ \mbox{to} \ v \ \mbox{not using} \ v_1\}.$$ We have assumed that $e$ is an isthmus, so $A_1 \cap A_2 = \emptyset$, otherwise $e$ would be part of a cycle. By removing $v_1$ we can disconnect graph $G$ into $A_1 - \{v_1\}$ and $A_2$ and the corresponding edges. If $A_1 - \{v_1\} = \emptyset$, then we can remove $v_2$ to obtain a disconnected graph. This is a contradiction with the fact that $G$ is 2-connected.
\end{proof}

\begin{thrm}
\label{hamilton}
Every $4$-connected, planar graph $G$ with at least $2$ edges is Hamiltonian. Furthermore if $G$ has no parallel edges and $e$, $e'$ are different edges of the same terminal cycle of $G$, then $G$ has a Hamiltonian cycle that contains $e$ and $e'$.
\end{thrm}
\begin{proof}
First we will prove the second part of the theorem, then this part is used to prove the first part of the theorem. So suppose $G$ has no parallel edges and let $e$ and $e'$ be different edges of the same terminal cycle $T$ of $G$. Edge $e$ is clearly not an isthmus, because it is part of cycle $T$, and $e'$ is part of the same terminal cycle, so we can apply Theorem \ref{stelling}. This theorem states that there exists a cycle $J$ that satisfies the three properties stated in Theorem \ref{stelling}. Now we make a case distinction on the number of edges in the bridges of $J$ in $G$.

Suppose every bridge of $J$ in $G$ consists of exactly one edge with its endpoints in cycle $J$. In this case $J$ is clearly a Hamiltonian cycle, because it visits every vertex exactly once. Furthermore $J$ contains edges $e$ and $e'$ because of the first property of Theorem \ref{stelling}.

Suppose that not every bridge of $J$ consists of exactly one edge. The second property of Theorem \ref{stelling} states that every bridge $B$ of $J$ satisfies $w(B) \leq 3$, so $w(B)<4$. This implies that we can apply Lemma \ref{connected}. In this case we have assumed that there is at least one bridge $B'$ that consists of more than one edge. Lemma \ref{connected} states that $\alpha_1(J) = w(B')$ and $B'$ is the only bridge of $J$ in $G$ having more than one edge. We know by the second property of Theorem \ref{stelling} that $w(B') \leq 3$ and we know that $\alpha_1(J) \geq 3$, because there are no parallel edges in $G$, so any cycle has at least $3$ edges. Combining these two facts we get that $\alpha_1(J) = w(B')=3$. Now we can conclude that $B'$ is the \emph{only} bridge of $J$ in $G$, because there are no parallel edges and all $3$ vertices of $J$ are already pairwise connected. By property $(i)$ of Theorem \ref{stelling} we know that $e, e' \in J$. Let $x$ and $y$ be the distinct endpoints of $e$, then we can apply Lemma \ref{bridge} to this situation. This lemma tells us that there exists a path $M$ from $x$ to $y$, not using any points of $J$ other than $x$ and $y$, such that $\{e\} \cup M$ is the bounding cycle of a terminal domain $T$ of $e$. Clearly this path $M$ is a subset of bridge $B'$, because this is the only bridge of $J$. So we can conclude that bridge $B'$ meets a terminal cycle of edge $e$. Property $(iii)$ of Theorem \ref{stelling} now tells us that $w(B')=2$. This contradicts the fact that $w(B')=3$. So apparently no bridge of $J$ can have more than one edge. We have already shown that $J$ is a Hamiltonian cycle if every bridge of $J$ has exactly one edge, so this completes the proof of the second part of the theorem.

Now we will prove the first part of the theorem. Suppose $G$ is a $4$-connected, planar graph with at least $2$ edges. Then $G$ is also $2$-connected, so by Lemma \ref{isthmus} we know that $G$ has no isthmus and this implies that $G$ has at least one cycle. If $\alpha_0(G)=2$, then we clearly have a Hamiltonian cycle, so from now on we assume that $\alpha_0(G) \geq 3$. Remove all parallel edges of graph $G$ to obtain subgraph $H$. This graph $H$ has the same vertices as graph $G$, so $H$ is still $4$-connected. By Lemma \ref{isthmus} we know that graph $H$ does not have an isthmus either. This implies that $H$ has at least one cycle and thus also a terminal cycle. Now we can apply the second part of this theorem, which we have proved already. From the second part of this theorem we can conclude that graph $H$ contains a Hamiltonian cycle. Graph $G$ has the same vertices as graph $H$ and the edges of graph $H$ are a subset of the edges of graph $G$, so graph $G$ contains the same Hamiltonian cycle. This completes the proof.
\end{proof}

\chapter{Barnette's Conjecture up to 64 vertices}
\label{Proofonvertices}
Holton, Manvel and McKay proved in 1984 that Barnette's conjecture is true for all graphs up to 64 vertices. In this chapter we will give an idea of the proof. The complete proof can be found in \cite{HMM1985}.

Before we can start with the actual proof, we use some definitions:
\begin{definition} A {\bf $k$-face} is a face bounded by $k$ edges and a {\bf $k$-cut} is a set of $k$ edges whose removal separates the graph $G$ into two parts, each with at least two vertices.  \end{definition}
\begin{definition} A 4-cut is called {\bf essential} if neither part is a 4-face and  {\bf major} if neither part is a 4-face, graph $R_1$ or graph $R_2$ of Figure \ref{r1r2}.   \end{definition}
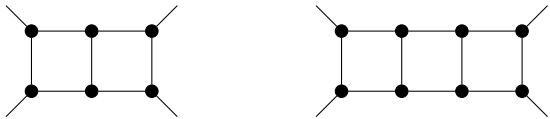
\begin{figure}[h!]%
	\centering
          \begin{minipage}{1.3in}
		\begin{tikzpicture}
 		 [scale=0.8,auto=left,every node/.style={circle, draw,fill=black,scale=0.5}]
  		\node (A) at (0,0) {};
		\node (B) at (1,0) {};
		\node (C) at (2,0) {};
		\node (E) at (0,1) {};
		\node (F) at (1,1) {};
		\node (G) at (2,1) {};
		\node [white] (I) at (-0.5,-0.5) {};
		\node [white] (J) at (-0.5,1.5) {};
		\node [white] (K) at (2.5,-0.5) {};
		\node [white] (L) at (2.5,1.5) {};
		\foreach \from/\to in {A/B,B/C,F/G,F/E,A/E,B/F,I/A,J/E,K/C,G/L}
		\draw (\from) -- (\to);
		\draw (C) -- (G);
		\end{tikzpicture}%
	\end{minipage}
\qquad
           \begin{minipage}{1in}%
           	\begin{tikzpicture}
 		 [scale=0.8,auto=left,every node/.style={circle, draw,fill=black,scale=0.5}]
  		\node (A) at (0,0) {};
		\node (B) at (1,0) {};
		\node (C) at (2,0) {};
		\node (D) at (3,0) {};
		\node (E) at (0,1) {};
		\node (F) at (1,1) {};
		\node (G) at (2,1) {};
		\node (H) at (3,1) {};
		\node[white]  (I) at (-0.5,-0.5) {};
		\node[white]  (J) at (-0.5,1.5) {};
		\node[white]  (K) at (3.5,-0.5) {};
		\node[white]  (L) at (3.5,1.5) {};
		\foreach \from/\to in {A/B,B/C,C/D,D/H,H/G,F/G,F/E,A/E,B/F,I/A,J/E,K/D,L/H}
		\draw (\from) -- (\to);
		\draw (C) -- (G);
		\end{tikzpicture}	
		\end{minipage}%
          \caption{Left: graph $R_1$. Right: graph $R_2$.}%
          \label{r1r2}%
\end{figure}

\begin{definition} A {\bf C3CBP} is a cubic, 3-connected, bipartite, planar graph.  A {\bf C3CBP4} is a C3CBP with no 3-cuts or essential 4-cuts, so that any 4-cut has at least one part that is a 4-face.  A {\bf C3CBP4*} is a C3CBP with no 3-cuts or major 4-cuts, so that any 4-cut has at least one part that is a 4-face, $R_1$ or $R_2$. (Some examples can be found in Figure \ref{c3cbp}). \end{definition}

\begin{figure}%
            \centering
            \parbox{1.3in}{ \begin{tikzpicture}
  [scale=.5,auto=left,every node/.style={circle, draw,fill=black,scale=0.5}]
  	\node (A) at (0,0) {};
  	\node (B) at (6,0) {};
  	\node (C) at (6,6) {};
  	\node (D) at (0,6) {};
  	\node (E) at (1.5,1.5) {};
  	\node (F) at (4.5,1.5) {};
  	\node (G) at (4.5,4.5) {};
  	\node (H) at (1.5,4.5) {};

  \foreach \from/\to in {A/B, B/C, C/D, D/A, E/F, F/G, G/H, H/E, A/E, B/F, C/G, D/H}
    	\draw (\from) -- (\to);
 
\end{tikzpicture}}%
            \qquad
            \begin{minipage}{1.5in}%
              \begin{tikzpicture}
  [scale=.37,auto=left,every node/.style={circle, draw,fill=black,scale=0.5}]
  	\node (A) at (2,0) {};
  	\node (B) at (6,0) {};
  	\node (C) at (8,4) {};
  	\node (F) at (0,4) {};
  	\node (D) at (6,8) {};
  	\node (E) at (2,8) {};
  	\node (G) at (2,4) {};
  	\node (H) at (3,2) {};
	\node (I) at (5,2) {};
  	\node (J) at (6,4) {};
  	\node (K) at (5,6) {};
  	\node (L) at (3,6) {};
  \foreach \from/\to in {A/B,B/C,C/D,D/E,E/F,F/A,G/H,H/I,I/J,J/K,K/L,L/G,A/H,B/I,C/J,D/K,L/E,F/G}
    	\draw (\from) -- (\to); 
\end{tikzpicture}

            \end{minipage}%
	  \begin{minipage}{1.2in}%
              \begin{tikzpicture}
  [scale=.6,auto=left,every node/.style={circle, draw,fill=black,scale=0.5}]
  	\node (A) at (0,5) {};
  	\node (B) at (5,0) {};
  	\node (C) at (0,3) {};
  	\node (D) at (2,3) {};
  	\node (E) at (2,5) {};
  	\node (F) at (1,3) {};
  	\node (G) at (1,4) {};
  	\node (H) at (2,4) {};
	\node (I) at (3,0) {};
  	\node (J) at (5,2) {};
  	\node (K) at (3,2) {};
  	\node (L) at (3,1) {};
  	\node (M) at (4,1) {};
  	\node (N) at (4,2) {};
 \foreach \from/\to in {A/C,A/E, A/G,F/G,G/H,C/F,F/D,E/H,H/D,D/K,C/I,E/J,K/N,N/J,K/L,L/I,L/M,M/N,I/B,J/B,B/M}
    	\draw (\from) -- (\to);
	\end{tikzpicture}

            \end{minipage}%
            \caption{From left to right: $C3CBP4$, $C3CBP4^*$ and $C3CBP$. }%
            \label{c3cbp}%
          \end{figure}
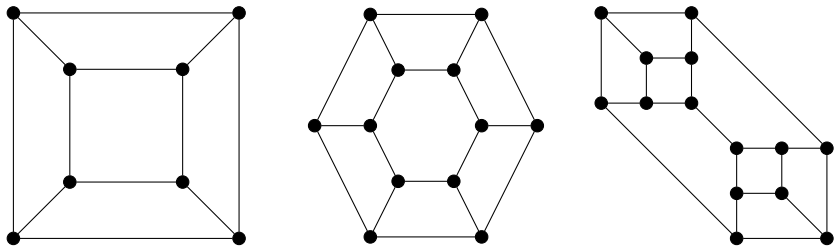

\begin{figure}[h!]%
	\centering
          \begin{minipage}{1.1in}
		\begin{tikzpicture}
 		 [scale=0.8,auto=left,every node/.style={circle, draw,fill=black!20,scale=0.5}]
  		\node (A) at (0,0) {A};
		\node (B) at (1,0) {B};
		\node (C) at (2,0) {C};
		\node (D) at (3,0) {D};
		\node (E) at (0,1) {E};
		\node (F) at (1,1) {F};
		\node (G) at (2,1) {G};
		\node (H) at (3,1) {H};
		\node[white] (I) at (-0.5,-0.5) {};
		\node[white]  (J) at (-0.5,1.5) {};
		\node[white]  (K) at (3.5,-0.5) {};
		\node[white]  (L) at (3.5,1.5) {};
		\foreach \from/\to in {A/B,B/C,C/D,D/H,H/G,F/G,F/E,A/E,B/F,I/A,J/E,K/D,L/H}
		\draw (\from) -- (\to);
		\draw (C) -- (G);
		\end{tikzpicture}%
	\end{minipage}
\qquad
	\begin{minipage}{0.1in}%
          $\Rightarrow$
          \end{minipage}%
\qquad
           \begin{minipage}{1in}%
           	\begin{tikzpicture}
 		 [scale=0.8,auto=left,every node/.style={circle, draw,fill=black!20,scale=0.5}]
  		\node (A) at (0,0) {A};
		\node (D) at (1,0) {D};
		\node (E) at (0,1) {E};
		\node (H) at (1,1) {H};
		\node[white]  (I) at (-0.5,-0.5) {};
		\node[white]  (J) at (-0.5,1.5) {};
		\node[white]  (K) at (1.5,-0.5) {};
		\node[white]  (L) at (1.5,1.5) {};
 		 \foreach \from/\to in {A/D,D/H,H/E,A/E,I/A,J/E,K/D,L/H}
    		\draw (\from) -- (\to);
		\end{tikzpicture}
	\end{minipage}%
          \caption{Reduction $R_2$ and central edges BF and CG.}%
          \label{centraledge}%
\end{figure}
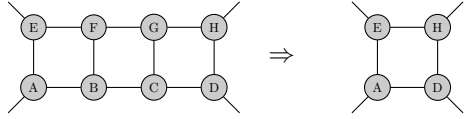

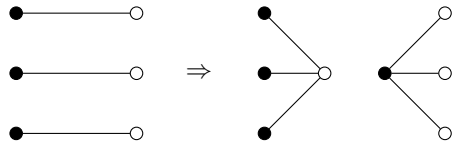
\begin{figure}[b]%
	\centering
	\begin{minipage}{0.6in}
		\begin{tikzpicture}
  		[scale=0.8,auto=left,every node/.style={circle, draw,fill=white,scale=0.5}]
  		\node[black] (A) at (0,0) {};
		\node[black] (B) at (0,1) {};
		\node[black] (C) at (0,2) {};
		\node (F) at (2,0) {};
		\node (G) at (2,1) {};
		\node (H) at (2,2) {};
		  \foreach \from/\to in {A/F,B/G,C/H}
    		\draw (\from) -- (\to);
		\end{tikzpicture}%
	\end{minipage}
\qquad
	\begin{minipage}{0.1in}%
              	$\Rightarrow$
	\end{minipage}%
\qquad
           \begin{minipage}{1.2in}%
            	\begin{tikzpicture}
  		[scale=0.8,auto=left,every node/.style={circle,draw,fill=white,scale=0.5}]
  		\node[black] (A) at (0,0) {};
		\node[black] (B) at (0,1) {};
		\node[black] (C) at (0,2) {};
		\node (D) at (1,1) {};
		\node[black] (E) at (2,1) {};
		\node (F) at (3,0) {};
		\node (G) at (3,1) {};
		\node (H) at (3,2) {};
		 \foreach \from/\to in {A/D,B/D,C/D,E/F,E/G,E/H}
    		\draw (\from) -- (\to);
		\end{tikzpicture}
	\end{minipage}%
	\caption{Split a graph with a 3-cut in two graphs.}%
	\label{cut}%
\end{figure}

\begin{definition} We say a C3CBP is H if it has a Hamiltonian cycle, H$^+$ (H$^-$) if it has a Hamiltonian cycle through (avoiding) any specified edge. We call a C3CBP H$^{+-}$ if any two edges can be specified such that we can find a Hamiltonian cycle that uses one edge and avoids the other. H$^*$ is defined later. N, N$^+$, N$^-$, N$^{+-}$ and N$^*$ are numbers corresponding with  the properties above. Thus $N$ is is the largest number for which every C3CBP on at most $N$ vertices is $H$. If no such $N$ exists, then all C3CBP's are H. 
\end{definition}

\begin{definition} A {\bf central edge} is an edge $BF$ or $CG$ in $R_2$ in Figure \ref{centraledge}
\end{definition}

\begin{definition} A graph $G$ is $H^*$ if it is $H^{+-}$ or if it satisfies the following condition. It contains a unique subgraph $R_2$ and any pair of edges, one to use, one to avoid, can be specified for a Hamiltonian cycle in $G$, unless the edge to be used is a central edge of $R_2$ and the edge to be avoided contains no vertex of that $R_2$ subgraph.
\end{definition}

Our goal is to show that $N \geq 64$. The actual proof consists of two parts:
\begin{enumerate}
\item Show that $N \geq \min \{ N^* + 20, 2N^*-6\}$;
\item Show that $N^* \geq 44$ by using the computer.
\end{enumerate}
In this way, if we have better computers to determine $N^*$, we still can use the relation between $N$ and $N^*$. It is clear that $N^* \geq 44$ implies $N \geq \min \{ 44 + 20, 2\cdot 44-6\}=64$.

\section{The relation between $N$ and $N^*$}
One can show that $N \geq \min \{ N^* + 20, 2N^*-6\}$. When $N^* \geq 26$, the lower bound will be given by $N^*+20$. In this section we will be really short and mostly refer to the proof in the original article \ref{HMM18985}. 

\begin{figure}[t]
\centering
\includegraphics[scale=0.5]{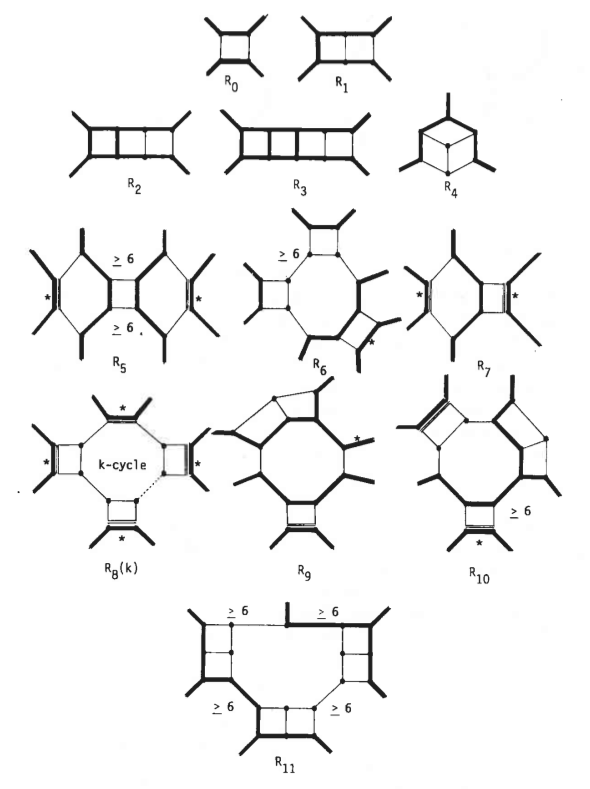}
\caption{Reductions}
\label{reductions}
\end{figure}

Here we use the reductions in Figure \ref{reductions}. In a reduction we remove the light edges.

In the proof three types of argument are used: 
\begin{enumerate}
\item If $G$ is the smallest non-Hamiltonian C3CBP, then:
\begin{enumerate}
\item $G$ has a 3-cut and $\#V(G) \geq 2N^*+2$, or
\item $G$ has no 3-cut, but has an essential 4-cut, and $\#V(G) \geq 2N^*-4$, or
\item $G$ is C3CBP4.
\end{enumerate}
\item A graph with a 3-cut or essential 4-cut can be separated at that cut. We introduce some new vertices for the edges in the cut, the example for the 3-cut is shown in Figure \ref{cut}. We can prove that both new graphs are C3CBP and combine the Hamiltonian cycles of both pieces to one Hamiltonian cycle for the whole graph. The prove for a 3-cut can be found in Lemma \ref{3cut}. The prove for a essential 4-cut can be found in \cite{HMM1985}.
\item A $C3CBP4$ can be reduced with one of the reductions $R_5, R_6$ or $R_8(k), k \geq 8$. This reduction reduces the number of vertices by at least 12. Now we still should be able to find a Hamiltonian cycle in the smaller graph which can be extended to one in the original graph. We can reduce the number of vertices even further (by 8) when we use the lemma's in the article of Holton, Manvel and McKay. This gives us a bound $N\geq N^* + 20$.
\end{enumerate}

When we combine these arguments we have the main idea of the proof. For the complete proof we refer to the original article of Holton, Manvel and McKay \cite{HMM1985}.

\section{Computer generation to determine $N^*$}
To find a lower bound on $N^*$ Holton, Manvel and McKay have generated all C3CBP's up to 40 vertices and all C3CBP's on 42 and 44 vertices that do not contain subgraphs $R_2$ or $R_4$. They proved and used the following lemma:

\begin{lemma}
Suppose all C3CBP's containing no subgraphs $R_2$ or $R_4$ on up to $n$ vertices are $H^{+-}$. Then $N^*=\min(n, 4+N^{+-})$.
\end{lemma}

Proof of this lemma can be found in the original article. 

Holton, Manvel and McKay tested for all  C3CBP's up to 40 vertices that all of them are $H^{+-}$, so $N^{+-} \geq 40$. For all  C3CBP's up to 44 vertices that do not contain subgraphs $R_2$ or $R_4$ they have checked that they are all $H^{+-}$, so $n\geq 44$. Now $N^*\geq \min(44, 4+40)=44$.

We will use the notation $G(R)G'$ to say that we reduce $G$ to $G'$ by reduction $R$. To generate all C3CBP's the following theorem is used:

\begin{thrm}
Let $G$ be a C3CBP of order greater then 8. Then, for some C3CBP $G'$ we have either $G(R_0)G'$ or $G(R_4)G'$. 
\end{thrm}

In other words:  we can get $G'$ from $G$ by using one of the reductions $R_0$ or $R_4$. This tells us that we can generate all C3CBP's by starting with the $C3CBP4$ in Figure \ref{c3cbp} and applying the reverses of reductions $R_0$ and $R_4$. We will prove this theorem in the next chapter.

Finally we want to give an idea about how many graphs Holton, Manvel and McKay have tested. On graphs with at most 40 vertices, there are 243.547 non-isomorphic Barnette graphs, from which 171.168 on exactly 40 vertices. Another 77.072 on 42 and 44 vertices that do not contain subgraph $R_2$ or $R_4$. The number of Barnette graphs grows exponentially in the number of vertices.

\chapter{Generating Barnette Graphs}
\label{Thetwooperations}
\begin{figure}[h]
\centering
\begin{minipage}{.49\textwidth}
\vspace{15pt}
\begin{center}
 \includegraphics[scale=.5]{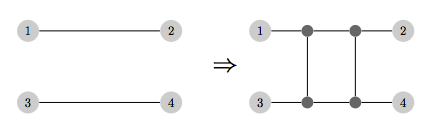}
\vspace{15pt}
\caption{Reversed reduction $R_0$}%
 \label{r0}%
\end{center}
\end{minipage}
\begin{minipage}{.49\textwidth}
\vspace{0pt}
\begin{center}
\includegraphics[scale=.5]{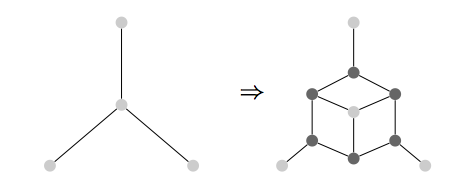}
\caption{Reversed reduction $R_4$}%
\label{r4}%
\end{center}
\end{minipage}
\end{figure}

Holton, Mavel and McKay \cite{HMM1985}  stated that we can generate every Barnette graph (from now on: C3CBP) by starting with the graph $C_1$ in Figure \ref{c1} and applying the reverses of the reductions $R_0$ and $R_4$ in Figure \ref{r0} and Figure \ref{r4}. In this chapter we will prove this theorem.

We use the following definition and Lemma's:
\begin{definition}  A graph is {\bf cyclically $k$-edge-connected} if at least $k$ edges must be removed to disconnect it into two components that each contain a cycle. \end{definition}

\begin{lemma} The smallest C3CBP has order 8, and graph $C_1$ in Figure \ref{c1} is the unique C3CBP on 8 vertices.  \end{lemma}
\begin{proof}
A C3CBP is cubic, thus every vertex has degree 3. As the sum of all degrees of all vertices must be even, the number of vertices of a C3CBP is even.
\begin{itemize}
\item A graph on two vertices cannot be cubic. 
\item The only cubic graph on four vertices is $K_4$, which is not bipartite.
\item The only bipartite cubic graph on 6 vertices is $K_{3,3}$, which is not planar.
\item There is only one graph on eight vertices that is cubic and bipartite, because it holds that every vertex in the upper set is not adjacent to exactly one vertex in the lower set. We get the graph in Figure \ref{c1flat}, which is $C_1$. According to the programs described in Chapter \ref{Programsinc}, this graph is also 3-connected. We can draw it as a planar graph as in Figure \ref{c1} and thus is is the unique C3CBP on 8 vertices.
\end{itemize}
Thus the smallest C3CBP has order 8, and graph $C_1$ in Figure \ref{c1} is the unique C3CBP on 8 vertices
\end{proof}


\begin{figure}[ht]
\begin{minipage}[b]{0.4\linewidth}
\centering
\begin{tikzpicture}
  [scale=.5,auto=left,every node/.style={circle, draw ,scale=0.5}]
  	\node[purple] (A) at (0,0) {1};
  	\node[blue] (B) at (6,0) {2};
  	\node[black] (C) at (6,6) {3};
  	\node[orange] (D) at (0,6) {4};
  	\node[black] (E) at (1.5,1.5) {5};
  	\node[orange] (F) at (4.5,1.5) {6};
  	\node[purple] (G) at (4.5,4.5) {7};
  	\node[blue] (H) at (1.5,4.5) {8};

  \foreach \from/\to in {A/B, B/C, C/D, D/A, E/F, F/G, G/H, H/E, A/E, B/F, C/G, D/H}
    	\draw (\from) -- (\to);
 
\end{tikzpicture}
\caption{Graph $C_1$}
\label{c1}
\end{minipage}
\hspace{0.1cm}
\begin{minipage}[b]{0.6\linewidth}
\centering
\begin{tikzpicture}
  [scale=2.5,auto=left,every node/.style={circle,draw,scale=0.5}]
  	\node[purple] (A) at (0,0) {1};
  	\node[blue] (B) at (1,0) {8};
  	\node[orange] (C) at (2,0) {6};
  	\node[black] (D) at (3,0) {3};
  	\node[purple] (E) at (0,1) {7};
  	\node[blue] (F) at (1,1) {2};
  	\node[orange] (G) at (2,1) {4};
  	\node[black] (H) at (3,1) {5};
  \foreach \from/\to in {A/F, A/G, A/H, B/E, B/G, B/H,C/E,C/F,C/H,D/E,D/F,D/G}
    	\draw (\from) -- (\to);
 
\end{tikzpicture}
\caption{Unique bipartite, cubic graph on 8 vertices.}
\label{c1flat}
\end{minipage}
\end{figure}

\begin{lemma} \label{sixfour} Let $G$ be a C3CBP of order at least 8. Then $G$ contains at least six 4-faces. \end{lemma}
\begin{proof}
According to the formula of Euler it holds that $F+V-E=2$, where $F$ is the number of faces, $V$ is the number of vertices and $E$ is the number of edges. $G$ is cubic, thus $E=\frac{3}{2}V$. 

We can compute the sum of all degrees $S$ in two ways. 
\begin{enumerate}
\item As $G$ is cubic, $S=3V$.
\item The sum of degrees $S$ equals $\sum_i i \cdot F_i$, where $F_i$ is the number of faces of size $i$.
\end{enumerate}

In a C3CBP all faces are even and of size $> 2$. We denote by $F_+$ the number of faces of order greater than 4, thus $F=F_4 + F_+$. Now it holds that:
\begin{eqnarray*}
3V & = & S\\
 & = & \sum_i i \cdot F_i\\
& = & 4F_4 + \sum_{i>4} i \cdot F_i \\ \label{line}
 & \geq & 4F_4 +6F_+ 
\end{eqnarray*}
Using Euler's formula and the equations above, we get:
\begin{eqnarray*}
F_4 + F_+ & = & F\\
 & = & 2-V+E\\
& = & 2+ \frac{1}{2}V\\
& \geq & 2 + \frac{1}{2} (\frac{1}{3} (4F_4 +6F_+))\\
& \geq & 2 + \frac{2}{3} F_4 + F_+
\end{eqnarray*}
And thus:
\begin{eqnarray*}
F_4 & \geq & 2 + \frac{2}{3} F_4 \\
\frac{1}{3} F_4 & \geq & 2 \\
F_4 & \geq & 6 \\
\end{eqnarray*}
Thus a C3CBP contains at least 6 4-faces.

\end{proof}

\begin{lemma} \label{3cut} Let $G$ be a C3CBP that has a 3-cut. When we split $G$ in $G_1$ and $G_2$ at this 3-cut as in Figure \ref{cut}, $G_1$ and $G_2$ are C3CBP. \end{lemma}
\begin{proof}
Without loss of generality, we only prove all properties for $G_1$. Let $a$ be the new vertex in $G_1$.
\begin{itemize}
\item {\bf Planarity}. One can check $G_1$ is planar from the construction and from the fact that $G$ is planar.
\item {\bf 3-Connectedness}. Suppose $G_1$ is not 3-connected, then there are two vertices $x_1\neq x_2$ that can be deleted to disconnect the graph.  One can check that whenever $x_1$ and $x_2$ would disconnect $G_1$, $G$ would not be 3-connected.
\item {\bf Cubic}. That $G_1$ is cubic is clear from the construction and from the fact that $G$ is cubic.
\item {\bf Bipartite}. We define some variables:
\begin{itemize}
\item $x_b$: the number of black vertices in $G_1 - {\{a\}}$ incident to the cut.
\item $x_w$: the number of white vertices in $G_1 - {\{a\}}$ incident to the cut.
\item $y_b$: the number of black vertices in $G_1 - {\{a\}}$ not incident to the cut.
\item $y_w$: the number of white vertices in $G_1 - {\{a\}}$ not incident to the cut.
\end{itemize}

A simple edge count in $G_1  \setminus{\{a\}}$ tells us that: $2x_b + 3y_b=2x_w + 3 y_w$. So $2x_b  \equiv 2x_w \mod 3$ and hence $x_b  \equiv  x_w \mod 3$. Thus the number of white vertices and black vertices in $G_1$ incident to the cut are the same modulo 3. Thus all vertices at one side of the cut are the same color. So we can keep the old coloring and without conflict color the last vertex $a$. Thus $G_1$ is bipartite.
\end{itemize}
Thus $G_1$ and $G_2$ are both C3CBP.
\end{proof}

\begin{lemma} \label{cyclicor3cut} Let $G$ be a C3CBP of order at least 8. Then $G$ is cyclically 4-edge connected or $G$ has a 3-cut. \end{lemma}
\begin{proof}
$G$ cannot have a 1-cut or 2-cut, because it is 3-connected. When $G$ has a 3-cut, we are done. So suppose $G$ does not contain a 3-cut. Then we know that $G$  will contain a 4-cut, as we can always separate two adjacent vertices in a cubic graph from all other vertices. We have to prove that $G$ has a 4-cut such that both parts contain a cycle. 

By Lemma \ref{sixfour} we know that $G$ has six 4-faces, and thus we can make a 4-cut around one of them. We denote that 4-face as $G_1$ and the rest of the graph as $G_2$. In the worst case we mess up four other 4-faces, but as we know there are at least six of them, at least one remains in $G_2$. Now $G_1$ and $G_2$ both contain a 4-face, and thus a cycle.

Thus every C3CBP has a 3-cut or is cyclically 4-edge connected.
\end{proof}

\begin{lemma} \label{cyclic} Let $G$ be a C3CBP of order greater than 8 that is cyclically 4-edge connected. Then at least one of the two possible applications of reduction $R_0$ to any 4-face produces a C3CBP \end{lemma}

\begin{figure}[h]%
            \centering
            \parbox{1.2in}{              \begin{tikzpicture}
  [scale=1.5,auto=left,every node/.style={circle,draw,fill=black!20,scale=0.5}]
  	\node (A) at (0,0) {3};
	\node (C) at (0,1) {1};
	\node (F) at (2,0) {4};
	\node (H) at (2,1) {2};
	\node[black!60]  (I) at (0.66,1){};
	\node[black!60]  (J) at (1.33,1){};
	\node[black!60]  (K) at (0.66,0){};
	\node[black!60]  (L) at (1.33,0){};
  \foreach \from/\to in {J/I,C/I,J/H,K/L,K/A,L/F,K/I,L/J}
    	\draw (\from) -- (\to);

	\end{tikzpicture}}
            \qquad
	\begin{minipage}{0.08in}
	$\Rightarrow$
	\end{minipage}
	\qquad
	\begin{minipage}{1.3in}%
                            \begin{tikzpicture}
  [scale=1.5,auto=left,every node/.style={circle,draw,fill=black!20,scale=0.5}]
  	\node (A) at (0,0) {3};
	\node (C) at (0,1) {1};
	\node (F) at (2,0) {4};
	\node (H) at (2,1) {2};
	\node[scale=0.01]  (I) at (0.66,1){};
	\node[scale=0.01]  (J) at (1.33,1){};
	\node[scale=0.01]  (K) at (0.66,0){};
	\node[scale=0.01]  (L) at (1.33,0){};
  \foreach \from/\to in {C/I,J/H,K/A,L/F,K/I,L/J}
    	\draw (\from) -- (\to);

	\end{tikzpicture}
            \end{minipage}%
	\qquad
	\begin{minipage}{0.08in}
	or
	\end{minipage}
	\qquad
            \begin{minipage}{1.2in}%
              \begin{tikzpicture}
  [scale=1.5,auto=left,every node/.style={circle,draw,fill=black!20,scale=0.5}]
  	\node (A) at (0,0) {3};
	\node (C) at (0,1) {1};
	\node (F) at (2,0) {4};
	\node (H) at (2,1) {2};
  \foreach \from/\to in {A/F,C/H}
    	\draw (\from) -- (\to);

	\end{tikzpicture}
            \end{minipage}%
            \caption{Two applications of reduction $R_0$.}%
            \label{applicationsr0}%
\end{figure}
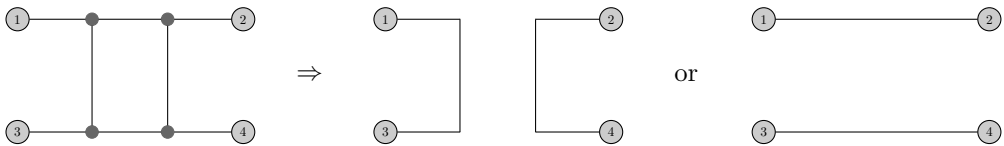
	
\begin{proof}
We have to prove that all properties of a C3CBP hold for $G'$, the graph obtained after the reduction.
\begin{itemize}
\item  {\bf Planarity}. This is trivial. By deleting edges a planar graph remains planar.
\item {\bf 3-connected} From Lemma 1 in \cite{Barnette1982} we know that at least one graph created by the one of the applications of the reduction $R_0$ will remain 3-connected.
\item  {\bf Cubic} From Figure \ref{applicationsr0} it is clear that vertices 1, 2, 3 and 4 still have degree 3 after reduction $R_0$. Thus $G'$ will remain cubic .
\item  {\bf Bipartite} As $G$ is bipartite, vertices 1 and 4 in Figure \ref{r0} have the same color color and vertices 2 and 3 another. In both possible reductions $R_0$ we only connect vertices with different colors, thus the old coloring is still correct. The C3CBP $G'$ obtained after the reduction remains bipartite.
\end{itemize} 
\end{proof}	
\begin{lemma} 
\label{minimalize} If a C3CBP, say graph $G$, has a 3-cut and we form $G_1$ and $G_2$ as in Figure \ref{cut} such that the number of vertices in $G_1$ is minimized. Then $G_1$ is cyclically 4-edge connected. \end{lemma}
\begin{proof} Suppose $G_1$ is not cyclically 4-edge connected. According to Lemma \ref{cyclicor3cut} $G_1$ must contain a 3-cut. But then $G_1$ was not minimized and we get a contradiction. Thus $G_1$ is cyclically 4-edge connected.
\end{proof}

We can now prove the following theorem:
\begin{thrm} Let $G$ be a C3CBP of order greater than 8. Then, there exists a C3CBP that we call $G'$, such that we have either $G(R_0)G'$ or $G(R_4)G'$. 
\label{twooperations}
\end{thrm}
	
\begin{figure}[h]
\centering
	 \begin{minipage}{0.4in}
             \begin{tikzpicture}
  [scale=0.8,auto=left,every node/.style={circle,fill=black!20,scale=0.5,draw}]
  	\node[{draw,orange!80}]  (A) at (1,1.3) {};
	\node[fill=orange!80]   (B) at (0.42,0.35) {};
	\node[fill=orange!80]  (C) at (1.58,0.35) {};
	\node (D) at (1,2) {};
	\node (E) at (0,0) {};
	\node  (F) at (2,0) {};
	\node[fill=violet!80]   (G) at (1,0.1) {};
	\node[fill=violet!80]   (H) at (0.42,1) {};
	\node[fill=violet!80]   (I) at (1.58,1) {};
	\node[fill=blue!80]  (J) at (1,0.75) {};
  \foreach \from/\to in {A/D,I/C,H/J,G/B}
    	\draw (\from) -- (\to);
	\draw (E) -- (B) -- (H) -- (A) --(I) -- (J) -- (G) -- (C) --(F);
	\end{tikzpicture}
	\end{minipage}
	 \qquad
	\begin{minipage}{0.2in}%
              $\Rightarrow$
            \end{minipage}%
            \begin{minipage}{1.2in}%
                           \begin{tikzpicture}
  [scale=0.8,auto=left,every node/.style={circle,draw,fill=black!20,scale=0.5}]
	\node (D) at (1,2) {};
	\node (E) at (0,0) {};
	\node  (F) at (2,0) {};
	\node[fill=blue!80]  (J) at (1,0.75) {};
  \foreach \from/\to in {D/J,E/J,F/J}
    	\draw (\from) -- (\to);
	\end{tikzpicture}            \end{minipage}%
		 \begin{minipage}{0.6in}
              \begin{tikzpicture}
  [scale=0.4,auto=left,every node/.style={circle,draw, fill=black!20 ,scale=0.5}]
  	\node[fill=orange!80]  (A) at (2,1) {};
	\node[fill=orange!80]   (B) at (3,3) {};
	\node[fill=orange!80]  (C) at (0,0) {};
	\node (D) at (3,0) {};
	\node (E) at (4,2) {};
	\node  (F) at (1,-1) {};
	\node[fill=violet!80]   (G) at (0,3) {};
	\node[fill=violet!80]   (H) at (2,2) {};
	\node[fill=violet!80]   (I) at (1,1) {};
	\node[fill=blue!80]  (J) at (1,2) {};
  \foreach \from/\to in {A/D,I/C,H/J,G/B}
    	\draw (\from) -- (\to);
	\draw  (E) -- (B) -- (H) -- (A) --(I) -- (J) -- (G) -- (C) --(F);
	\end{tikzpicture}
	\end{minipage}
	 \qquad
	\begin{minipage}{0.3in}%
              $\Rightarrow$
            \end{minipage}%
            \begin{minipage}{0.6in}%
              \begin{tikzpicture}
  [scale=0.5,auto=left,every node/.style={circle,draw, fill=black!20 ,scale=0.5}]
	\node (D) at (3,0) {};
	\node (E) at (4,2) {};
	\node  (F) at (1,-1) {};
	\node[fill=blue!80]  (J) at (1,2) {};
  \foreach \from/\to in {E/J,D/J,F/J}
    	\draw (\from) -- (\to);
	\end{tikzpicture}
            \end{minipage}%
	\caption{Reduction $R_4$, when a minimized $G_1$ by cutting at a 3-cut is $C_1$.}
	\label{r4used}
	\end{figure}
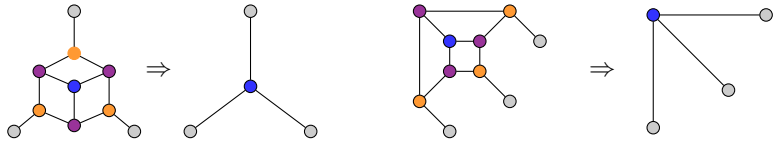

\begin{proof}
According to Lemma \ref{cyclicor3cut} there are two cases. Either $G$ is cyclically 4-edge connected or  it has a 3-cut.
\begin{enumerate}
\item Suppose $G$ is cyclically 4-edge connected. According to Lemma \ref{cyclic} we can now apply the reduction $R_0$.
\item Suppose $G$ has a 3-cut. Now form $G_1$ and $G_2$ as in Figure \ref{cut} such that $G_1$ is minimized. According to Lemma \ref{minimalize} $G_1$ is now cyclically 4-edge connected. If $G_1$ is the graph in Figure \ref{c1}, then we could apply reduction $R_4$ to $G$ as in Figure \ref{r4used}. If not, $G_1$ must have a 4-face which does not use any of the edges created by separating at the 3-cut (because a $G_1$ is a C3CBP according to Lemma \ref{3cut} and C3CBP contains at least 6 4-faces according to Lemma \ref{sixfour}) and we can apply at least one of the two possible applications of $R_0$ to one of the 4-faces by Lemma \ref{cyclic}.
\end{enumerate}
Thus for every C3CBP  $G$ of order greater than 8 there exists a $G'$ that is C3CBP such that we have either $G(R_0)G$' or $G(R_4)G'$.
\end{proof}

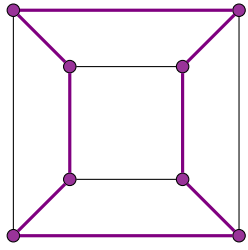
\begin{figure}
\centering
\begin{tikzpicture}
  [scale=.5,auto=left,every node/.style={circle,draw,fill=violet!80,scale=0.5}]
  	\node (A) at (0,0) {};
  	\node (B) at (6,0) {};
  	\node (C) at (6,6) {};
  	\node (D) at (0,6) {};
  	\node (E) at (1.5,1.5) {};
  	\node (F) at (4.5,1.5) {};
  	\node (G) at (4.5,4.5) {};
  	\node (H) at (1.5,4.5) {};

  \foreach \from/\to in {B/C, D/A, E/F, G/H}
    	\draw (\from) -- (\to);

\draw[very thick, violet] (A) -- (B) -- (F) --(G)--(C)--(D)--(H)--(E)--(A);

\end{tikzpicture}
\caption{Hamiltonian cycle in graph $C_1$}
\label{h}
\end{figure}

\begin{corollary} \label{gevolg1} The graph $C_1$ in Figure \ref{c1} is the only C3CBP on 8 or fewer vertices and all repeated reductions of other C3CBP eventually lead to $C_1$. And thus all Barnette graphs can be generated from $C_1$ by applying the reversed reductions $R_0$ and $R_4$ \end{corollary}


Now, if we can prove that $C_1$ contains a Hamiltonian cycle and the two operations do maintain Hamiltonicity, then we know according to Corollary \ref{gevolg1} that every C3CBP contains a Hamiltonian cycle and thus that Barnette's conjecture is true.

\begin{figure}[h]%
            \centering
	\begin{minipage}{0.4in}
           \begin{tikzpicture}
  [scale=0.7,auto=left,every node/.style={circle,draw,fill=black!20,scale=0.5}]
  	\node[fill=violet!80] (A) at (0,0) {};
	\node[fill=violet!80] (C) at (0,1) {};
	\node[fill=violet!80] (F) at (2,0) {};
	\node[fill=violet!80] (H) at (2,1) {};
  \foreach \from/\to in {A/F,C/H}
    	\draw[very thick,violet!80] (\from) -- (\to);
	\end{tikzpicture}
	\end{minipage}
            \qquad
	\begin{minipage}{0.2in}%
              $\Rightarrow$
            \end{minipage}%
            \begin{minipage}{0.8in}%
              \begin{tikzpicture}
  [scale=0.7,auto=left,every node/.style={circle,draw,fill=black!20,scale=0.5}]
  	\node[fill=violet!80]  (A) at (0,0) {};
	\node[fill=violet!80]  (C) at (0,1) {};
	\node[fill=violet!80]  (F) at (2,0) {};
	\node[fill=violet!80]  (H) at (2,1) {};
	\node[fill=violet!80]  (I) at (0.66,1){};
	\node[fill=violet!80]  (J) at (1.33,1){};
	\node[fill=violet!80]  (K) at (0.66,0){};
	\node[fill=violet!80]  (L) at (1.33,0){};
  \foreach \from/\to in {K/I,L/J}
    	\draw (\from) -- (\to);
 \foreach \from/\to in {J/I,C/I,J/H,K/L,K/A,L/F}
    	\draw[very thick,violet!80]  (\from) -- (\to);
	\end{tikzpicture}
            \end{minipage}%
	\qquad
	 \begin{minipage}{0.4in}
	\begin{tikzpicture}
  [scale=0.7,auto=left,every node/.style={circle,draw,fill=black!20,scale=0.5}]
  	\node (A) at (0,0) {};
	\node[fill=violet!80]  (C) at (0,1) {};
	\node (F) at (2,0) {};
	\node[fill=violet!80]  (H) at (2,1) {};
  \foreach \from/\to in {A/F,C/H}
    	\draw (A) -- (F);
	\draw[very thick,violet!80] (C) -- (H);
	\end{tikzpicture}
	\end{minipage}
	  \qquad
	\begin{minipage}{0.2in}%
              $\Rightarrow$
            \end{minipage}%
            \begin{minipage}{0.4in}%
             \begin{tikzpicture}
  [scale=0.7,auto=left,every node/.style={circle,draw,fill=black!20,scale=0.5}]
  	\node (A) at (0,0) {};
	\node[fill=violet!80]  (C) at (0,1) {};
	\node (F) at (2,0) {};
	\node[fill=violet!80]  (H) at (2,1) {};
	\node[fill=violet!80]   (I) at (0.66,1){};
	\node[fill=violet!80]  (J) at (1.33,1){};
	\node[fill=violet!80] (K) at (0.66,0){};
	\node[fill=violet!80]   (L) at (1.33,0){};
  \foreach \from/\to in {C/I,J/H,K/L,K/I,L/J}
    	\draw[very thick,violet!80]  (\from) -- (\to);
	\draw (K) -- (A);
	\draw (L) -- (F);
	\draw (I) -- (J);
	\end{tikzpicture}            \end{minipage}%
	\qquad 
	$$ $$
	 \begin{minipage}{0.4in}
	\begin{tikzpicture}
  [scale=0.7,auto=left,every node/.style={circle,draw,fill=black!20,scale=0.5}]
  	\node[fill=violet!80] (A) at (0,0) {};
	\node (C) at (0,1) {};
	\node[fill=violet!80] (F) at (2,0) {};
	\node (H) at (2,1) {};
  \foreach \from/\to in {A/F,C/H}
    	\draw[very thick,violet!80]  (A) -- (F);
	\draw (C) -- (H);
	\end{tikzpicture}
	\end{minipage}
	 \qquad
	\begin{minipage}{0.2in}%
              $\Rightarrow$
            \end{minipage}%
            \begin{minipage}{0.8in}%
              \begin{tikzpicture}
  [scale=0.7,auto=left,every node/.style={circle,draw,fill=black!20,scale=0.5}]
  	\node[fill=violet!80] (A) at (0,0) {};
	\node (C) at (0,1) {};
	\node[fill=violet!80] (F) at (2,0) {};
	\node (H) at (2,1) {};
	\node[fill=violet!80]  (I) at (0.66,1){};
	\node[fill=violet!80]  (J) at (1.33,1){};
	\node[fill=violet!80] (K) at (0.66,0){};
	\node[fill=violet!80]  (L) at (1.33,0){};
  \foreach \from/\to in {J/I,K/A,L/F,K/I,L/J}
    	\draw[very thick,violet!80]  (\from) -- (\to);

  \foreach \from/\to in {C/I,J/H,K/L}
    	\draw (\from) -- (\to);

	\end{tikzpicture}
	\end{minipage}%
	\qquad
	 \begin{minipage}{0.4in}
	\begin{tikzpicture}
  [scale=0.7,auto=left,every node/.style={circle,draw,fill=black!20,scale=0.5}]
  	\node (A) at (0,0) {};
	\node (C) at (0,1) {};
	\node (F) at (2,0) {};
	\node (H) at (2,1) {};
  \foreach \from/\to in {A/F,C/H}
    	\draw (\from) -- (\to);
	\end{tikzpicture}
	\end{minipage}
	 \qquad
	\begin{minipage}{0.2in}%
              $\Rightarrow$
            \end{minipage}%
            \begin{minipage}{0.4in}%
 	 ?
            \end{minipage}%
            \caption{Extended Hamiltonian cycle in $R_0$.}%
            \label{r0reduction}%
          \end{figure}
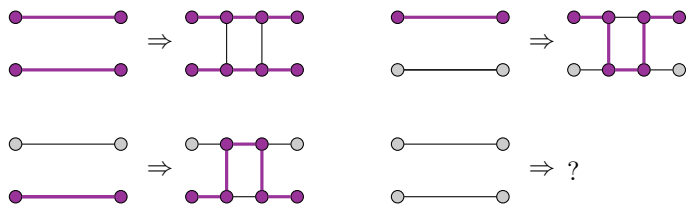

\begin{itemize}
\item
$C_1$ contains a Hamiltonian cycle. One is depicted in Figure \ref{h}.

\item Operation $R_4$ maintains Hamiltonicity. If we can find a Hamiltonian cycle in graph $G$ and we apply the reversed reduction $R_4$ to get graph $G'$, then we can find a Hamiltonian cycle in $G'$. We can use the Hamiltonian cycle found in $G$ and adapt it as in Figure \ref{h1}.


\begin{figure}[h]%
            \centering
	\begin{minipage}{0.4in}
            \begin{tikzpicture}
  [scale=0.8,auto=left,every node/.style={circle,draw,fill=black!20,scale=0.5}]
  	\node[fill=violet!80] (A) at (1,0.85) {};
	\node[fill=violet!80]  (C) at (1,2) {};
	\node[fill=violet!80]  (F) at (0,0) {};
	\node (H) at (2,0) {};

    	\draw[very thick, violet] (A) -- (C);
    	\draw[very thick, violet] (A) -- (F);
    	\draw (A) -- (H);

	\end{tikzpicture}%
	\end{minipage}
            \qquad
	\begin{minipage}{0.2in}%
              $\Rightarrow$
            \end{minipage}%
            \begin{minipage}{1in}%
              \begin{tikzpicture}
  [scale=0.8,auto=left,every node/.style={circle,draw,fill=black!20,scale=0.5}]
  	\node[fill=violet!80]   (A) at (1,1.3) {};
	\node[fill=violet!80]   (B) at (0.42,0.35) {};
	\node[fill=violet!80]  (C) at (1.58,0.35) {};
	\node[fill=violet!80]  (D) at (1,2) {};
	\node[fill=violet!80]   (E) at (0,0) {};
	\node  (F) at (2,0) {};
	\node[fill=violet!80]   (G) at (1,0.1) {};
	\node[fill=violet!80]  (H) at (0.42,1) {};
	\node[fill=violet!80]   (I) at (1.58,1) {};
	\node[fill=violet!80]  (J) at (1,0.75) {};
  \foreach \from/\to in {C/F,A/I,H/B,G/J}
    	\draw (\from) -- (\to);

	\draw[very thick,violet!80]  (D) -- (A) -- (H) -- (J) -- (I) -- (C) -- (G) -- (B) -- (E);

	\end{tikzpicture}
            \end{minipage}%
	\qquad
	 \begin{minipage}{0.4in}
	\begin{tikzpicture}
  [scale=0.8,auto=left,every node/.style={circle,draw,fill=black!20,scale=0.5}]
  	\node[fill=violet!80] (A) at (1,0.85) {};
	\node[fill=violet!80] (C) at (1,2) {};
	\node (F) at (0,0) {};
	\node[fill=violet!80] (H) at (2,0) {};
    	\draw[very thick,violet!80] (A) -- (C);
	\draw (A) -- (F);
	\draw[very thick,violet!80] (A) -- (H);
	\end{tikzpicture}
	\end{minipage}
	  \qquad
	\begin{minipage}{0.2in}%
              $\Rightarrow$
            \end{minipage}%
            \begin{minipage}{1in}%
              \begin{tikzpicture}
  [scale=0.8,auto=left,every node/.style={circle,draw,fill=black!20,scale=0.5}]
  	\node[fill=violet!80]  (A) at (1,1.3) {};
	\node[fill=violet!80]  (B) at (0.42,0.35) {};
	\node[fill=violet!80]   (C) at (1.58,0.35) {};
	\node[fill=violet!80]  (D) at (1,2) {};
	\node (E) at (0,0) {};
	\node[fill=violet!80]  (F) at (2,0) {};
	\node[fill=violet!80]  (G) at (1,0.1) {};
	\node[fill=violet!80]   (H) at (0.42,1) {};
	\node[fill=violet!80]   (I) at (1.58,1) {};
	\node[fill=violet!80]  (J) at (1,0.75) {};
  \foreach \from/\to in {B/E,A/H,I/C,G/J}
    	\draw (\from) -- (\to);
	\draw[very thick,violet!80] (D) -- (A) -- (I) -- (J) --(H) -- (B) -- (G) -- (C) -- (F);
	\end{tikzpicture}
            \end{minipage}%
	\qquad
	 \begin{minipage}{0.4in}
	\begin{tikzpicture}
  [scale=0.8,auto=left,every node/.style={circle,draw,fill=black!20,scale=0.5}]
  	\node[fill=violet!80]  (A) at (1,0.85) {};
	\node (C) at (1,2) {};
	\node[fill=violet!80]  (F) at (0,0) {};
	\node[fill=violet!80]  (H) at (2,0) {};
    	\draw (A) -- (C);
	\draw[very thick,violet!80] (A) -- (F);
	\draw[very thick,violet!80] (A) -- (H);

	\end{tikzpicture}
	\end{minipage}
	 \qquad
	\begin{minipage}{0.2in}%
              $\Rightarrow$
            \end{minipage}%
            \begin{minipage}{0.6in}%
              \begin{tikzpicture}
  [scale=0.8,auto=left,every node/.style={circle,draw,fill=black!20,scale=0.5}]
  	\node[fill=violet!80]  (A) at (1,1.3) {};
	\node[fill=violet!80]   (B) at (0.42,0.35) {};
	\node[fill=violet!80]  (C) at (1.58,0.35) {};
	\node (D) at (1,2) {};
	\node[fill=violet!80]  (E) at (0,0) {};
	\node[fill=violet!80]  (F) at (2,0) {};
	\node[fill=violet!80]   (G) at (1,0.1) {};
	\node[fill=violet!80]   (H) at (0.42,1) {};
	\node[fill=violet!80]   (I) at (1.58,1) {};
	\node[fill=violet!80]  (J) at (1,0.75) {};
  \foreach \from/\to in {A/D,I/C,H/J,G/B}
    	\draw (\from) -- (\to);
	\draw[very thick,violet!80]  (E) -- (B) -- (H) -- (A) --(I) -- (J) -- (G) -- (C) --(F);
	\end{tikzpicture}
            \end{minipage}%
            \caption{Extended Hamiltonian cycle in $R_4$.}%
            \label{h1}%
          \end{figure}
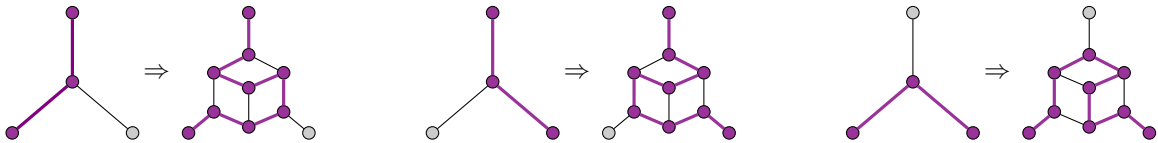


\item We cannot prove that the reversed reduction of $R_0$ maintains Hamiltonicity. Suppose we have a C3CBP $G$ with a Hamiltonian cycle and we apply the reversed reduction $R_0$ to get $G'$. Then there are four ways the Hamiltonian cycle could have passed the two edges used by reversed reduction $R_0$. In three of those cases we can extend a cycle in $G$ to a Hamiltonian cycle in $G'$ (see Figure \ref{r0reduction}). But in one case we cannot.
\end{itemize}


If we could prove that even in that last case a Hamiltonian cycle can be found, then Barnette's conjecture would be true. This question remains unsolved. 

\chapter{Programs in C++}
\label{Programsinc}
When researching Barnette's conjecture we wrote some programs in C++. They answer the following questions:
\begin{itemize}
\item Is the graph cubic?
\item Is the graph bipartite?
\item Is the graph 3-connected?
\item Does the graph contain a Hamiltonian cycle?
\item Does the graph contain a Hamiltonian path?
\end{itemize}

In this chapter we explain the functions used in the program and determine their complexity. The whole program can be found in Appendix \ref{c++}. The programs below are not optimal for their goals. Instead we tried to keep them simple and comprehensible. 

\section{Preliminaries}
Here we define the structures that are used in the program.

\begin{lstlisting}[frame=single]
typedef vector<int> Line;
typedef int Vertex;
typedef vector<Vertex> vertices;
typedef vector<Line> Lines;
struct Node
{
    Vertex Vertex;
    vector<Vertex> neighbors;
};
typedef vector<Node> Nodes;

\end{lstlisting}

An edge ({\it Line}) is a vector of two integers. We denote all vertices and edges in the vectors {\it Nodes} and {\it Lines}. Below you find an example about how to denote a graph.

\clearpage

\begin{lstlisting}[frame=single]
const int number_of_vertices=38;
const int number_of_lines=57;
int LINES [number_of_lines][2]=  
{{1,2},{1,4},{1,18},{2,3},{2,6},{3,8},{3,20},{4,5},{4,11},{5,6},{5,9},{6,7},
{7,8},{7,10},{8,14},{9,10},{9,11},{10,13},{11,12},{12,13},{12,15},{13,14},
{14,15},{15,16},{16,17},{16,19},{17,18},{17,22},{18,21},{19,20},{19,23},
{20,38},{21,22},{21,36},{22,23},{23,24},{24,25},{24,27},{25,26},{25,28},
{26,27},{26,30},{27,32},{28,29},{28,33},{29,30},{29,34},{30,31},{31,32},
{31,35},{32,38},{33,34},{33,36},{34,35},{35,37},{36,37},{37,38}};
\end{lstlisting}

In this specific case we defined the Barnette-Bos\'ak-Lederberg Graph (see Figure \ref{testcase}). From {\it LINES} we define {\it Lines} by adding for every edge $\{v,w\}$ in {\it LINES} the edges $\{v,w\}$ and $\{w,v\}$ in {\it Lines}. 

\section{Cubic}
Testing whether or not a graph is cubic can be done by counting the neighbors of each vertex. At the same time it stores those neighbors for further use (for example in the Hamiltonian cycle program). If there is a vertex that has not exactly three neighbors, the Boolean function returns false and immediately stops. In short, the algorithm is:
\begin{enumerate}[]
\item For every point $v \in V$:
\begin{enumerate}[]
\item If (not (degree($v$)=3)), return false;
\end{enumerate}
\item Return true;
\end{enumerate}

The whole algorithm can be found below:

\begin{lstlisting}[frame=single]
bool is_cubic(Lines lines, Nodes& nodes)
{   for(int i=0; i<number_of_vertices;i++)
    {   int degree=0;
        for(int j=0;j<static_cast<int>(lines.size());j++)
        {   if(lines[j][0]==i+1)
            {   nodes[i].neighbors.push_back(lines[j][1]);
                degree++;
            }
        }
        if(degree!=3)
        {   return false;}
    }
    return true;
}
\end{lstlisting}

When we have $n$ vertices and  $m$ lines, the number of iterations in the first for-loop is bounded by $O(n)$. There are $O(m)$ operations in each iteration. So the complexity of this algorithm is $O(nm)$. This means checking whether or not a graph is cubic can be done in polynomial time.

An optimal program for this problem has complexity $O(m)$. 
 
\clearpage

\section{Bipartite}

Checking whether a graph is bipartite consists of two parts:
\begin{enumerate}
\item {\bf Find a possible coloring}. To find a possible coloring, we define two sets: red and blue. We color one of the vertices red, and then we are going to iterate. The iterations consist of finding all uncolored neighbors to vertices in the red set, and color them blue. Then we are going to look for all uncolored neighbors to vertices in the blue set, and color them red. As in each iteration at least one point gets colored, the number of iterations is bounded by $O(n)$. In every iteration we look at all points in a set ($O(n)$) and check whether or not it has uncolored neighbors ($O(nm)$). Thus we have a total complexity of $O(n^2m)$.
\item {\bf Validate the coloring}. For every edge ($O(m)$) we check if both endpoints have different colors ($O(n)$), and this has a complexity of $O(nm)$.
\end{enumerate}
A total complexity of this function is thus $O(n^2m+nm)=O(n^2m)$. So testing whether or not a graph is bipartite can be done in polynomial time. An optimal algorithm can check whether or not a graph is bipartite within $O(n \log n)$.\\

\begin{lstlisting}[frame=single]
bool is_bipartite(Lines lines, Nodes& blue, Nodes& red, Nodes nodes)
{   Nodes uncolored;
    uncolored=nodes;
    red.push_back(uncolored[static_cast<int>(uncolored.size())-1]);
    uncolored.pop_back();
    while(static_cast<int>(uncolored.size())>0)
    {   for(int i=0;i<static_cast<int>(uncolored.size());i++)
        {   for(int j=0;j<static_cast<int>(red.size());j++)
            {   if(exist_path(uncolored[i],red[j].vertex,lines))
                {   blue.push_back(uncolored[i]);
                    move_back(i,uncolored);
                    i--;
                }
            }
        }
        for(int i=0;i<static_cast<int>(uncolored.size());i++)
        {   for(int j=0;j<static_cast<int>(blue.size());j++)
            {   if(exist_path(uncolored[i],blue[j].vertex,lines))
                {   red.push_back(uncolored[i]);
                    move_back(i,uncolored);
                    i--;
                }
            }
        }
    }
    for(int i=0;i<2*number_of_lines;i=i+2)
    {   if((is_red(lines[i][0],red) && is_red(lines[i][1],red)) || 
           (is_blue(lines[i][0],blue) && is_blue(lines[i][1],blue)) )
        {   cout<<"Error found at "<<lines[i][0]<< " - "<<lines[i][1]<<endl;
            return false; 
        }
    }
    return true;
}
\end{lstlisting}

\section{3-Connected}

A graph is called 3-connected when you have to delete at least three vertices to disconnect the graph. So a way to test whether a graph is 3-connected is to delete all combinations of two vertices end check whether the graph remains connected. Now every iteration consists of two steps:
\begin{enumerate}
\item Deleting two vertices ($O(n)$) and the adjacent edges ($O(m)$), which has a total complexity of $O(n+m)$.
\item Checking whether the remaining graph remains connected. This we do by making two sets: colored and uncolored. We delete the last vertex of the uncolored set (vector) and we add it to the colored one. Now we iterate: in every iteration we color all uncolored neighbors of vertices in the colored set ($O(nm)))$. When the graph is connected, we can add at least one vertex to the colored set in every iteration, so when the graph is connected, the number of iterations needed to color every vertex is bounded by $n$. If not every vertex is colored by then, the graph is not connected. Thus we keep a counter for the number of iterations. When the graph is disconnected, the counter will exceed $n$ and the function stops and returns false. This gives a total complexity of $O(n^2m)$. 
\end{enumerate} 

As the possible combinations of two vertices is bounded by $O(n^2)$, the total complexity of this function is $O(n^2(n+m+n^2m))=O(n^4m)$.\\

The core of the program, that given a graph, tests whether the graph is connected, is denoted below.

\begin{lstlisting}[frame=single]
bool connected(Nodes nodes, Lines lines)
{   Nodes uncolored=nodes;
    Nodes colored;
    colored.push_back(uncolored[static_cast<int>(uncolored.size())-1]);
    uncolored.pop_back();
    int count=1;
    while(count<number_of_vertices && static_cast<int>(uncolored.size())>0)
    {   count++;
        for(int i=0;i<static_cast<int>(uncolored.size());i++)
        {   for(int j=0;j<static_cast<int>(colored.size());j++)
            {   if(exist_path(uncolored[i],colored[j].vertex,lines) 
			&& i>=0 && static_cast<int>(uncolored.size())>i)
                {   colored.push_back(uncolored[i]);
                    move_back(i,uncolored,hamiltonian_cycle);
                }
            }
        }
    }
    if(count<number_of_vertices-1)
    {   return true;   }
    return false;
}
\end{lstlisting}

There are better algorithms nowadays. Checking whether or not a graph is 3-connected can be done in $O(mn^{\frac{2}{3}})$ according to Algorithm 9 in \cite{Abdol}.

\section{Hamiltonian cycle}

To find a Hamiltonian cycle we use a recursive depth-first search algorithm. The function works on a vertex (the vertex we will look at in that step), the vertices and the edges of the graph, a path that contains the vertices we have passed (and in which order), and finally the vertex where we started in the first place. This function has some base cases and a recursive part:
\begin{itemize}
\item The base cases are:
\begin{enumerate}
\item We are at the vertex where we started and we have visited all vertices. In this case we have found a Hamiltonian cycle and we return the value true.
\item We are at a vertex that we have visited before, but it is not the previous case. Now we are at a dead end, so we return the value false.
\end{enumerate}
\item In the recursive part we are at a vertex that we have not visited before. We add the vertex to the visited vertices and call the function again, but now with each of the adjacent vertices. 
\end {itemize}

This algorithm results in a depth-first search. When it has found a Hamiltonian cycle, it stops. When it is stuck in some path, it goes back to the previous vertex and tries another possibility. 

In the worst case, when there is no Hamiltonian cycle, the whole search tree will be searched. As every leaf in the tree will lead to 3 new leaves, the number of nodes in the search tree will be in $O(n\cdot 3^n)$. It is clear that this algorithm does not run in polynomial time. (But we expected this, because finding a Hamiltonian cycle is even NP-complete.) An optimal algorithm works in $O(2^n)$. \\

\begin{lstlisting}[frame=single]
bool contains_hamiltonian_cycle(Node node,Nodes nodes, Lines lines, 
				Nodes& already_visited, int start)
{   // Base cases:
    if (node.vertex==start && visited_all(nodes, already_visited))
    {   return true;}
    else
    {   if (visited(node,already_visited))
        {   return false;}
    //Recursion
        else
        {   already_visited.push_back(node);
            for(int i=0;i<static_cast<int>(node.neighbors.size());i++)
            {   if(contains_hamiltonian_cycle(nodes[node.neighbors[i]-1],
				nodes, lines, already_visited, start))
                {   return true;}
            }
            already_visited.pop_back();
            return false;
        }
    }
}
\end{lstlisting}

\clearpage
 
\section{Hamiltonian path}
The algorithm to find a Hamiltonian path is almost equal to finding a Hamiltonian cycle. Only two changes are made:

\begin{enumerate}
\item In the first base case, we only check whether or not we visited all vertices. So we do not check whether or not we are at the starting vertex. This new function has the name {\it `contains\_hamiltonian\_path1'}.
\item We have a for-loop over all possible starting vertices. We cannot take an arbitrary vertex as starting vertex.
\end{enumerate}

\begin{lstlisting}[frame=single]
bool contains_hamiltonian_path(Nodes nodes, Lines lines, 
						Nodes& already_visited)
{   for(int start=0; start<number_of_vertices; start++)
    {   if(contains_hamiltonian_path1(nodes[start],nodes,lines,
				already_visited, nodes[start].vertex))
        {   return true;}
    }
    return false;
}
\end{lstlisting}

\section{Results}

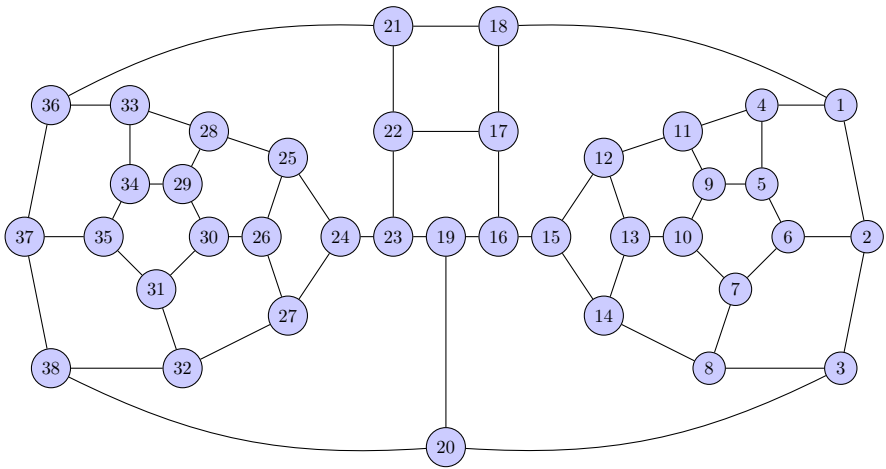
\begin{figure}[h]
\centering
\begin{tikzpicture}
  [scale=.35,auto=left,every node/.style={circle,fill=blue!20,scale=0.7,draw}]
  	\node (A1) at (4,9) {35};
  	\node (B1) at (6,7) {31};
  	\node (C1) at (8,9) {30};
  	\node (D1) at (7,11) {29};
  	\node (E1) at (5,11) {34};
  	\node (G1) at (7,4) {32};
  	\node (I1) at (8,13) {28};
  	\node (J1) at (5,14) {33};
  	\node (K1) at (2,4) {38};
  	\node (L1) at (1,9) {37};
  	\node (M1) at (2,14) {36};
  	\node (N1) at (11,6) {27};
  	\node (O1) at (10,9) {26};
  	\node (P1) at (11,12) {25};
  	\node (Q1) at (13,9) {24};
  	\node (R1) at (15,9) {23};

  \foreach \from/\to in {A1/B1, B1/C1, C1/D1, D1/E1, E1/A1, A1/L1, B1/G1, C1/O1, D1/I1, E1/J1, O1/P1, O1/N1, M1/J1, J1/I1, I1/P1, P1/Q1, Q1/N1, N1/G1, G1/K1, K1/L1, L1/M1, Q1/R1}
    	\draw (\from) -- (\to);

	\node (A2) at (30,9) {6};
  	\node (B2) at (28,7) {7};
  	\node (C2) at (26,9) {10};
  	\node (D2) at (27,11) {9};
  	\node (E2) at (29,11) {5};
  	\node (G2) at (27,4) {8};
  	\node (I2) at (26,13) {11};
  	\node (J2) at (29,14) {4};
  	\node (K2) at (32,4) {3};
  	\node (L2) at (33,9) {2};
  	\node (M2) at (32,14) {1};
  	\node (N2) at (23,6) {14};
  	\node (O2) at (24,9) {13};
  	\node (P2) at (23,12) {12};
  	\node (Q2) at (21,9) {15};
  	\node (R2) at (19,9) {16};

  \foreach \from/\to in {A2/B2, B2/C2, C2/D2, D2/E2, E2/A2, A2/L2, B2/G2, C2/O2, D2/I2, E2/J2, O2/P2, O2/N2, M2/J2, J2/I2, I2/P2, P2/Q2, Q2/N2, N2/G2, G2/K2, K2/L2, L2/M2, Q2/R2}
    	\draw (\from) -- (\to);

  	\node (X1) at (15,13) {22};
  	\node (X2) at (19,13) {17};
  	\node (Y1) at (15,17) {21};
  	\node (Y2) at (19,17) {18};
  	\node (Z1) at (17,9) {19};
  	\node (Z2) at (17,1) {20};

  \foreach \from/\to in {R1/Z1, R1/X1, X1/Y1, Y1/Y2, X2/Y2, R2/X2, Z1/R2, Z1/Z2, X1/X2}
    	\draw (\from) -- (\to);

  	\draw (M1) to [bend left=15] (Y1); 
  	\draw (M2) to [bend right=15] (Y2); 
  	\draw (K1) to [bend right=15] (Z2); 
  	\draw (K2) to [bend left=15] (Z2); 
     
\end{tikzpicture}
\caption{Barnette-Bos\'ak-Lederberg Graph.}
\label{testcase}
\end{figure}

We tried the program on the smallest counterexample for Tait's conjecture, the Barnette-Bos\'ak-Lederberg graph (Figure \ref{testcase}).

\clearpage

The result is copied below:
\begin{lstlisting}[frame=single]
Graph is cubic

Proposed coloring:
Red : 38 - 3 - 19 - 36 - 35 - 31 - 27 - 1 - 18 - 17 - 22 - 15 - 14 -
      28 - 29 - 6 - 7 - 25 - 9 - 11.
Blue: 20 - 37 - 32 - 2 - 34 - 8 - 26 - 16 - 24 - 23 - 33 - 21 - 30 - 
      13 - 12 - 10 - 5 - 4.
Validate coloring
Error found at 1 - 18
Graph is not bipartite

Graph is 3-connected

Hamiltonian path found:
1 - 4 - 5 - 6 - 2 - 3 - 8 - 7 - 10 - 9 - 11 - 12 - 13 - 14 - 15 - 16 - 
17 - 18 - 21 - 22 - 23 - 19 - 20 - 38 - 32 - 27 - 24 - 25 - 26 - 30 -
31 - 35 - 34 - 29 - 28 - 33 - 36 - 37.

No Hamiltonian cycle found
\end{lstlisting}

As the graph is  a Tait-graph, we know the graph is planar (not checked, but clearly true when we look at the graph), cubic and 3-connected. According to our program it is not bipartite and there is a Hamiltonian path, but no Hamiltonian cycle could be found.

Speed:
\begin{itemize}
\item Cubic: $<$ 1 second
\item Bipartite: $<$ 1 second
\item 3-Connected: $<$ 10 seconds
\item Hamiltonian path: 38 minutes
\item Hamiltonian cycle: 2 hours
\end{itemize}

For Barnette's conjecture we also have to check planarity, but that program is more complicated than the ones described above. There exist programs that decide in polynomial time whether or not a graph is planar, but we will not describe such a program here.

When Barnette's conjecture is true, it will take less time to check whether a graph contains a Hamiltonian cycle. If the graph is bipartite, cubic, 3-connected and planar (which can all be checked in polynomial time), the graph must contain a Hamiltonian cycle.

\chapter{Constructing a couterexample}
\label{Counterexample}
Getting closer to the end of our project we had the feeling it would be usefull to think about possible ways to construct a counterexample to Barnette's conjecture. As a starting point we considered the counterexamples for Tait's and Tutte's conjecture. The main ingredient for these counterexamples were required edge fragments, so we tried to find a required edge fragment that is cubic, bipartite, planar and 3-connected.  These fragments can then be, dependent on the coloring, connected as in Figure \ref{connectedfragments} or Figure \ref{label}. It is clear that these graphs can not contain a Hamiltonian cycle.  

\begin{figure}[h]
\begin{minipage}{3in}
\begin{center}
\begin{tikzpicture}
  [scale=1,auto=left,every node/.style={circle,draw, fill=white, scale=1}]
  	\node (R) at (0,0) {};
  	\node (Q1) at (0,-1.2) {$\uparrow$};
  	\node (Q2) at (1,1) {$\swarrow$};
  	\node (Q3) at (-1,1) {$\searrow$};

  \foreach \from/\to in {Q1/R, Q2/R, Q3/R}
    	\draw (\from) -- (\to);

	\draw (Q1) to [bend right=45] (Q2); 
  	\draw (Q2) to [bend right=30] (Q3); 
	\draw (Q3) to [bend right=45] (Q1); 

\end{tikzpicture}
\caption{Three required edge fragments connected like in Chapter \ref{Tait}.}
\label{connectedfragments}
\end{center}
\end{minipage}
\qquad
\begin{minipage}{3in}
\begin{center}
\begin{tikzpicture}
  [scale=0.4,auto=left,every node/.style={circle, draw, fill=white,  scale=.7}]

\draw plot [smooth] coordinates { (4,6) (4.5,6.5) (7,4) (8,1) (7,1)};
\draw plot [smooth] coordinates { (3,1) (3.3,1) (4,1.5) (4.6,2) (4.9,2)};
\draw plot [smooth] coordinates { (1,1)(0.7,0.5) (3.5,0) (6,0.2) (6,0.7)};

  	\node (aP1) [black] at (4,6) {};
  	\node (aF2) [black] at (3,4) {};
  	\node (aB2) [black] at (5,4) {};

  	\node (bP1) at (2,3) {};
  	\node (bF2) at (1,1) {};
  	\node (bB2) at (3,1) {};

  	\node (cP1) at (6,3) {};
  	\node (cF2) at (5,1) {};
  	\node (cB2) at (7,1) {};

  	\node (cQ) [black] at (4.9,2) {};
  	\node (cR) [black] at (6,0.7) {};

 	\draw (aP1) to [bend right=50] (aF2); 
 	\draw (aF2) to [bend right=30] (aB2); 
 	\draw (aB2) to [bend right=50] (aP1); 

	\draw (bP1) to [bend right=50] (bF2); 
 	\draw (bF2) to [bend right=30] (bB2); 
 	\draw (bB2) to [bend right=50] (bP1); 

	\draw (cB2) to [bend right=50] (cP1); 

 	\draw (cQ) to [bend right=15] (cF2); 
	\draw (cF2) to [bend right=15] (cR); 

 	\draw (cQ) to [bend left=25] (cP1); 
	\draw (cB2) to [bend left=15] (cR); 

 	\draw (aF2) to (bP1); 
 	\draw (aB2) to (cP1); 

\draw plot [smooth] coordinates {(4.5,1.5)(5.5,0.5)};
\draw plot [smooth] coordinates {(4.5,0.5)(5.5,1.5)};

\end{tikzpicture}
\end{center}
\caption{Three required edge fragments connected like in Chapter \ref{Tutte}.}
\label{label}
\end{minipage}
\end{figure}

\begin{itemize}
\item For Tait's conjecture a counterexample was made by looking at the properties of a pentogonal prism. We described in Chapter \ref{Tait} why some edges can not be used in a Hamiltonian cycle. This argument can not be used for Barnette's conjecture, because  a pentogonal prism contains a 5-cycle, and thus is not bipartite. Trying to copy this construction with a cycle of even length as its base does not work, as it loses the property that made the construction work in the first place.
\item For Tutte's conjecture a counterexample was made with Horton fragments. We thought it would not be usefull to search for circles with fewer vertices, because Horton would have found them already. Circles with more vertices will keep the $K_{3,3}$ as a minor, and thus will not be planar.
\end{itemize}


We tried some other possible starting points, using the operations treated in Chapter \ref{Thetwooperations}. In this search we formulated a condition that would imply Barnette's conjecture: \emph{The fourth case of  Figure \ref{r0reduction} can be solved}. If we are able to prove that in this fourth case Hamiltonian cycles are also preserved after applying the reversed reductions, we will be able to prove Barnette's conjecture. One can prove that all Barnette graphs can be generated by starting with graph $C_1$ in Figure \ref{c1} and applying the reversed reductions $R_0$ and $R_4$ given in Figure \ref{r0} and Figure \ref{r4}. So if Hamiltonian cycles are indeed preserved after these reversed reductions, all Barnette graphs contain a Hamiltonian cycle, because it is clear that graph $C_1$ (Figure \ref{c1}) contains a Hamiltonian cycle.

After formulating this condition that would imply Barnette's conjecture, we thought about a way of proving this condition. For this way of thinking we use the following definitions on graphs.

\begin{definition}
A {\bf matching} $M$ in a graph $G$ is a set of edges in which no two different edges share a common vertex.
\end{definition}

\begin{definition}
A vertex is {\bf matched} by a matching $M$ if it is an endpoint of one of the edges in $M$.
\end{definition}

\begin{definition}
A {\bf perfect matching} is a matching with the property that all vertices of the graph are matched.
\end{definition}

All Barnette graphs are cubic, so whenever we have found one Hamiltonian cycle and remove the edges of this cycle, every vertex has degree one. Thus the remaining edges are a perfect matching of the graph. The two edges given in the fourth case of Figure \ref{r0reduction} are not part of the intitial Hamiltonian cycle, so those two edges are part of the corresponding matching. The question is whether or not we can make a second Hamiltonian cycle using all edges of the matching. If this would be possible, then we can solve the fourth case of Figure \ref{r0reduction}, because we can always change this case into the first case of Figure \ref{r0reduction} in which both edges are used in the Hamiltonian cycle. By handling the fourth case like this, we see that Hamiltonicity is preserved after the reversed reductions. Possibly we can use the following theorem, which is stated as Theorem 1 in \cite{chiaong}, to prove this.

\begin{thrm}
Every edge of a cubic graph lies on an even number of Hamiltonian cycles. Consequently a cubic Hamiltonian graph has at least three Hamiltonian cycles.
\end{thrm}

This theorem tells us something about the number of Hamiltonian cycles in a Barnette graph, because all Barnette graphs are cubic. We were not able to prove that the matching that is left after removing the edges of one Hamiltonian cycle is contained in a second Hamiltonian cycle. We did not succeed in finding a counterexample and neither we got an idea how to prove Barnette's conjecture.

\bibliography{_Graphtheory}
\bibliographystyle{plain}

\appendix

\chapter{Code of the programms in C++}

\begin{lstlisting}

#include <iostream>
#include <vector>
#include <cassert>

using namespace std;

// Variables
typedef vector<int> Line;
typedef int Vertice;
typedef vector<Vertice> Vertices;
typedef vector<Line> Lines;
struct Node
{
    Vertice vertice;
    vector<Vertice> neighboors;
};
typedef vector<Node> Nodes;

// Graph
const int number_of_vertices=38;
const int number_of_lines=57;
int LINES [number_of_lines][2]=  {{1,2},{1,4},{1,18},{2,3},{2,6},
{3,8},{3,20},{4,5},{4,11},{5,6},{5,9},{6,7},{7,8},{7,10},{8,14},
{9,10},{9,11},{10,13},{11,12},{12,13},{12,15},{13,14},{14,15},
{15,16},{16,17},{16,19},{17,18},{17,22},{18,21},{19,20},{19,23},
{20,38},{21,22},{21,36},{22,23},{23,24},{24,25},{24,27},{25,26},
{25,28},{26,27},{26,30},{27,32},{28,29},{28,33},{29,30},{29,34},
{30,31},{31,32},{31,35},{32,38},{33,34},{33,36},{34,35},{35,37},
{36,37},{37,38}};

/////// HAMILTONCYCLE /////////////

bool visited(Node node, Nodes already_visited)
{
    //Precondition:
    assert(static_cast<int>(already_visited.size())>=0);
    //Postcodition: Determines if the vertice "Vertice" is aready visitid before
    for(int i=0;i<static_cast<int>(already_visited.size());i++)
    {
        if (already_visited[i].vertice==node.vertice)
        {
            return true;
        }
    }
    return false;
}

bool visited_all(Nodes nodes, Nodes already_visited)
{
    //Precondition:
    assert(static_cast<int>(already_visited.size())>=0 && static_cast<int>
							(nodes.size())>=0);
    //Postcodition: Determines if all the vertices in "vertices" are visited
    if(static_cast<int>(already_visited.size())==0)
    {
        return false;
    }
    for(int i=0;i<number_of_vertices;i++)
    {
        if (!visited(nodes[i],already_visited))
        {
            return false;
        }
    }
    return true;
}


bool exist_path(Node node, int vertice2, Lines lines)
{
    //Precondition:
    assert(static_cast<int>(lines.size())>=0);
    assert(static_cast<int>(node.neighboors.size())>=0);
    //Postcodition: Determines if the exist an edge between punt1 and punt2.
    for(int i=0; i<static_cast<int>(node.neighboors.size());i++)
    {
        if(node.neighboors[i]==vertice2)
        {
            return true;
        }
    }
    return false;
}

void print (Nodes nodes)
{
    //Precondition:
    assert(static_cast<int>(nodes.size())>=0);
    //Postcodition: Determines if the vertice "punt" is aready visitid before
    if(static_cast<int>(nodes.size())!=0)
    {
        cout<<nodes[0].vertice;
        for (int i=1;i<static_cast<int>(nodes.size());i++)
        {
            cout<<" - "<<nodes[i].vertice;
        }
        cout<<"."<<endl;
    }
}

bool contains_hamiltoncycle(Node node,Nodes nodes, Lines lines, Nodes& 
						already_visited, int start)
{
    //Precondition:
    assert(static_cast<int>(nodes.size())>=0);
    assert(static_cast<int>(lines.size())>=0);
    assert(static_cast<int>(already_visited.size())>=0);
    assert(start<=number_of_vertices);
    //Postcodition: Looks for a hamiltoncycle. Recursive & deptfirst. 
		Given a sequence of vertices, tries all possible vertices.
    //// Base cases
    if (node.vertice==start && visited_all(nodes, already_visited))
    {
        return true;
    }
    else
    {
        if (visited(node,already_visited))
        {
            return false;
        }
    //Recursion
        else
        {
            already_visited.push_back(node);
            for(int i=0;i<static_cast<int>(node.neighboors.size());i++)
            {

                if(contains_hamiltoncycle(nodes[node.neighboors[i]-1],
			nodes, lines, already_visited, start)==true)
                {
                    return true;
                }
            }
            already_visited.pop_back();
            return false;
        }
    }
}

/////// HAMILTONCYCLE /////////////

bool contains_hamiltonpath1(Node node,Nodes nodes, Lines lines, Nodes& 
						already_visited, int start)
{
    //Precondition:
    assert(static_cast<int>(nodes.size())>=0);
    assert(static_cast<int>(lines.size())>=0);
    assert(static_cast<int>(already_visited.size())>=0);
    //Postcodition: Looks for a path. Recursive & deptfirst. 
	Given a sequence of vertices, tries all possible vertices.
    //// base cases
    already_visited.push_back(node);
    if (visited_all(nodes, already_visited))
    {
        return true;
    }
    else
    {
        already_visited.pop_back();
        if (visited(node,already_visited))
        {
            return false;
        }
    //Recursion
        else
        {
            already_visited.push_back(node);
            for(int i=0;i<static_cast<int>(node.neighboors.size());i++)
            {
                if(contains_hamiltoncycle(nodes[node.neighboors[i]-1],nodes, 
					lines, already_visited, start)==true)
                {
                    return true;
                }
            }
            already_visited.pop_back();
            return false;
        }
    }
}

//////////////// Cubic //////////////////////
bool is_cubic(Lines lines, Nodes& nodes)
{
    //Precondition:
    assert(static_cast<int>(lines.size())>=0);
    //Postcodition: Determines if the graph is cubic
    for(int i=0; i<number_of_vertices;i++)
    {
        int degree=0;
        for(int j=0;j<static_cast<int>(lines.size());j++)
        {
            if(lines[j][0]==i+1)
            {
                nodes[i].neighboors.push_back(lines[j][1]);
                degree++;
            }
        }
        if(degree!=3)
        {
            return false;
        }
    }
    return true;
}

/////////////////////////// Bipartite  //////////////////////

void move_back(int i, Nodes& uncoloured)
{
    //Precondition:
    assert(static_cast<int>(uncoloured.size())>=0);
    assert(i<static_cast<int>(uncoloured.size()));
    assert(i>=0);
    //Postcodition: Switches the element on position i with the end of the vector
    int j=static_cast<int>(uncoloured.size())-1;
    Node temp=uncoloured[i];
    uncoloured[i]=uncoloured[j];
    uncoloured[j]=temp;
    uncoloured.pop_back();

}

bool is_red(int vertice, Nodes red)
{
    //Precondition:
    assert(static_cast<int>(red.size())>=0);
    //Postcodition: Is the vertice coloured red?
    for(int i=0;i<static_cast<int>(red.size());i++)
    {
        if (vertice==red[i].vertice)
        {
            return true;
        }
    }
    return false;
}

bool is_blue(int vertice, Nodes blue)
{
    //Precondition:
    assert(static_cast<int>(blue.size())>=0);
    //Postcodition: Is the vertice coloured blue?
    for(int i=0;i<static_cast<int>(blue.size());i++)
    {
        if (vertice==blue[i].vertice)
        {
            return true;
        }
    }
    return false;
}

bool is_bipartite(Lines lines, Nodes& blue, Nodes& red, Nodes nodes)
{
    Nodes uncoloured;
    uncoloured=nodes;
    red.push_back(uncoloured[static_cast<int>(uncoloured.size())-1]);
    uncoloured.pop_back();
    while(static_cast<int>(uncoloured.size())>0)
    {
        ////// Blauw aanvullen ///////
        for(int i=0;i<static_cast<int>(uncoloured.size());i++)
        {
            for(int j=0;j<static_cast<int>(red.size());j++)
            {
                if(exist_path(uncoloured[i],red[j].vertice,lines))
                {
                    blue.push_back(uncoloured[i]);
                    move_back(i,uncoloured);
                    i--;
                }
            }
        }
        /////// Rood aanvullen /////////
        for(int i=0;i<static_cast<int>(uncoloured.size());i++)
        {
            for(int j=0;j<static_cast<int>(blue.size());j++)
            {
                if(exist_path(uncoloured[i],blue[j].vertice,lines))
                {
                    red.push_back(uncoloured[i]);
                    move_back(i,uncoloured);
                    i--;
                }
            }
        }
    }
    //// Printen //////
    cout<<"Proposed colouring:"<< endl;
    cout<< "Red : ";
    print(red);
    cout<< "Blue: ";
    print(blue);
    cout<<"Vadidate colouring"<<endl;
    ////// Validate coloring
    for(int i=0;i<2*number_of_lines;i=i+2)
    {
        if((is_red(lines[i][0],red) && is_red(lines[i][1],red)) || 
		(is_blue(lines[i][0],blue) && is_blue(lines[i][1],blue)) )
        {
            cout<<"Error found at "<<lines[i][0]<< " - "<<lines[i][1]<<endl;
            return false; 
        }
    }
    return true;
}

    ////// 3-Connected  /////////////////////////

bool connected(Nodes nodes, Lines lines)
{
    //Precondition:
    assert(static_cast<int>(nodes.size())>=0);
    assert(static_cast<int>(lines.size())>=0);
    //Postcodition: Determines is the graph is connected. 'colors' one point, 
	and then adds all adjacent points. If there are uncolored points after 
	number_of_vertices-1 steps, then the graps is disconnected.
    Nodes uncoloured=nodes;
    Nodes colored;
    colored.push_back(uncoloured[static_cast<int>(uncoloured.size())-1]);
    uncoloured.pop_back();
    int count=1;
    while(count<number_of_vertices && static_cast<int>(uncoloured.size())>0)
    {
        count++;
        for(int i=0;i<static_cast<int>(uncoloured.size());i++)
        {
            for(int j=0;j<static_cast<int>(colored.size());j++)
            {
                if(exist_path(uncoloured[i],colored[j].vertice,lines) && i>=0 
					&& static_cast<int>(uncoloured.size())>i)
                {
                    colored.push_back(uncoloured[i]);
                    move_back(i,uncoloured);
                }
            }
        }
    }
    if(count<number_of_vertices-1)
    {
        return true;
    }
    return false;
}

void move_back_and_delete(int index, Lines& lines)
{
    //Precondition:
    assert(static_cast<int>(lines.size())>index);
    assert(static_cast<int>(lines.size())>=0);
    //Postcodition: Determines if there is a hamiltonpath in the graph.
    Line swap=lines[static_cast<int>(lines.size())-1];
    lines[static_cast<int>(lines.size())-1]=lines[index];
    lines[index]=swap;
    lines.pop_back();
}

void print_lines(Lines lines)
{
    for(int i=0; i<static_cast<int>(lines.size());i++)
    {
        cout<< lines[i][0]<<","<< lines[i][1]<<" - ";
    }
    cout<<endl;
}

void remove_neighboor(Node& node,int k)
{
    int swap=node.neighboors[static_cast<int>(node.neighboors.size())-1];
    node.neighboors[static_cast<int>(node.neighboors.size())-1]=node.neighboors[k];
    node.neighboors[k]=swap;
    node.neighboors.pop_back();
}

void Remove_lines_and_vertices(Node i,Node j,Nodes& nodes, Lines& lines)
{
//    cout<<"vertices: "<<i<<", "<<j<<" weghalen"<<endl;
    for(int index=0; index<static_cast<int>(lines.size())-1; index++)
    {
        if(lines[index][0]==i.vertice || lines[index][1]==i.vertice || 
		lines[index][0]==j.vertice || lines[index][1]==j.vertice)
        {
            move_back_and_delete(index, lines);
            index--;
        }
    }
    for(int index=0; index<static_cast<int>(nodes.size()); index++)
    {
        if(nodes[index].vertice==i.vertice || nodes[index].vertice==j.vertice)
        {
            move_back(index, nodes);
            index--;
        };
        for(int k=1;k<static_cast<int>(nodes[index].neighboors.size());k++)
            if(nodes[index].neighboors[k]==i.vertice || 
			nodes[index].neighboors[k]==j.vertice)
            {
                remove_neighboor(nodes[index],k);
            }
    }

}

bool connected1(Node i, Node j, Nodes nodes, Lines lines)
{
    Remove_lines_and_vertices(i,j,nodes,lines);
    if(connected(nodes,lines))
    {
        return true;
    }
    return false;
}

bool three_connected(Nodes nodes, Lines lines)
{
    for(int i=0; i<number_of_vertices-1; i++)
    {
        for(int j=i+1; j<number_of_vertices; j++)
        {
            if(!connected1(nodes[i],nodes[j],nodes,lines))
            {
                return false;
            }
        }
    }
    return true;
}


int main()
{
    cout<<endl<<endl<<endl;
    ///////////// GRAAF MAKEN //////////////////////////////
    Nodes nodes;
    for(int i=1;i<number_of_vertices+1;i++)
    {
        vector<int> neigboors;
        Node node;
        node.vertice=i;
        node.neighboors=neigboors;
        nodes.push_back(node);
    }
    Lines lines;
    Line line;
    for(int i=0;i<number_of_lines;i++)
    {
        line.push_back(LINES[i][0]);
        line.push_back(LINES[i][1]);
        lines.push_back(line);
        line.pop_back();
        line.pop_back();
        line.push_back(LINES[i][1]);
        line.push_back(LINES[i][0]);
        lines.push_back(line);
        line.pop_back();
        line.pop_back();
    }

////////// FUNCTIES //////////////////////////////////////

    Nodes already_visited;

    if(is_cubic(lines, nodes))
    {
        cout<<"Graph is cubic"<<endl;
        cout<<endl;
    }
    else
    {
        cout<<"Graph is not cubic"<<endl;
        cout<<endl;
    }

    Nodes blue;
    Nodes red;
    if(is_bipartite(lines,blue,red, nodes))
    {
        cout<<"Graph is bipartite"<<endl;
        cout<<endl;
    }
    else
    {
        cout<<"Graph is not bipartite"<<endl;
        cout<<endl;
    }


     if(three_connected(nodes,lines))
    {
        cout<<"Graph is 3-connected"<<endl;
        cout<<endl;
    }
    else
    {
        cout<<"Graph is not 3-connected"<<endl;
        cout<<endl;
    }

    if(contains_hamiltonpath(nodes, lines, already_visited))
    {
        cout<<"Hamiltonpath found:"<<endl;
        print(already_visited);
        cout<<endl;
    }
    else
    {
        cout<<"No Hamiltonpath found"<<endl;
        cout<<endl;
    }

    already_visited.clear();

    if(contains_hamiltoncycle(nodes[0], nodes, lines, already_visited, 
							nodes[0].vertice))
    {
        cout<<"Hamiltoncycle found:"<<endl;
        print(already_visited);
        cout<<endl;
    }
    else
    {
        cout<<"No Hamiltoncycle found"<<endl;
        cout<<endl;
    }
    return 0;
}

\end{lstlisting}

\end{document}